\documentclass[a4paper, 11pt]{article}

\usepackage{inputenc}
\usepackage{amsmath, amsthm, amssymb}
\usepackage{color}
\usepackage{graphicx}
\usepackage{hyperref}
\usepackage{a4wide}
\hypersetup{colorlinks=true,linkcolor=blue}

\usepackage{lipsum}
\usepackage{amsfonts}
\usepackage{graphicx}
\usepackage{epstopdf}
\usepackage{algorithmic}
\ifpdf
\DeclareGraphicsExtensions{.eps,.pdf,.png,.jpg}
\else
\DeclareGraphicsExtensions{.eps}
\fi

\usepackage{graphicx} 
\graphicspath{ {./images/} }
\usepackage[font={small,it}]{caption}
\usepackage{subcaption} 
\usepackage{cases}
\usepackage{mathtools}
\usepackage{xfrac}
\usepackage{rotating}
\usepackage{hyperref}
\usepackage{pdflscape}
\usepackage{inputenc}
\usepackage{amsmath, amssymb}
\usepackage{upgreek}
\usepackage{nicefrac}

\theoremstyle{theorem}
\newtheorem{theorem}{Theorem} [section]

\newtheorem{proposition}[theorem]{Proposition}
\theoremstyle{remark}
\newtheorem{remark}[theorem]{Remark}

\numberwithin{equation}{section}

\theoremstyle{definition}

\newenvironment{abst}{\begin{minipage}[c]{0.9\textwidth} \footnotesize \textbf{Abstract.}}
{\end{minipage}\\[2ex]}

\newenvironment{key}{\begin{minipage}[c]{0.9\textwidth} \footnotesize \textbf{Keywords.}}
{\end{minipage}\\[0ex]}

\newenvironment{amsclass}{\begin{minipage}[c]{0.9\textwidth} \footnotesize \textbf{AMS Subject Classification.}}
{\end{minipage}\\[2ex]}

\newcommand{\RNum}[1]{\uppercase\expandafter{\romannumeral #1\relax}}

\newcommand{\Sphere}{\mathbb{S}}
\newcommand{\Sphered}{\Sphere^{\scriptscriptstyle d-1}}
\newcommand{\R}{\mathbb{R}}

\newcommand{\HS}{\mathrm{HS}}

\newcommand{\polynd}{\widehat{\Pi}_n^d}

\DeclareFontFamily{U}{mathx}{}
\DeclareFontShape{U}{mathx}{m}{n}{<-> mathx10}{}
\DeclareSymbolFont{mathx}{U}{mathx}{m}{n}
\DeclareMathAccent{\widecheck}{0}{mathx}{"71}
\newcommand{\polyaltnd}{\widecheck{\Pi}_n^d}

\newcommand{\ba}{\begin{align*}}
	\newcommand{\ea}{\end{align*}}
\newcommand{\be}{\begin{equation}}
	\newcommand{\ee}{\end{equation}}

\newcommand{\LetterForDim}{{\mathcal{N}}}
\newcommand{\nd}{\LetterForDim_n^d}
\newcommand{\ndl}{\LetterForDim_\bw^d}

\newcommand{\nnd}{\widehat{\LetterForDim}_{n}^d}
\newcommand{\nmnd}{\widetilde{\LetterForDim}_{m,n}^d}
\newcommand{\nmd}{\widetilde{\LetterForDim}_{m,n}^d}
\newcommand{\nnnd}{\widecheck{\LetterForDim}_{n}^d}
\newcommand{\nds}{\LetterForDim_n(\Sphere^{d-1})}

\newcommand{\N}{\mathbb{N}}

\newcommand{\nollZ}{\mathsf{Z}}
\newcommand{\nollR}{\mathsf{R}}

\newcommand{\supscriptInCFourOl}{\upvarrho}
\newcommand{\BBCFour}[1]{\mathbb{B}^{d}_{#1}}

\newcommand{\ZernikePolynomials}{Z}

\newcommand{\degb}{\ell}

\newcommand{\bn}{\mathcal{B}_{L}}
\newcommand{\td}{\mathcal{S}_{D}}

\newcommand{\btbl}{\mathcal{BSB}_{D,L}}
\newcommand{\tbtl}{\mathcal{SBS}_{D,L}}
\newcommand{\limif}{\lim_{n\to\infty}}
\newcommand{\limifl}{\lim_{L\to\infty}}
\newcommand{\liminfif}{\liminf_{n\to\infty}}
\newcommand{\limsupif}{\limsup_{n\to\infty}}

\newcommand{\B}{\mathbb{B}}

\newcommand{\BB}{\mathbb{B}}

\newcommand\scalemath[2]{\scalebox{#1}{\mbox{\ensuremath{\displaystyle #2}}}}

\newcommand{\Kn}{\widehat{\mathcal{K}}_{n}}
\newcommand{\Klo}{\mathcal{K}_{L,\Omega}}
\newcommand{\Kmn}{\widetilde{\mathcal{K}}_{m,n}}

\newcommand{\nkkMone}{\mathrm{dim}(V_k^{d,\mu})}

\newcommand{\bw}{L}

\newcommand{\parzernike}{j}
\newcommand{\diffsymbol}{\mathrm{d}}

\begin{document}
\begin{center}
\Large\bfseries  Spatiospectral localization within the ball\\ -- studies on the influence of the spectral shape -- \normalsize\mdseries
\\[3ex]
{Christian Gerhards\footnotemark[1]}\footnotetext[1]{TU Bergakademie Freiberg, Institute of Geophysics and Geoinformatics, Germany, \texttt{christian.gerhards@geophysik.tu-freiberg.de}},
{Xinpeng Huang\footnotemark[2]}\footnotetext[2]{Central South University, School of Geoscience and Info-Physics, China,
	\texttt{xinpeng.huang@csu.edu.cn} (corresponding author)}
\\[3ex]
\today
\end{center}

\begin{abst}
We investigate the Slepian spatiospectral localization problem within subdomains of the $d$-dimensional ball. Opposed to the more classical setups of the Euclidean space or the sphere, the ball lacks a standard or universally accepted definition of bandwidth. Here, we consider a Fourier-Jacobi function system, decoupling the spherical and radial contributions via spherical harmonics and Jacobi polynomials. Special cases of this setup are of interest for various inverse problems in geophysics and medical imaging, since they relate to the underlying non-uniqueness, as well as in optics, where they represent the widely used Zernike polynomials. Bandwidth can be prescribed separately for the spherical and the radial contributions, where the particular choice of coupling between the two contributions determines the spectral shape, i.e., the overall notion of bandlimit. Understanding the effects of the spectral shape on the eigenvalue distribution of the Slepian spatiospectral localization problem can provide hints on particularly suitable notions of bandwidth for different applications. We provide rigorous asymptotic results for the spectral shape being defined via the overall polynomial degree as well as for being defined via sequential limits for the spherical and radial contributions. For various other spectral shapes, we provide numerical illustrations of the asymptotic eigenvalue distribution. Furthermore, we demonstrate a direct connection of the spectral shape to common indexing schemes for Zernike polynomials.
\end{abst}

\begin{key}
Slepian spatiospectral concentration, bandwidth, Zernike polynomials, Jacobi polynomials, spherical harmonics, reproducing kernel, eigenvalue distribution 
\end{key}

\begin{amsclass}
26B99, 42C10, 42C25, 47B34, 78M99, 86A99
\end{amsclass}

\section{Introduction}
Spatiospectral localization represents a fundamental concept in mathematical signal processing, concerned with the simultaneous concentration of signals in both spatial and spectral domains. The Landau-Pollak-Slepian-type spatiospectral concentration problem addresses a particularly important aspect: to what extent can a spectrally bandlimited signal be concentrated within a predefined spatial subdomain? This problem was initially investigated on the real line \cite{Slepian1961,Landau1961,Landau1962,Slepian1978} and has spawned extensive further analysis \cite{bonkar17,hoglac12,osirok13,xiarok01} and applications \cite{hoorey22,ledand20,peidin05,xuhay99}.
	
While these works focus on the Euclidean space $\R^d$, the Slepian concentration problem has more recently been extended to the sphere \cite{sneeuw99,Wieczorek2005,Simons2006b}, finding widespread applications in geophysics \cite{dashaj09,harsim12,handit07,Kim2017,plattner15}. Subsequent work has led to vectorial and tensorial extensions \cite{Plattner2014,Plattner2017,Michel2022}. These have been further incorporated into inverse problems, such as the downward continuation of potential field data \cite{Plattner2017,plamazger24,michelsimons18}.
	
Ball-shaped objects are ubiquitous in geophysics and astrophysics, from planetary bodies to small rock samples. For such objects, understanding the signal behavior requires consideration of both surface variations and variations in the internal radial structure, naturally leading to the question of signal localization within subdomains of the ball \cite{degroot22,lanras12,Michel_2008,simlor11}. However, in contrast to the well-developed theories for the Euclidean space and the sphere, a comprehensive treatment for the ball remains largely unexplored. Existing works \cite{Khalid2016a,leweke18} provide numerical evidence but lack (i) a unified framework for different bandwidth notions, (ii) rigorous proofs of asymptotic eigenvalue distributions, and (iii) an explicit characterization of the Shannon number's dependence on spectral shape.
	
This gap is compounded by the absence of a canonical approach to defining bandwidth on the ball. 
In Euclidean space $\mathbb{R}^d$, the Fourier transform provides a natural definition. On the sphere $\mathbb{S}^{d-1}$, spherical harmonic degree serves this purpose. For the ball, however, no such consensus exists. Different approaches employ Fourier-Laguerre and Fourier-Bessel functions \cite{Khalid2016a} or Fourier-Jacobi functions \cite{leweke18}. We adopt in this paper the latter framework, as Fourier-Jacobi functions naturally relate to non-uniqueness aspects in gravimetry, inverse magnetization problems, and EEG/MEG imaging \cite{Baratchart_2021,Leweke_2020,Michel_2008}, and encompass the widely used Zernike polynomials from optics \cite{Dai2007-ju,herled21,Noll:76,Tatian:74}.
	
A common feature of these function systems is their decomposition into spherical and radial components, permitting separate prescription of spherical and radial bandwidths. However, the effect of this bandwidth specification on spatiospectral localization properties has not been systematically investigated. \textit{The primary contributions of this work are:} (1) to establish a general framework for defining bandwidth on the ball and investigate its influence on the Shannon number, (2) to provide rigorous asymptotic analysis for bandwidth defined via polynomial degree and via sequential spherical-radial limits, (3) to demonstrate how different bandwidth notions---which we term \emph{spectral shapes}---optimize concentration for different regions of the ball, and (4) to connect these theoretical findings to practical applications by showing how specific spectral shapes correspond to established indexing schemes for Zernike polynomials \cite{Niu2022}. As an example, certain spectral shapes favor the concentration of near-surface structures (e.g., Earth's crust) over deeper structures (e.g., the mantle or core).
	
In the next section, we specify the precise setup and general type of results of this paper before proceeding to the more detailed derivations. In particular, the type of results that the reader can expect are condensed in Section \ref{sec:typeresults}.


\section{General framework and results} \label{sec:framework}
To start with some basic notation: by $\BBCFour{r}=\{x\in\R^d:\|x\|\leq r\}$ we mean the closed $d$-dimensional ball of radius $r>0$ and we use $\BB^d=\BBCFour{1}$ to abbreviate the closed unit ball. The $(d-1)$-dimensional unit sphere is denoted by $\Sphered=\{x\in\R^d:\|x\|=1\}$. Meanwhile, we use $\mathrm{vol}(\Sphered)=\tfrac{2\pi^{d/2}}{\Gamma({d}/{2})}$ to denote the volume of the ($d-1$)-dimensional unit sphere, where $\Gamma$ denotes
the Gamma function. By $L^2(\BB^d)$ we denote the space of square-integrable functions on the unit ball, with inner product $\langle f,g \rangle=\int_{\BB^d}f(x){g(x)}\diffsymbol x$ and norm $\|f\|=\langle f,f\rangle^{\frac{1}{2}}$. The dimension $d$ is supposed to satisfy $d\geq 2$. We use $\mathbb{N}$ for the set of positive integers and $\mathbb{N}_{0}$ for the non-negative integers.

\subsection{Orthogonal spherical Fourier-Jacobi basis}
The space of orthonormalized spherical harmonics $Y_{j,\degb}$ of degree $j$ on the sphere $\Sphered$ will be denoted by $H_j^d=\textnormal{span}\{Y_{j,\degb}\}_{\degb=1,\ldots,\dim(H_j^d)}$ (see, e.g., \cite{Mller1997AnalysisOS}). It holds  $\dim(H_j^d)=\binom{j+d-1}{j}-\binom{j+d-3}{j-2}$. By $P_i^{\alpha,\beta}$ we mean the Jacobi polynomial of degree $i$ with weight parameters $\alpha,\beta>-1$ .
Now, let $\{\supscriptInCFourOl_j\}_{j\in\N_0}$ be a sequence of real numbers such that $\inf_{j\in\N_{0}}\{{\supscriptInCFourOl}_j+\tfrac{d-2}{2}\}>-1$. Then we define, for $i,j\in \mathbb{N}_{0}$,  $\degb=1,\ldots,\mathrm{dim}(H_j^d)$, the Fourier-Jacobi functions
\begin{align}\label{def:ONS}
	{\ZernikePolynomials}_{i,j,\degb}(x)=\gamma_{ij}P_i^{0,{\supscriptInCFourOl}_j+\tfrac{d-2}{2}}(2r^2-1)\,r^{{\supscriptInCFourOl}_j}Y_{j,\degb}(\xi),\,
	& r=\|x\|\in[0,1], \ \xi=\frac{x}{\|x\|} \in \mathbb{S}^{d-1}, 
\end{align}
with $\gamma_{ij}$ being the normalization constant $\gamma_{ij}=\sqrt{4i+2{\supscriptInCFourOl}_j+d}$,
such that $\|\ZernikePolynomials_{i,j,\degb}\|=1$. The set $\{{\ZernikePolynomials}_{i,j,\degb}\}_{i,j\in \mathbb{N}_{0},\degb=1,\ldots,\mathrm{dim}(H_j^d)}$ forms a complete orthonormal system in $L^2(\BB^d)$. This has been shown, e.g., in \cite[Chap. 3]{MichelOrzlowsk2015} for the 3-dimensional ball and holds analogously for any $d$-dimensional setting, with $d\geq 2$. Note that for the particular choice ${\supscriptInCFourOl}_j=j$, the functions ${\ZernikePolynomials}_{i,j,\degb}$ are polynomials of degree $2i+j$. For $d=2,3$, they coincide with the Zernike polynomials that have been used in various optical imaging applications (e.g., \cite{Dai2007-ju,herled21,Noll:76,Tatian:74}). In dimension $d=3$, the functions ${\ZernikePolynomials}_{i,j,\degb}$ are suitable for modeling certain gravimetric and magnetic problems, with corresponding parameters ${\supscriptInCFourOl}_j=j$ and ${\supscriptInCFourOl}_j=j-1$, respectively (e.g., \cite{MichelOrzlowsk2015}). 
To stay notationally concise in this paper, we discuss only the polynomial setup, although some results also hold for general $\rho_j$. Therefore, if not mentioned otherwise, we fix ${\supscriptInCFourOl}_j=j$ for the remainder of the paper.

\subsection{Bandlimited Spaces associated with Fourier-Jacobi functions}
The basis function ${\ZernikePolynomials}_{i,j,\degb}$ essentially have the two parameters $i$ and $j$ that can be used to prescribe radial and spherical bandwidths, respectively (the parameter $\degb$ is typically assumed to cover the full range $\degb=1,\ldots,\mathrm{dim}(H_j^d)$). A very general notion of bandlimit can be achieved by relying on a predefined compact set $\Omega\subset[0,\infty)^2$, which may be considered a low-pass if it satisfies that $(x,y)\in\Omega$ implies $[0,x]\times[0,y] \subset \Omega$ (basically meaning that $\Omega$ covers the region bounded by the positive coordinate axes and a non-increasing function that intersects both axes
). Then, for some $L>0$, we define the linear space of functions
\begin{align}
	\Pi_{\bw,\Omega}^{d}=\mathrm{span}\left\{\ZernikePolynomials_{i,j,\degb}: \left(\frac{i}{\bw},\frac{j}{\bw}\right)\in\Omega, \degb =1,\ldots, \mathrm{dim}(H_{j}^{d})\right\}.\label{eqn:pilom}
\end{align}
We call $\bw$ and $\Omega$ the \emph{bandwidth} and \emph{spectral shape} of $\Pi_{\bw,\Omega}^{d}$, respectively. The spectral shape links the radial and spherical parameters $i$ and $j$ to define the bandlimited space $	\Pi_{\bw,\Omega}^{d}$. It also determines what we call \emph{notion of bandlimit}. The following particular choices of $\Omega$ provide some intuition on this link:
\begin{align}
\widetilde{\Pi}_{m,n}^d=\,&\text{span}\{{\ZernikePolynomials}_{i,j,\degb}: i\leq m, j\leq   n,\degb=1,\ldots,\mathrm{dim}(H_j^d)\}=\Pi_{\bw,\Omega}^{d},\label{eqn:nmspace}
	\\&\textnormal{for }\Omega=\{(x,y)\in[0,\infty)^2:x\leq 1,y\leq \kappa\}, \,L=m \textnormal{ and }\kappa=\tfrac{n}{m},\nonumber
	\\[1.5ex]	\polynd=&\text{span}\{{\ZernikePolynomials}_{i,j,\degb}: 2i+j\leq   n,\degb=1,\ldots,\mathrm{dim}(H_j^d)\}=\Pi_{\bw,\Omega}^{d},\label{eqn:hatpi}
	\\&\textnormal{for }\Omega=\{(x,y)\in[0,\infty)^2:x\leq \frac{1}{2},y\leq 1-2x\}, \,L=n,\nonumber
	\\[1.5ex]	\polyaltnd=&\text{span}\{{\ZernikePolynomials}_{i,j,\degb}: i+j\leq   n,\degb=1,\ldots,\mathrm{dim}(H_j^d)\}=\Pi_{\bw,\Omega}^{d},\label{eqn:ckeckpi}
	\\&\textnormal{for }\Omega=\{(x,y)\in[0,\infty)^2:x\leq 1,y\leq 1-x\}, \,L=n.\nonumber
\end{align}

\begin{figure}
	\centering
	\includegraphics[scale=0.25]{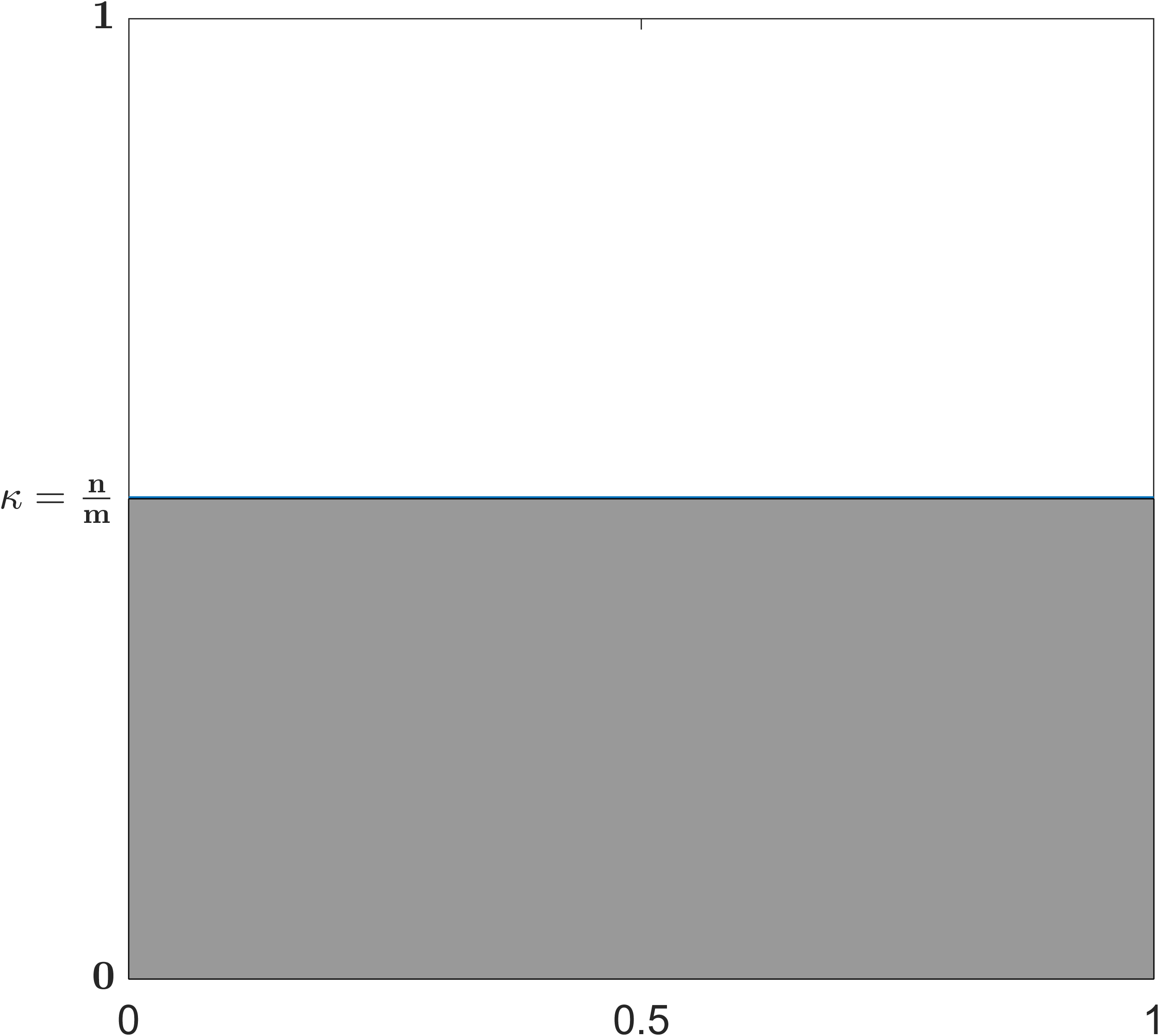}\qquad\includegraphics[scale=0.25]{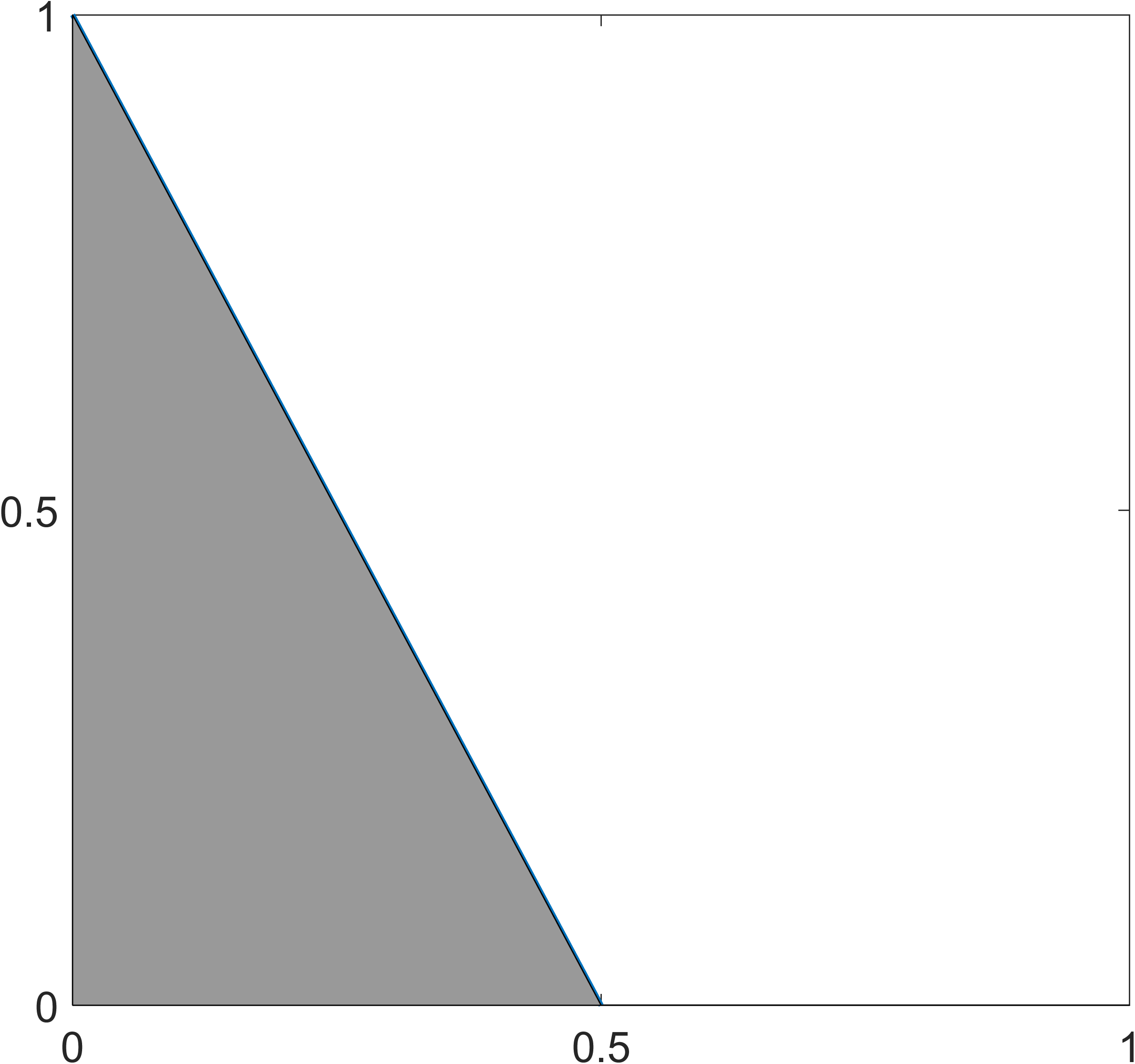}\qquad\includegraphics[scale=0.25]{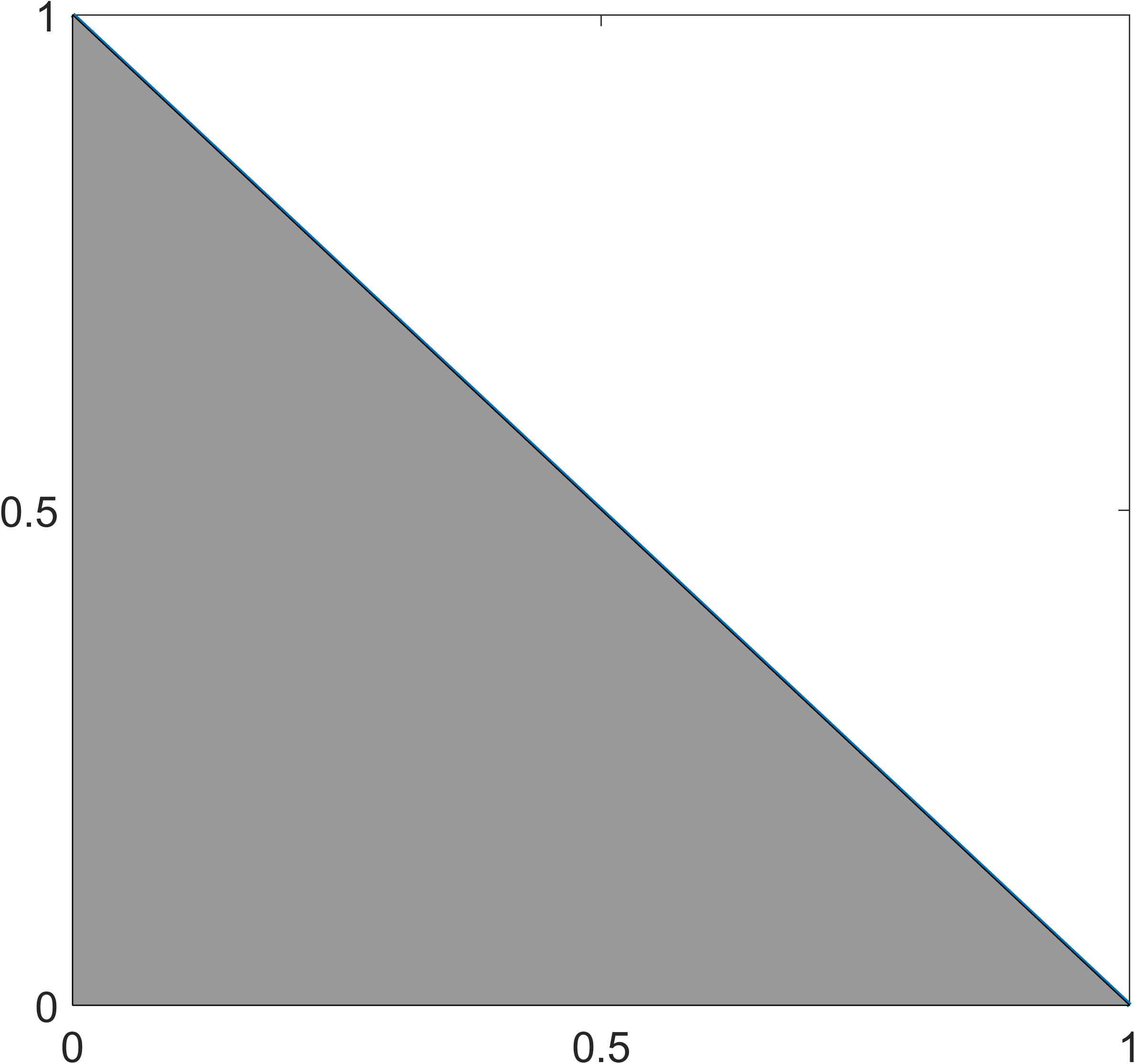}
	\caption{Illustration of spectral shapes $\Omega$ associated with notions of bandlimit~\eqref{eqn:nmspace},~\eqref{eqn:hatpi}, and~\eqref{eqn:ckeckpi} (from left to right).}\label{fig:omegaillust}
\end{figure}

Specifically, the definition in~\eqref{eqn:nmspace} involves two parameters $m$ and $n$. From the perspective of spectral shape, the family of spaces $\widetilde{\Pi}_{m,n}^d$ corresponds to a rectangular $\Omega$ whose extent is classified by the parameter $\kappa=\tfrac{n}{m}$ (Figure \ref{fig:omegaillust} illustrates the spectral shapes from above). Alternatively, if $\kappa$ is fixed, this linearly couples the parameters $m$ and $n$ and leads to the spaces
\begin{align}\label{eqn:tildekappapi}
	\widetilde{\Pi}_{n,\kappa n}^d=\,&\text{span}\{{\ZernikePolynomials}_{i,j,\degb}: i\leq n, j\leq \kappa n,\degb=1,\ldots,\mathrm{dim}(H_j^d)\}.
\end{align}
We rigorously analyze the concentration problem for the sequential limit case of $\widetilde{\Pi}_{m,n}^d$ (loosely speaking, the case $m\gg n$) in Section \ref{sec:seq} and numerically discuss the case $\widetilde{\Pi}_{n,\kappa n}^d$ and the role of $\kappa$ in Section \ref{sec:kappa}.

Moreover, for $\polynd$ we obtain
\begin{align}
	\polynd=\left\{\sum_{|\alpha|\leq n}c_{\alpha}\,x^{\alpha}:c_{\alpha}\in\mathbb{R}\right\},\label{eqn:polyeq}
\end{align}
where $\alpha$ is a $d$-dimensional non-negative multi-index, i.e., $\alpha=(\alpha_1,\alpha_2,...,\alpha_d)\in \N_{0}^{d},x^{\alpha}=x_{1}^{\alpha_1}x_{2}^{\alpha_2}\cdots x_{d}^{\alpha_d}$, and $|\alpha|=\sum_{i=1}^{d}\alpha_{i}$. In other words, $\polynd$ denotes the space of all polynomials with degree at most $n$ (equality \eqref{eqn:polyeq} is proven in appendix \ref{app:auxproofs} by a simple counting argument). The case of $\polynd$ is discussed in more detail in Section \ref{sec:poly}. 
The relation of $\polynd$ and $\polyaltnd$ to common indexing schemes for Zernike polynomials is discussed in Section \ref{sec:zernike}. Further examples of $\Pi_{\bw,\Omega}^{d}$ are illustrated in the appendix.

For further reference, we denote by
\begin{align}
	\ndl=\textnormal{dim}(\Pi_{\bw,\Omega}^{d})
\end{align}
the dimension of the bandlimited space $\Pi_{\bw,\Omega}^{d}$. For the special cases \eqref{eqn:nmspace}-\eqref{eqn:tildekappapi} we use the corresponding notations $\nmnd$, $\nnd$, $\nnnd$, respectively.

\subsection{Spatiospectral localization} Once we have defined the previous bandlimited function spaces, we can ask how well functions from these spaces can be concentrated in spatial subdomains $D\subset\BB^d$ of the unit ball. We seek functions $f$ in $\Pi_{\bw,\Omega}^{d}$ that maximize the Rayleigh quotient
\begin{align}\label{eqn:rayleighcoeff}
	\lambda(f)=\frac{\int_D|f(x)|^2\diffsymbol x}{\int_{\B^d}|f(x)|^2\diffsymbol x}.
\end{align}
The above closely relates to studying the eigenvalue distribution of the spatiospectral concentration operator $\td\bn\td$, with  $\td$ denoting the restriction to the spatial subdomain $D$ and $\bn$ the projection onto $\Pi_{\bw,\Omega}^{d}$ with bandwidth $L$. If $f$ from $\Pi_{\bw,\Omega}^{d}$ and $\lambda\geq0$ satisfy	
\begin{align}
	\td\bn\td f=\lambda f,\label{eqnC4:eigen1}
\end{align} 
then $\lambda$ matches the Rayleigh quotient for $f$ from \eqref{eqn:rayleighcoeff}. In other words, the largest eigenvalue of $\td\bn\td$ reflects the best possible spatial concentration within $D$ for a bandlimited function of bandwidth $L$. The second largest eigenvalue reflects the next best concentration achievable by a function orthogonal to the first one, and so on. We order the eigenvalues
\begin{align}
	\lambda_i=\lambda_i(D,L), \qquad 1\leq i\leq \ndl,
\end{align}
in a nonincreasing manner $1\geq \lambda_1\geq\lambda_2\geq\ldots\geq\lambda_{\ndl}\geq 0$. The notation $\lambda_i(D,L)$ is used to highlight the dependence on the choice of the spatial subdomain $D$ and the bandwidth $L$, although we typically just write $\lambda_i$ for brevity. A dependence on the spectral shape $\Omega$ is not indicated explicitly but should be clear from the context. In fact, the dependence on the spectral shape is what lies at the core of our discussion.

To be more precise, the spatial restriction operator takes the form
\begin{align}
	\td f=\textnormal{proj}_{L^2(D)}f=\chi_{D}\,f,
\end{align}
where $\chi_{D}$ denotes the characteristic function on $D$. The bandlimiting operator can be expressed as
\begin{align}
	\bn f=\textnormal{proj}_{\Pi_{\bw,\Omega}^{d}}f=\int_{\B^d}\mathcal{K}_{L,\Omega}(\cdot,x)f(x)\diffsymbol x,
\end{align}
where $\mathcal{K}_{L,\Omega}$ denotes the reproducing kernel of $\Pi_{\bw,\Omega}^{d}$. For brevity, we often write $\btbl=\bn\td\bn$ and $\tbtl=\td \bn\td$. 
Both operators are self-adjoint and non-negative definite, with operator norms bounded by one. In addition, since $\bn$ is of finite rank (and thus compact) and both $\bn$ and $\td$ are bounded, we have that the compositions $\btbl$ and $\tbtl$ are also finite-rank (and compact) and, therefore, allow an eigenvalue decomposition. From the fact $(\td\bn)^{*}=\bn\td$, one obtains that $\btbl$ and $\tbtl$ share the same non-zero eigenvalues. In consequence, for our purposes, it suffices to focus on $\tbtl$. It can be expressed as
\begin{align}
	\tbtl\, f=\int_{\B^d}\chi_D(\cdot)\,\mathcal{K}_{L,\Omega}(\cdot,x)\,\chi_D(x)f(x)\diffsymbol x.
\end{align}
For the particular cases \eqref{eqn:tildekappapi}-\eqref{eqn:nmspace} we use the corresponding notations  $\widetilde{\mathcal{SBS}}_{D,m,n}$, $\widehat{\mathcal{SBS}}_{D,n}$, $\widecheck{\mathcal{SBS}}_{D,n}$, $\widetilde{\mathcal{SBS}}_{D,n,\kappa n}$, respectively.

\subsection{Type of results}\label{sec:typeresults}

Studying the eigenvalue distribution of $\tbtl$ for different spectral shapes $\Omega$ is the main tool to obtain insights into spatiospectral localization properties with respect to the bandlimited space $\Pi_{\bw,\Omega}^{d}$. The two types of results that are of interest in the paper at hand are the following (with $\sharp$ denoting the cardinality of a finite set):

\paragraph{Eigenvalue property 1} For any $0<\varepsilon<1/2$, it holds
\begin{equation}\label{eigdis11}
	\limifl \frac{\sharp\{i:\varepsilon<\lambda_i(D;L)<1-\varepsilon\}}{\ndl}=0.
\end{equation}

\paragraph{Eigenvalue property 2} For a particular function $W=W_\Omega$ in $L^2(\B^d)$ and any $0<\tau<1$, it holds
\begin{equation}\label{eigdis21}
	\limifl \frac{\sharp\{i:\tau<\lambda_i(D;L)\leq 1\}}{\ndl}=\int_{D}W(x)\,\diffsymbol x.
\end{equation}

\begin{figure}
	\centering
	\includegraphics[scale=0.25]{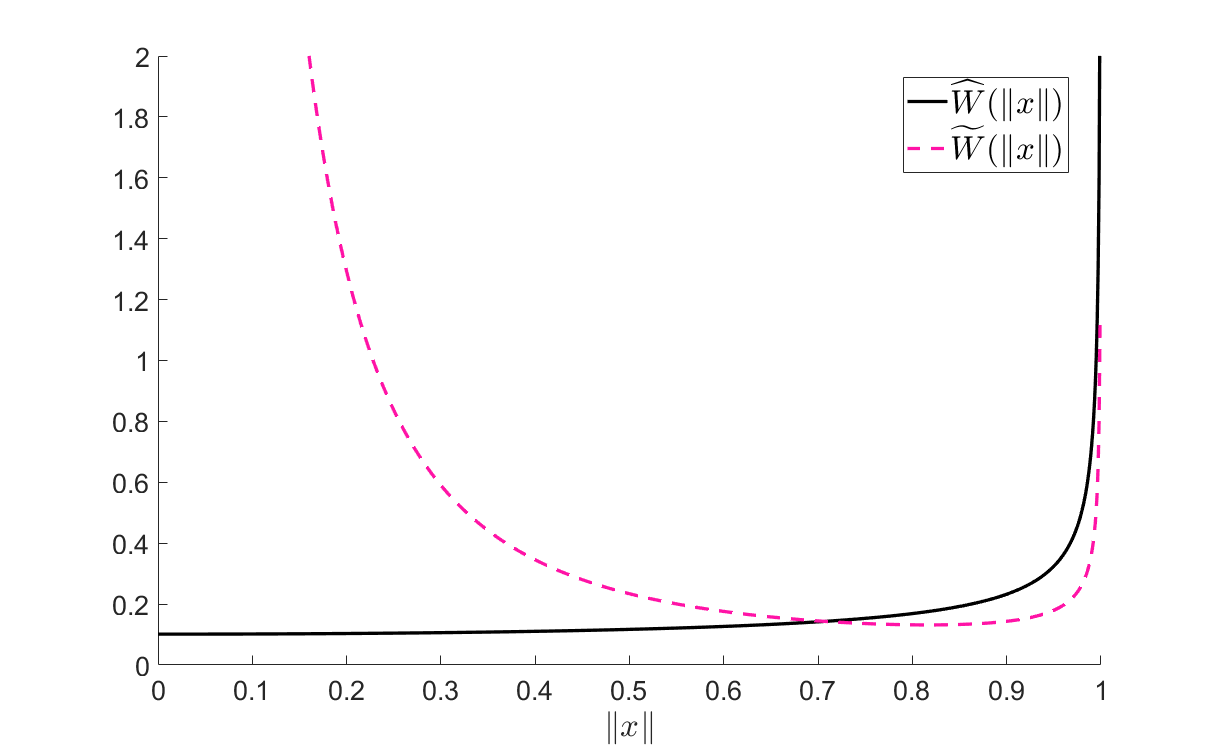}
	\caption{Illustration of the weight functions $W=\widehat{W}$ for $\widehat{\Pi}_n^3$ and $W=\widetilde{W}$ for $\widetilde{\Pi}_{m,n}^3$.}\label{fig:tildew1}
\end{figure}

The first one states that asymptocially all eigenvalues are clustered near one and zero, i.e., an eigenfunction is mostly either strongly spatially concentrated in the domain $D$ (eigenvalues close to one) or in the exterior $\B^d\setminus D$ (eigenvalues close to zero). This is a well-known phenomenon for the Euclidean case and for the case of the sphere. It is also expected for the ball and has been illustrated numerically in \cite{Khalid2016a,leweke18}, but not with respect to different notions of bandwidth and not with rigorous asymptotic results. The more interesting second property—which we call the relative Shannon number—provides estimates on the transition point between the relevant eigenvalues (those close to one) and the non-relevant eigenvalues (those close to zero).. Informally, it provides information on how many eigenfunctions in $\Pi_{\bw,\Omega}^{d}$ are efficient to reasonably approximate a signal in $D$:
\begin{equation}\label{eqn:shannongen}
	\sharp\{i:\tau<\lambda_i(D;L)\leq 1\}\simeq \ndl\int_{D}W(x)\,\diffsymbol x.
\end{equation}	
Again, this is well known for the real line (cf. \cite{slep83}):
\begin{align}
	{\sharp\{i:\tau<\lambda_i\leq 1\}}\simeq 2L T,\label{eqnC4:shannon11}
\end{align}
with $[-L,L]$ describing the interval of bandwidth and $D=[-T,T]$ the interval of spatial concentration, and for the 2-sphere (cf. \cite{Simons2006b}):
\begin{align}
	{\sharp\{i:\tau<\lambda_i\leq 1\}}\simeq (L+1)^2 \frac{A}{4\pi},\label{eqnC4:shannon21}
\end{align}
with $L$ denoting the maximal spherical harmonic degree that defines the bandwidth and $A$ the area of the subdomain $D$ of the 2-sphere. As before, for the case of the ball, \cite{Khalid2016a,leweke18} have provided formulae of the Shannon number for certain bandlimits, but they have neither provided asymptotic analysis nor have they investigated the auxiliary function $W$ in \eqref{eigdis21} and its dependence on the spectral shape $\Omega$.
Understanding the latter can not only be the mathematical foundation for the construction of appropriate systems of Slepian functions on the ball, but also link to computational and applied properties of these spaces (for example, the sampling density and taper method).
An interesting aspect of \eqref{eqnC4:shannon11} and \eqref{eqnC4:shannon21} is that the Shannon number possesses a multiplicative structure. One factor depends on the bandwidth, the other on the size of the subdomain D, but neither depends on the shape or location of D within the ambient space. The latter implies that the weight function W is constant in these two classical examples, reflecting the isotropy of their respective bandlimits. Since Fourier bandlimit and spherical harmonic bandlimit are translate and rotate invariant, respectively. However, this is not true anymore for the considered framework in the paper at hand, where $W$ may give different weights to different subregions of the ball. More interestingly, we expect $W=W_\Omega$ may differ for different spectral shapes $\Omega$, i.e., for different notions of bandwidth. Figure \ref{fig:tildew1} illustrates the two different auxiliary functions $W$ obtained here for the notions of bandwidth defined by the spaces $\polynd$ as in \eqref{eqn:hatpi} and the sequential-limit case for $\widetilde{\Pi}_{m,n}^d$ as in \eqref{eqn:tildekappapi}, respectively.

To be more precise, we expect that $W$ is related to the diagonal part of reproducing kernel $\mathcal{K}_{L,\Omega}$ of $\Pi_{\bw,\Omega}^{d}$ (also termed the Christoffel function in the study of polynomials) via 
\begin{align}
	\lim_{L\to\infty}\frac{\Klo(x,x)}{\ndl}=W(x),\label{eqn:reprolim}
\end{align}
for an adequate $W=W_\Omega$, which one can subsequently use to prove 
\begin{align}
	\limifl \frac{\mathrm{tr}(\tbtl)}{\ndl}&=\int_{D}W(x)\diffsymbol x,\label{asytr1}
	\\
	\limifl \frac{\|\tbtl\|_{\HS}^{2}}{\ndl}&=\int_{D}W(x)\diffsymbol x,\label{asyhs1}
\end{align}
with 
\begin{align}
	\mathrm{tr}(\tbtl)&=\sum_{i=1}^{\ndl}\lambda_i=\int_{D}\Klo(x,x)\diffsymbol x,\label{eqn:trace1}
	\\\|\tbtl\|_{\HS}^{2}&=\sum_{i=1}^{\ndl}\lambda_i^{2}=\int_{D}\int_{D}|\Klo(x,y)|^{2}\diffsymbol x\,\diffsymbol y\label{eqn:hsn1}
\end{align}
denoting the trace and the Hilbert-Schmidt norm of the operator $\tbtl$, respectively (see, e.g., \cite[Chap. VI.6]{reed1981functional}). Once \eqref{asytr1} and \eqref{asyhs1} are confirmed, this automatically implies the desired eigenvalue properties 1 and 2 from \eqref{eigdis11} and \eqref{eigdis21}. A proof of the latter is provided in the appendix \ref{app:eigdist}. Thus, it suffices to rigorously investigate \eqref{eqn:reprolim}-\eqref{asyhs1}. For the special cases $\polynd$ and $\widetilde{\Pi}_{m,n}^d$, we do this in Sections \ref{sec:poly} and \ref{sec:seq}. In particular, studying the asymptotic behaviour of ${\Klo(x,x)}/{\ndl}$ is useful to get an intuition about $W$. In classical setups, such limits have long been investigated, e.g., in \cite{Bos1994,Xu1996,Kroo2013a}, 
but there seem to be fairly few discussions that specifically focus on the Fourier-Jacobi framework.
For the special case $\widetilde{\Pi}_{n,\kappa n}^d$, we numerically illustrate the convergence of \eqref{eqn:reprolim} 
in Section \ref{sec:kappa} but cannot provide rigorous proofs. A related setup for Fourier-Bessel functions has been investigated systematically in \cite{huang23}. For notational convenience, we denote the reproducing kernels of the special cases \eqref{eqn:nmspace}-\eqref{eqn:tildekappapi} by $\widetilde{\mathcal{K}}_{m,n}$, $\widehat{\mathcal{K}}_{n}$, $\widecheck{\mathcal{K}}_{n}$, $\widetilde{\mathcal{K}}_{n,\kappa n}$, respectively. 

\paragraph{Relation to indexing schemes of Zernike polynomials} (2-d) Zernike polynomials are widely applied in optical imaging and various indexing schemes are proposed to arrange the countably many functions. The Noll indices (cf. \cite{Noll:76}) are just one example. We motivate that some commonly used indexing schemes (see, e.g., \cite{Niu2022} for an overview) basically fall into two categories that can be identified with the spaces $\polynd$ and $\polyaltnd$, respectively, indicating a connection to the notion of bandlimit. More details are provided in Section \ref{sec:zernike}.


\section{Spatiospectral concentration for the polynomial space $\polynd$}\label{sec:poly}

In this section, we consider the function space $\polynd$ of multivariate polynomials with maximal degree $n$. Polynomials are frequently used as building blocks in multivariate approximation (e.g., in \cite{Petrushev2008} for the construction of a needlet system on the ball), in this sense, the maximal polynomial degree can be considered a fairly natural notion of bandlimit.

Althogh we mainly focus on the $L^2(\BB^d)$ inner product (the unweighted/Lebesgue weighted setup), it would be helpful to start up with weighted polynomial spaces in the section, which relates to the behaviour of the radial Jacobi functions in our latter discussion. For any positive weight function $w:\BB^d\to\R$, the space $\polynd$ is a Hilbert space under the inner product $\langle f,g \rangle_{w}=\int_{\BB^d}f(x){g(x)}w(x)\diffsymbol x$. Here, we only focus on the special case of Jacobi weights $w_{\mu}$, with index $\mu\geq 0$, whose definition is given by (e.g., \cite{Xu1999})
\begin{align}\label{eqnC4:Wmu}
	w_{\mu}(x)=w_{\mu,d}(x)=\omega_{\mu}(1-\|x\|^2)^{\mu-\frac{1}{2}},
\end{align}
where $\omega_{\mu}$ is the normalization constant such that $\int_{\BB^d} w_{\mu}(x)\diffsymbol x=1$. 
In the following, we frequently use the notation $\diffsymbol \sigma_\mu(x)=w_{\mu}(x)\diffsymbol x$ for brevity. As a special case, the Jacobi weight with index $\mu=\frac{1}{2}$ leads to the Lebesgue measure, up to a multiplicative constant.


\subsection{Reproducing kernel and universality limit}\label{sec:polyunilim}

The reproducing kernel is of fundamental importance for the study of spatiospectral localization.
For polynomial spaces, the reproducing kernel associated with univariate polynomial admits closed-form expressions via the Christoffel-Darboux formula. The generalization to multi-variate polynomials is much more complicated (see, e.g., the discussion in \cite[Chap. 3.6]{Dunkl_Xu_2014}). Nevertheless, for $\polynd$ equipped with the inner product $\langle\cdot,\cdot\rangle_{w_\mu}$ on the ball $\BB^d$, a closed-form representation has been derived in \cite{Xu1999}. Their explicit expressions can be very useful but also very tedious to work with. Instead, the properties that are essential for our purposes are the upcoming two asymptotic properties of $\Kn^{\mu}$.

The first property is on the limit of $\Kn^{\mu}(x,x)$. It can be found, under varying conditions, in \cite{Bos1994,Xu1996,Kroo2013a}. Note the dimension of $\polynd$ is given by $\nnd=\binom{n+d}{d}=\frac{(n+d)!}{n!d!}=n^d/d!+\mathcal{O}(n^{d-1})$ (therefore, it holds $\lim_{n\to \infty}\nnd/n^d=\frac{1}{d!}$ for fixed $d$). 

\begin{proposition}\label{convergenceofChristoffelfunction}
	For any $x$ in the interior of $\BB^d$ (i.e., for $\|x\|<1$), the following pointwise limit holds true:
	\begin{equation}\label{chrf}
		\lim_{n\to\infty}\frac{\Kn^{\mu}(x,x)}{\nnd}=\frac{w_{0}(x)}{w_{\mu}(x)}.
	\end{equation}
	Furthermore, the above limit holds uniformly on any compact set $D$ contained in the interior of $\BB^d$.
\end{proposition}

The other required property, established in \cite[Thm. 1.7]{Kroo2013a}, is a multi-variate analogue of the universality limit \cite{Lubinsky2009}, which is a characterization of the asymptotic behaviour of $\Kn^{\mu}(x,y)$ near the diagonal line $x=y$.

\begin{proposition}\label{prop:universallimit}
	Let $D$ be a compact set contained in the interior of $\BB^d$. Then, uniformly for $x\in D$ and for $w,v$ in a compact subset of $\R^d$, it holds
	\begin{equation}\label{unilimit}
		\lim_{n\to\infty}\frac{\Kn^{\mu}(x+\frac{w}{n},x+\frac{v}{n})}{\Kn^{\mu}(x,x)}=\frac{J_{d/2}^{*}(\sqrt{G(x,w,v)})}{J_{d/2}^{*}(0)},
	\end{equation}
	where $J^{*}_{\alpha}(z)=z^{-\alpha}J_{\alpha}(z)$, with $J_{\alpha}$ denoting the Bessel function of the first kind, and $G$ being defined as
	\begin{equation}\label{gf}
		G(x,w,v)=\|w-v\|^{2}+\frac{|x\cdot(w-v)|^2}{1-\|x\|^2}.
	\end{equation}
\end{proposition}

\begin{remark}
	The invariance with respect to $\mu$ of the right-hand side in \eqref{unilimit} and of the numerator of the right-hand side in \eqref{chrf} is a useful property. In the proof in \cite{Kroo2013a}, the expression of \eqref{unilimit} is calculated by applying Mehler-Heine's asymptotic formula to the closed-form expression of $\Kn^{\mu}$ from \cite{Xu1999} for the case $\mu=0$. Proving \eqref{unilimit} for $\mu>0$ requires the so-called regularity condition 
	(sometimes also called Bernstein–Markov condition) to be satisfied by the measure $\sigma_\mu$, which facilitates the adaptation of Lubinsky's approach \cite{Lubinsky2009} for univariate polynomials to the multivariate setting. We also note that in the univariate case, a counterpart of Proposition \ref{prop:universallimit} exists (cf. \cite[Theorem 1.1]{Lubinsky2009}) and is of fundamental importance in the next section.
\end{remark}

\subsection{Asymptotics for the trace and Hilbert-Schmidt norm}

With the tools from the previous subsection, we can prove \eqref{asytr1} and \eqref{asyhs1} for the specific multivariate polynomial setup.

\begin{theorem}\label{thm:EssentialEstimation}Let $D\subset\B^d$ be a Lipschitz domain. Then, the following asymptotic relations hold true:
	\begin{align}\label{asytr}
		\limif \frac{\mathrm{tr}(\widehat{\mathcal{SBS}}_{D,n})}{\nnd}&=\int_{D}\widehat{W}(x)\diffsymbol x,
		\\\label{asyhs}
		\limif \frac{\|\widehat{\mathcal{SBS}}_{D,n}\|_{\HS}^{2}}{\nnd}&=\int_{D}\widehat{W}(x)\diffsymbol x,
	\end{align}
	with 
	\begin{align}
		\widehat{W}(x)=w_0(x)=\frac{\Gamma(\frac{d+1}{2})}{\pi^{\frac{d+1}{2}}}\frac{1}{\sqrt{1-\|x\|^2}}.
	\end{align}
\end{theorem}

\begin{proof}
	The proof of \eqref{asytr} follows directly from Proposition \ref{convergenceofChristoffelfunction}. But some care near the boundary of $\BB^d$ are necessary, since the convergence \eqref{chrf} only holds in the interior of the ball. Details are provided in the appendix \ref{app:auxproofs}.

The proof of \eqref{asyhs} follows from a two-sided estimate. The easier part is the upper bound. Due to the reproducing property, it holds
\begin{equation}\label{repro4}
	\Kn^{\mu}(x,x)=\int_{\BB^d}  \Kn^{\mu}(y,x) \Kn^{\mu}(x,y) \diffsymbol \sigma_\mu(y)=\int_{\BB^d} |\Kn^{\mu}(x,y)|^2\diffsymbol \sigma_\mu(y)\geq 0.
\end{equation}
Therefore,
	\begin{align}\frac{\|\widehat{\mathcal{SBS}}_{D,n}\|_{\HS}^{2}}{\nnd}& =\int_{D}\int_{D}\frac{|\Kn^{\mu}(x,y)|^{2}}{\nnd}\diffsymbol \sigma_\mu(y)\diffsymbol \sigma_\mu(x)\label{eqnC4:uuethm1} \\&\leq \int_{D}\int_{\BB^d}\frac{|\Kn^{\mu}(x,y)|^{2}}{\nnd}\diffsymbol \sigma_\mu(y)\diffsymbol \sigma_\mu(x) = \int_{D}\frac{\Kn^{\mu}(x,x)}{\nnd}\diffsymbol \sigma_\mu(x),\nonumber
	\end{align}
for which one can then argue as before in \eqref{asytr}. We now turn to the estimate from below. We will show that, for any $\epsilon_1, \epsilon_2>0$, it holds
	\begin{align}
		\liminf_{n\to\infty}\int_{D}\int_{D}\frac{|\Kn^{\mu}(x,y)|^{2}}{\nnd}\diffsymbol \sigma_\mu(x)\diffsymbol \sigma_\mu(y) \geq (1-\epsilon_1) \int_{D_{\epsilon_2}}w_0(x)\diffsymbol x,\label{eqnC4:llethm1}
	\end{align}
	with $D_{\epsilon_2}=\{x\in D: \|x-y\|\geq\epsilon_2 \text{ for all } y\in \partial D\}\subset \BBCFour{1-\epsilon_2}$. Assuming \eqref{eqnC4:llethm1} to be true and taking $\epsilon_1,\epsilon_2$ to zero, we get the desired lower bound, which finishes the proof of \eqref{asyhs}.
	
	Thus, it remains to prove \eqref{eqnC4:llethm1}. We shall make use of the universality limit to guarantee that the energy of $\Kn^{\mu}(x,y)$ concentrates close enough to the diagonal line $x=y$ as $n$ tends to infinity. Before proceeding to that, we need two auxiliary observations. First, we fix a point $x$ in the interior of $\BB^d$ and set $w=0$ in~\eqref{gf}. The set $G(x,L)=\{v\in \BB^d: \,G(x,0,v)\leq L^2\}$, for some $L>0$, is then bounded by an ellipsoid centered at $x$ with the shortest axis of length $L\sqrt{1-\|x\|^2}$ and all other axes of length $L$. This allows the integral transformation
	\begin{align}\label{eqnC4:integraloverG}
		\int_{G(x,L)}\left|\frac{J_{d/2}^{*}(\sqrt{G(x,0,v)})}{J_{d/2}^{*}(0)}\right|^2\diffsymbol v=\sqrt{1-\|x\|^2}\int_{\BBCFour{L}}\left|\frac{J_{d/2}^{*}(\|t\|)}{J_{d/2}^{*}(0)}\right|^2\diffsymbol t.
	\end{align}
	In the following, for $n \in \mathbb{N}$, $x$ in the interior of $\BB^d$, and $L>0$, we define the set
	\begin{align}
		G_n(x,L)=\{y\in \BB^d:G(x,0,n(y-x))\leq L^2\}.
	\end{align}
	It is a subset of $x+\BBCFour{\nicefrac{L}{n}}$ and, therefore, if additionally $x$ is in the interior of $D$ and $n$ sufficiently large, a subset of $D$. The second auxiliary observation (detailed calculations are moved to the appendix \ref{app:auxproofs}) concerns the quantity
	\begin{align}
		e_d=\int_{\R^d}\left|\frac{J_{d/2}^{*}(\|t\|)}{J_{d/2}^{*}(0)}\right|^2\diffsymbol t=\frac{\omega_0}{d!}.	\label{eqnC4:ednorm}
	\end{align}
	Since $e_d$ is finite, the following holds true: given any $\epsilon_1,\epsilon_2>0$, one can choose a constant $L>0 $ large enough such that
	\begin{align}
		\int_{\BBCFour{L}}\left|\frac{J_{d/2}^{*}(\|t\|)}{J_{d/2}^{*}(0)}\right|^2\diffsymbol t>e_d-\epsilon_1.\label{eqnC4:jgeqeps}
	\end{align}
	And given such an $L$, it holds that $G_n(x,L)\subset x+\BBCFour{\nicefrac{L}{n}}\subset D$ for all $n>\nicefrac{L}{\epsilon_2}$ and $x\in D_{\epsilon_2}$.

	With these auxiliary observations, we can proceed to showing that \eqref{eqnC4:llethm1} holds true. For an $L>0$ such that \eqref{eqnC4:jgeqeps} is satisfied and for $n>\nicefrac{L}{\epsilon_2}$, we have
	\begin{align}\label{esths}
		&\int_{D}\int_{D}\frac{|\Kn^{\mu}(x,y)|^{2}}{\nnd} \diffsymbol \sigma_\mu(x)\diffsymbol \sigma_\mu(y)\nonumber
		\\&\geq \int_{D_{\epsilon_2}}\int_{G_n(x,L)}\frac{ |\Kn^{\mu}(x,y)|^{2}}{\nnd}\diffsymbol \sigma_\mu(y)\diffsymbol \sigma_\mu(x) \nonumber
		\\&=\int_{D_{\epsilon_2}}  \int_{G_n(x,L)} \left|\frac{\Kn^{\mu}(x,y)}{\Kn^{\mu}(x,x)}\right|^{2}\frac{\Kn^{\mu}(x,x)^2}{\nnd}\diffsymbol \sigma_\mu(y)\diffsymbol \sigma_\mu(x) 
		\\&=\int_{D_{\epsilon_2}} \int_{G(x,L)} \frac{\Kn^{\mu}(x,x)}{n^d} \left|\frac{\Kn^{\mu}(x,x+\frac{v}{n})}{\Kn^{\mu}(x,x)}\right|^{2}\frac{\Kn^{\mu}(x,x)}{\nnd}w_\mu\left(x+\frac{v}{n}\right)\diffsymbol v\,\diffsymbol \sigma_\mu(x),\nonumber
	\end{align}
	where the substitution $v=n(y-x)$ has been used for the last equality. Since $\lim_{n\to\infty}{\nnd}/{n^d}=\frac{1}{d!}$, Proposition \ref{convergenceofChristoffelfunction} implies uniform convergence of ${\Kn^{\mu}(x,x)}/{n^d}$ and ${\Kn^{\mu}(x,x)}/{\nnd}$ with respect to $x\in D_{\epsilon_2}$ as $n$ tends to infinity. Additionally, Proposition \ref{prop:universallimit} implies uniform convergence of ${\Kn^{\mu}(x,x+\frac{v}{n})}/{\Kn^{\mu}(x,x)}$ with respect to $x\in D_{\epsilon_2}$ and $v\in G(x,L)$ as $n$ tends to infinity. Thus, we may interchange the order of limit and integration in the last line of \eqref{esths} and get
	\begin{align}\label{eqnC4:thm2finest}
		&\liminf_{n\to\infty}\int_{D}\int_{D}\frac{|\Kn^{\mu}(x,y)|^{2}}{\nnd} \diffsymbol \sigma_\mu(x)\diffsymbol \sigma_\mu(y)\nonumber
		\\&\geq \int_{D_{\epsilon_2}}\int_{G(x,L)} \lim_{n\to\infty}  \frac{\Kn^{\mu}(x,x)}{n^d} \left|\frac{\Kn^{\mu}(x,x+\frac{v}{n})}{\Kn^{\mu}(x,x)}\right|^{2}\frac{\Kn^{\mu}(x,x)}{\nnd}w_\mu\left(x+\frac{v}{n}\right)\diffsymbol v\,\diffsymbol \sigma_\mu(x)\nonumber
		\\&= \int_{D_{\epsilon_2}}\int_{G(x,L)} \frac{w_0(x)}{d!w_\mu(x)}  \left|\frac{J_{d/2}^*(\sqrt{G(x,0,v)})}{J_{d/2}^*(0)}\right|^2  \frac{w_0(x)}{w_\mu(x)}w_\mu(x)\diffsymbol v\,\diffsymbol \sigma_\mu(x)
		\\&=\int_{D_{\epsilon_2}} \frac{w_0(x)}{d!}\sqrt{1-\|x\|^2} \int_{\BBCFour{L}} \left|\frac{J_{d/2}^*(\|v\|)}{J_{d/2}^*(0)}\right|^2 \diffsymbol v\,w_0(x)\diffsymbol x\nonumber
		\\&\geq\frac{\omega_0}{d!}\left(e_d-\epsilon_1\right)\int_{D_{\epsilon_2}} w_0(x)\diffsymbol x.\nonumber
	\end{align}
	The second to last line follows from \eqref{eqnC4:integraloverG} and the last line from \eqref{eqnC4:jgeqeps}. Finally, observing \eqref{eqnC4:ednorm} and possibly modifying $\epsilon_1$ by a factor $\frac{d!}{\omega_0}$, \eqref{eqnC4:thm2finest} provides the desired proof for \eqref{eqnC4:llethm1}.
\end{proof}

\subsection{Asymptotic eigenvalue distribution}

With Theorem \ref{thm:EssentialEstimation} established, the desired eigenvalue properties for the spatiospectral localization operator $\widehat{\mathcal{SBS}}_{D,n}$ (reflecting the notion of bandlimit defined via the maximal polynomial degree) follow directly from the general statement of Theorem \ref{thm:geneigdist}.

\begin{theorem}\label{thm:1}
	Let $D\subset\B^d$ be a Lipschitz domain. Then, for any $0<\varepsilon<1/2$ and any $0<\tau<1$, it holds
	\begin{align}\label{eigdis1}
		\limif \frac{\sharp\{i:\varepsilon<\widehat{\lambda}_i(D;n)<1-\varepsilon\}}{\nnd}&=0,
		\\\label{eigdis2}
		\limif \frac{\sharp\{i:\tau<\widehat{\lambda}_i(D;n)\leq 1\}}{\nnd}&=\int_{D}\widehat{W}(x)\diffsymbol x,
	\end{align}
	with  
	\begin{align}
		\widehat{W}(x)=\frac{\Gamma(\frac{d+1}{2})}{\pi^{\frac{d+1}{2}}}\frac{1}{\sqrt{1-\|x\|^2}}.
	\end{align}
\end{theorem}


\section{Spatiospectral concentration for the space $\widetilde{\Pi}_{m,n}^d$: the sequential-limit case}\label{sec:seq}
This section we consider the concentration problem for spaces $\widetilde{\Pi}_{m,n}^d$ from \eqref{eqn:nmspace}, this type of bandlimit tries to constrain the radial and spherical freedom separately, which naturally originates from applications (e.g., in \cite{herled21,Schneider25}). However, the formal analysis presented in this section is limited to the case of $\widetilde{\Pi}_{m,n}^d$ under a sequential limit condition, where the radial bandlimit $m\to\infty$ is taken before the spherical bandlimit $n\to\infty$. While sequential limits may be of limited applicability, this approach serves two important purposes in our investigation of spectral shapes. First, it provides a foundational benchmark case where the asymptotic effects of radial and spherical bandlimiting can be cleanly separated, yielding an explicit expression for the weight function $\widetilde{W}(x)$. Second, as demonstrated numerically in Section \ref{sec:kappa}, this sequential limit represents the extreme case of linear coupling $n=\kappa m$ where radial resolution dominates ($\kappa\to 0$), thereby establishing one boundary of the spectrum of possible concentration behaviors governed by spectral shape.

The dimension of $\widetilde{\Pi}_{m,n}^d$ is given by
\begin{align}\label{eqnC4:nmd} 
	\nmd=\dim(\widetilde{\Pi}_{m,n}^d)=(m+1)\nds,
\end{align}
with $\nds$ denotes the dimension of spherical harmonics up to degree $n$ 
\begin{align}\label{eqnC4:nsd}
	\nds=\sum_{j=1}^n\dim(H_j^d)=\binom{n+d-1}{n}+\binom{n+d-2}{n-1}.
\end{align}
Thus, we get $\nds=\frac{2n^{d-1}}{(d-1)!}+\mathcal{O}(n^{d-2})$ and $\nmd=\frac{(m+1)2n^{d-1}}{(d-1)!}+\mathcal{O}(mn^{d-2})$.

\subsection{Analysis of the reproducing kernel}

Substituting the explicit expression of ${\ZernikePolynomials}_{i,j,\degb}$ from \eqref{def:ONS}, we get for the reproducing kernel of  $\widetilde{\Pi}_{m,n}^d$ that
\begin{align}\label{eqnC4:reproducingkerneloftildePi}
&{\Kmn}(x,y)=\sum_{i\leq m}\sum_{j\leq n}\sum_{\degb \leq \mathrm{dim}(H_j^d)} {\ZernikePolynomials}_{i,j,\degb}(x){{\ZernikePolynomials}_{i,j,\degb}(y)}\\
\nonumber	&=\scalemath{0.90}{ \sum_{j\leq n}\left(\sum_{i\leq m} \gamma_{ij}^2 \,(r_xr_y)^{\parzernike} P_i^{0,\parzernike+\tfrac{d-2}{2}}(2r_x^2-1) P_i^{0,\parzernike+\tfrac{d-2}{2}}(2r_y^2-1)\right) \left(\sum_{\degb \leq \mathrm{dim}(H_j^d)}Y_{j,\degb}(\xi_x){Y_{j,\degb}(\xi_y)}\right),}
\end{align}
where $r_x=\|x\|, r_y=\|y\|$ and $\xi_x=\frac{x}{\|x\|},\xi_y=\frac{y}{\|y\|}$ denote the radial and spherical components of $x$ and $y$, respectively. One notices that the radial and spherical contributions decouple with respect to sums over $i$ and $\degb$, but they are coupled by summation over the degree $j$. Asymptotic estimations as in Section \ref{sec:polyunilim} of this kernel are not readily available here. We first have to put some effort into deriving related properties for the current setup. The sum over $\degb$ can be treated fairly easily via the addition theorem for spherical harmonics (see, e.g., \cite{Mller1997AnalysisOS}), namely, 
\begin{align}\label{eqnC4:reproducingkernelofSH}
\sum_{\degb \leq \mathrm{dim}(H_j^d)}Y_{j,\degb}(\xi_x){Y_{j,\degb}(\xi_y)} =\frac{\mathrm{dim}(H_j^d)}{\mathrm{vol}(\Sphered)} P_{j}^{(d)}(\xi_x\cdot\xi_y),
\end{align}
where $P_{j}^{(d)}$ denotes the Legendre polynomial of degree $j$ for the $d$-dimensional setup, normalized to satisfy $P_{j}^{(d)}( 1)=1$. Thus, the main effort lies in the analysis of the radial contribution. 

The following observation can be seen as an analogue of Proposition \ref{convergenceofChristoffelfunction}. Its proof is provided in the appendix \ref{app:univlimit}.

\begin{proposition}\label{prop:convergenceofChristoffelfunctionoftildePi}
For any fixed $n\in\N_0$ and for any $x$ in the interior of $\BB^d\setminus\{0\}$ (i.e., for $0<\|x\|<1$), the following pointwise limit holds true:
\begin{align}
	\lim_{m\to\infty} \frac{\Kmn(x,x)}{\nmd} = \widetilde{W}(x),
\end{align}
with $\widetilde{W}$ given by
\begin{align}
	\widetilde{W}(x)=\frac{2}{\pi \mathrm{vol}(\Sphered)}\frac{1}{\|x\|^{d-1}\sqrt{1-\|x\|^2}}.\label{eqn:wtilde1}
\end{align}
Furthermore, the above limit holds uniformly for any compact subset $D$ in the interior of $\BB^d\setminus\{0\}$. 
\end{proposition}

\begin{remark} 
The previous proposition implies $\lim_{n\to\infty}\lim_{m\to\infty} {{\Kmn}(x,x)}/{\nmd}$ $=\widetilde{W}(x)$, but only in the sense of taking the sequential limit (this is because, in the upcoming proof, we rely on a result for one-dimensional universality limits \cite{Lubinsky2009} that does not provide further information on the behaviour of the limit with respect to the ambient parameter $n$). Actually, the limit with coupled bandwidths $m$ and $n$ can be very different, as we illustrate numerically in Section \ref{sec:kappa}.
\end{remark}

We are not able to provide a universality limit fully analogous to that of Proposition \ref{prop:universallimit}, but we can show the following closely related although  weaker result that suffices for the proof of the desired estimate on the Hilbert-Schmidt norm of $\widetilde{\mathcal{SBS}}_{D,m,n}$. Its proof is provided in the appendix \ref{app:univlimit}.

\begin{proposition}\label{prop:unilimit2}
Let $D$ be a compact set contained in the interior of $\BB^d\setminus\{0\}$. Furthermore, for some given $x\in D$, we define $x_{t,\xi}^m=(\|x\|+\frac{t}{m+1})\xi$, with $\xi\in\Sphered$ and $t\in\R$. Then, uniformly for $x\in D$, for $\xi\in\Sphered$, and for $t$ in a compact subset of $\R$, it holds
\begin{equation}\label{unilimit2}
	\lim_{m\to\infty}\frac{{\Kmn}(x,x_{t,\xi}^m)}{{\Kmn}(x,x)}=\frac{\mathrm{vol}(\Sphered)}{\nds}\,\textnormal{sinc}\left(\frac{2t}{\sqrt{1-\|x\|^2}}\right)\mathcal{K}_{\mathrm{Harm}_{n}}\left(\frac{x}{\|x\|},\xi\right),
\end{equation}
where $\textnormal{sinc}(x)=\sin(x)/x$ and
\begin{equation}\label{gf2}
	\mathcal{K}_{\mathrm{Harm}_{n}}\left(\eta,\xi\right)=\sum_{j\leq n}\frac{\mathrm{dim}(H_j^d)}{\mathrm{vol}(\Sphered)}P_{j}^{(d)}( \eta\cdot \xi),
\end{equation}	
which is the reproducing kernel of $\mathrm{Harm}_{n}(\Sphered)=\oplus_{0\leq j \leq n}H_{j}^{d}$.
\end{proposition}

\subsection{Asymptotics for the trace and Hilbert-Schmidt norm}

With the tools from the previous subsection, we can prove \eqref{asytr1} and \eqref{asyhs1} for the setup of $\widetilde{\Pi}_{m,n}^d$.

\begin{theorem}\label{thm:EssentialEstimation2}Let $D\subset\B^d$ be a Lipschitz domain. Then, the following asymptotic relations hold true:
\begin{align}\label{asytr2}
	\limif \lim_{m\to\infty}\frac{\mathrm{tr}(\widetilde{\mathcal{SBS}}_{D,m,n})}{\nmd}&=\int_{D}\widetilde{W}(x)\diffsymbol x,
	\\\label{asyhs2}
	\limif \lim_{m\to\infty}\frac{\|\widetilde{\mathcal{SBS}}_{D,m,n}\|_{\HS}^{2}}{\nmd}&=\int_{D}\widetilde{W}(x)\diffsymbol x,
\end{align}
with 
\begin{align}
	\widetilde{W}(x)=\frac{2}{\pi \mathrm{vol}(\Sphered)}\frac{1}{\|x\|^{d-1}\sqrt{1-\|x\|^2}}.
\end{align}
It needs to be emphasized that the limits in \eqref{asytr2} and \eqref{asyhs2} should be interpreted as sequential limits, i.e., first taking the limit over $m$  and then taking the limit over $n$.
\end{theorem}

\begin{proof}
Having Proposition \ref{prop:convergenceofChristoffelfunctionoftildePi} at hand, the proof of \eqref{asytr2} is identical to the corresponding proof in Theorem \ref{thm:EssentialEstimation}. The same holds true for the upper bound
\begin{equation}
	\limif \lim_{m\to\infty}\frac{\|\widetilde{\mathcal{S}\mathcal{B}\mathcal{S}}_{D,m,n}\|_{\HS}^{2}}{\nmd}\leq\int_{D}\widetilde{W}(x)\diffsymbol x.\label{eqnC4:ubths}
\end{equation}
However, the corresponding lower bound
\begin{equation}
	\limif \lim_{m\to\infty}\frac{\|\widetilde{\mathcal{S}\mathcal{B}\mathcal{S}}_{D,m,n}\|_{\HS}^{2}}{\nmd}\geq\int_{D}\widetilde{W}(x)\diffsymbol x\label{eqnC4:lbths}
\end{equation}
requires some more effort. This will be elaborated in the remainder of the proof.

We start with some auxiliary notation. The spherical cap with center $\xi\in\Sphered$ and (polar) radius $\epsilon_2>0$ is denoted by $\mathcal{C}_{\epsilon_2}(\xi)=\{\eta\in\Sphered:1-\xi\cdot\eta<\epsilon_2\}$. Furthermore, by $\mathcal{U}_m(x,L,\epsilon_2)=\{y\in \BB^d: r_x-\tfrac{L}{m+1}<r_y<r_x+\tfrac{L}{m+1}, \xi_y\in\mathcal{C}_{\epsilon_2}(\xi_x)\}$ we denote a truncated spherical cone with additional parameters $x\in\BB^d$ and $L>0$. We assume $\epsilon_2>0$ to be arbitrary but fixed. For any $x\in D$ with $\mathcal{U}_m(x,L,\epsilon_2)\subset D\setminus\{0\}$, we can then estimate
\begin{align}\label{eqnC4:concentrationoftildeK}
	& \int_{D}\frac{|{\Kmn}(x,y)|^2 }{\nmd}\diffsymbol y  \geq  \int_{\mathcal{U}_m(x,L,\epsilon_2)}\frac{|{\Kmn}(x,y)|^2}{|{\Kmn}(x,x)|^2}\nmd\,  \,\frac{|{\Kmn}(x,x)|^2}{|\nmd|^2}\diffsymbol y 
	\\&= \int_{-L}^L\int_{\mathcal{C}_{\epsilon_2}(\xi_x)}\frac{|{\Kmn}(x,x_{t,\xi}^m)|^2}{|{\Kmn}(x,x)|^2}\frac{1}{m+1}\left(\|x\|+\frac{t}{m+1}\right)^{d-1}\nmd \,\frac{|{\Kmn}(x,x)|^2}{|\nmd|^2}\, \diffsymbol \omega(\xi)\, \diffsymbol t\nonumber
	\\&= \int_{-L}^L\int_{\mathcal{C}_{\epsilon_2}(\xi_x)}\frac{|{\Kmn}(x,x_{t,\xi}^m)|^2}{|{\Kmn}(x,x)|^2}\left(\|x\|+\frac{t}{m+1}\right)^{d-1}\nds \,\frac{|{\Kmn}(x,x)|^2}{|\nmd|^2}\, \diffsymbol \omega(\xi)\, \diffsymbol t,\nonumber
\end{align}
where $\diffsymbol \omega$ denotes the surface measure on the unit sphere $\Sphered$. For the equality in the second line, we have used the co-area formula and  $\mathcal{U}_m(x,\epsilon_2)=\{x_{t,\xi}^m: t\in[-L,L], \xi\in \mathcal{C}_{\epsilon_2}(\xi_x)\}$, with $x_{t,\xi}^m=(\|x\|+\frac{t}{m+1})\xi$ as in Proposition \ref{prop:unilimit2}.
When taking the limit $m$ to infinity, we are allowed to interchange the order of integral and limit in the last line of \eqref{eqnC4:concentrationoftildeK} due to the uniform  convergence of the integrand for $\xi\in\Sphered$ and for $t$ in a compact subset of $\R$, as indicated in Proposition \ref{prop:unilimit2}. Additionally remembering Proposition \ref{prop:convergenceofChristoffelfunctionoftildePi}, we can continue from there to further estimate 
\begin{align} \label{eqnC4:esths2}
	& \lim_{m\to\infty}\int_{D}\frac{|{\Kmn}(x,y)|^2 }{\nmd}\diffsymbol y  
	\\&\geq  \int_{-L}^L\int_{\mathcal{C}_{\epsilon_2}(\xi_x)} \left|\textnormal{sinc} \left(\frac{2t}{\sqrt{1-\|x\|^2}}\right)\frac{\mathrm{vol}(\Sphered)}{\nds}\mathcal{K}_{\mathrm{Harm}_{n}}(\xi_x,\xi)\right|^2 \|x\|^{d-1}\nds\,|\widetilde{W}(x)|^2\,\diffsymbol \omega(\xi)\, \diffsymbol t\nonumber
	\\&=  \int_{-L}^L \left|\textnormal{sinc} \left(\frac{2t}{\sqrt{1-\|x\|^2}}\right)\right|^2\mathrm{vol}(\Sphered)\|x\|^{d-1}|\widetilde{W}(x)|^2\, \diffsymbol t\,\int_{\mathcal{C}_{\epsilon_2}(\xi_x)}\frac{\mathrm{vol}(\Sphered)}{\nds}\left|\mathcal{K}_{\mathrm{Harm}_{n}}(\xi_x,\xi)\right|^2\diffsymbol \omega(\xi)\nonumber \\
	&=  \int_{-\frac{2L}{\sqrt{1-\|x\|^2}}}^{\frac{2L}{\sqrt{1-\|x\|^2}}} \left|\textnormal{sinc} \left(s\right)\right|^2\frac{\sqrt{1-\|x\|^2}}{2}\mathrm{vol}(\Sphered)\|x\|^{d-1}|\widetilde{W}(x)|^2\,\diffsymbol s\, \int_{\mathcal{C}_{\epsilon_2}(\xi_x)}\frac{\mathrm{vol}(\Sphered)}{\nds}\left|\mathcal{K}_{\mathrm{Harm}_{n}}(\xi_x,\xi)\right|^2\diffsymbol \omega(\xi)\nonumber
	\\&=  \int_{-\frac{2L}{\sqrt{1-\|x\|^2}}}^{\frac{2L}{\sqrt{1-\|x\|^2}}} \left|\textnormal{sinc} \left(s\right)\right|^2\diffsymbol s\,\frac{\widetilde{W}(x)}{\pi}\, \int_{\mathcal{C}_{\epsilon_2}(\xi_x)}\frac{\mathrm{vol}(\Sphered)}{\nds}\left|\mathcal{K}_{\mathrm{Harm}_{n}}(\xi_x,\xi)\right|^2\diffsymbol \omega(\xi),\nonumber
\end{align}
where we have substituted $s=2t/\sqrt{1-\|x\|^2}$ in the fourth line. Next, we choose an arbitrary but fixed $\epsilon_1>0$. It is well-known that $\int_{\R}|\textnormal{sinc}(s)|^2ds=\pi$. So, we can choose a sufficiently large $L>0$ such that 
\begin{align*}
	\int_{-L}^L|\textnormal{sinc}(s)|^2 \diffsymbol s\geq(1-\epsilon_1)\pi.
\end{align*}
Since it holds $2L/\sqrt{1-\|x\|^2}\geq L$, we can apply the above to the last line in \eqref{eqnC4:esths2} and get
\begin{align}\label{eqnC4:esths3}
	& \lim_{m\to\infty}\int_{D}\frac{|{\Kmn}(x,y)|^2 }{\nmd}\diffsymbol y \geq  (1-\epsilon_1)\,\widetilde{W}(x) \int_{\mathcal{C}_{\epsilon_2}(\xi_x)}\frac{\mathrm{vol}(\Sphered)}{\nds}\left|\mathcal{K}_{\mathrm{Harm}_{n}}(\xi_x,\xi)\right|^2\diffsymbol \omega(\xi).
\end{align}
It remains to investigate the integral over the spherical cap $\mathcal{C}_{\epsilon_2}(\xi_x)$ as we take the limit $n$ to infinity. In fact, it holds uniformly for $\xi_x\in\Sphere$ that
\begin{align}\label{eqnC4:concentrationonsphere}
	\lim_{n\to\infty}\int_{\mathcal{C}_{\epsilon_2}(\xi_x)}|\mathcal{K}_{\mathrm{Harm}_{n}}(\xi_x,\xi)|^2 \frac{\mathrm{vol}(\Sphered) }{\nds} \diffsymbol \omega(\xi)=1.
\end{align}
The above seems like a rather expected localization property of the kernel $\mathcal{K}_{\mathrm{Harm}_{n}}$ for spherical harmonics. However, we are not aware of an explicit reference. Therefore, we prove it in the appendix \ref{app:auxproofs}, based on a modified calculation from \cite{Marzo2007}.

Finally, combining \eqref{eqnC4:esths2} and \eqref{eqnC4:concentrationonsphere}, and integrating over $x\in D$, we get
\begin{align}\label{eqnC4:lowesttwoeps}
	\lim_{n\to\infty}\lim_{m\to\infty}\int_{D}\int_{D} \frac{|{\Kmn}(x,y)|^2\diffsymbol y\, \diffsymbol x }{\nmd} & \geq (1-\epsilon_1) \int_{D_{\epsilon_2}} \widetilde{W}(x) \diffsymbol x,
\end{align}
where $D_{\epsilon_2}=\{x\in D: \|x-y\|\geq\epsilon_2 \textnormal{ for all }y\in\partial D\}$. The restriction to this subset on the right-hand side is necessary in order to guarantee the interchangeability of the limit and the outer integral (uniform convergence in Propositions \ref{prop:convergenceofChristoffelfunctionoftildePi} and \ref{prop:unilimit2} is only guaranteed for $x$ in compact subsets of the interior of $\BB^d\setminus\{0\}$). Letting $\epsilon_1, \epsilon_2$ tend to zero leads to the desired lower bound \eqref{eqnC4:lbths}.
\end{proof}

\subsection{Asymptotic eigenvalue distribution}

Now, having Theorem \ref{thm:EssentialEstimation2} established, the desired eigenvalue properties for the spatiospectral localization operator $\widetilde{\mathcal{SBS}}_{D,m,n}$ (reflecting the notion of bandlimit defined via separate limits on the radial and spherical contributions, with limits taken sequentially) follow directly from the general statement of Theorem \ref{thm:geneigdist}.

\begin{theorem}\label{thm:eigendistributionfortildePi}
Let $D\subset\B^d$ be a Lipschitz domain. Then, for any $0<\varepsilon<1/2$ and any $0<\tau<1$, it holds
\begin{align}\label{eqnC4:eigdis1_tildePi}
	\limif \lim_{m\to\infty} \frac{\sharp\{i:\varepsilon<\widetilde{\lambda}_i(D;m,n)<1-\varepsilon\}}{\nmd}&=0,
	\\\label{eqnC4:eigdis2_tildePi}
	\limif \lim_{m\to\infty}\frac{\sharp\{i:\tau<\widetilde{\lambda}_i(D;m,n)\leq 1\}}{\nmd}&=\int_{D}\widetilde{W}(x)\diffsymbol x,
\end{align}
with
\begin{align}
	\widetilde{W}(x)=\frac{2}{\pi \mathrm{vol}(\Sphered)}\frac{1}{\|x\|^{d-1}\sqrt{1-\|x\|^2}}.
\end{align}
It needs to be emphasized that the limits in \eqref{eqnC4:eigdis1_tildePi} and \eqref{eqnC4:eigdis2_tildePi} should be interpreted as sequential limits, i.e., first taking the limit over $m$  and then taking the limit over $n$.
\end{theorem} 


\section{Numerical examples for the space $\widetilde{\Pi}_{n,\kappa n}^d$: the linear coupling case}\label{sec:kappa}
The rigorous analysis of the sequential limit in the previous section provides clarity, but a scenario of greater practical relevance is the simultaneous constraint of radial and spherical bandwidths. This section numerically investigates such cases by focusing on linearly coupled bandwidths, defined by the relation
\begin{align}\label{eqn:coulingkappa}
	n_k=\kappa m_k.
\end{align}
for sequences $(m_k)_{k\in\mathbb{N}},(n_k)_{k\in\mathbb{N}}\subset\N_0$ with $\lim_{k\to\infty}m_k=\lim_{k\to\infty}n_k=\infty$ and a fixed ratio $\kappa>0$. As outlined in Section \ref{sec:framework}, this choice corresponds to spaces $\widetilde{\Pi}_{m_k, n_k}^d$ with a rectangular spectral shape.

\begin{figure}[t!]
	\centering
	\includegraphics[scale=0.28,angle=0]{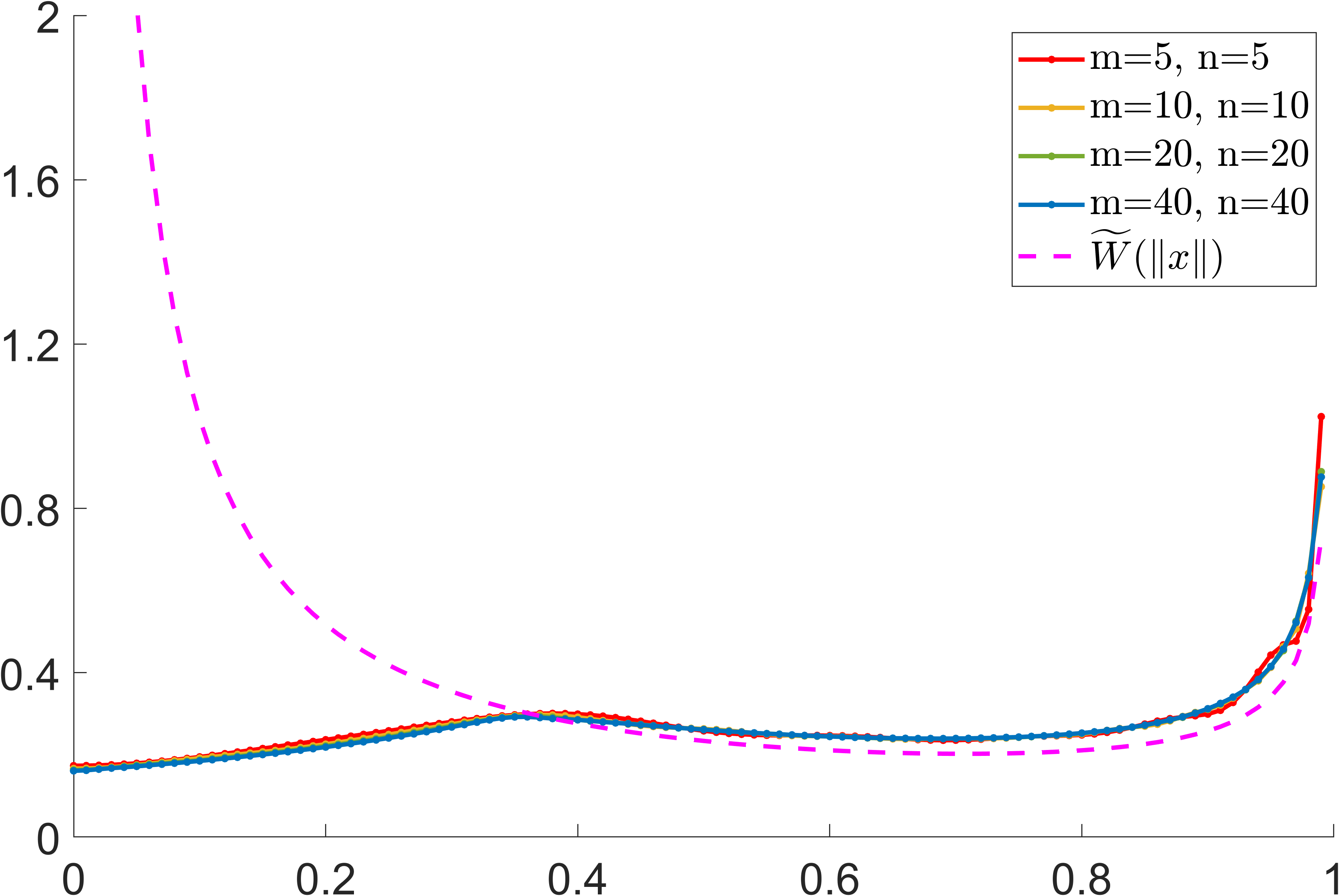}\quad		\quad
	\includegraphics[scale=0.28,angle=0]{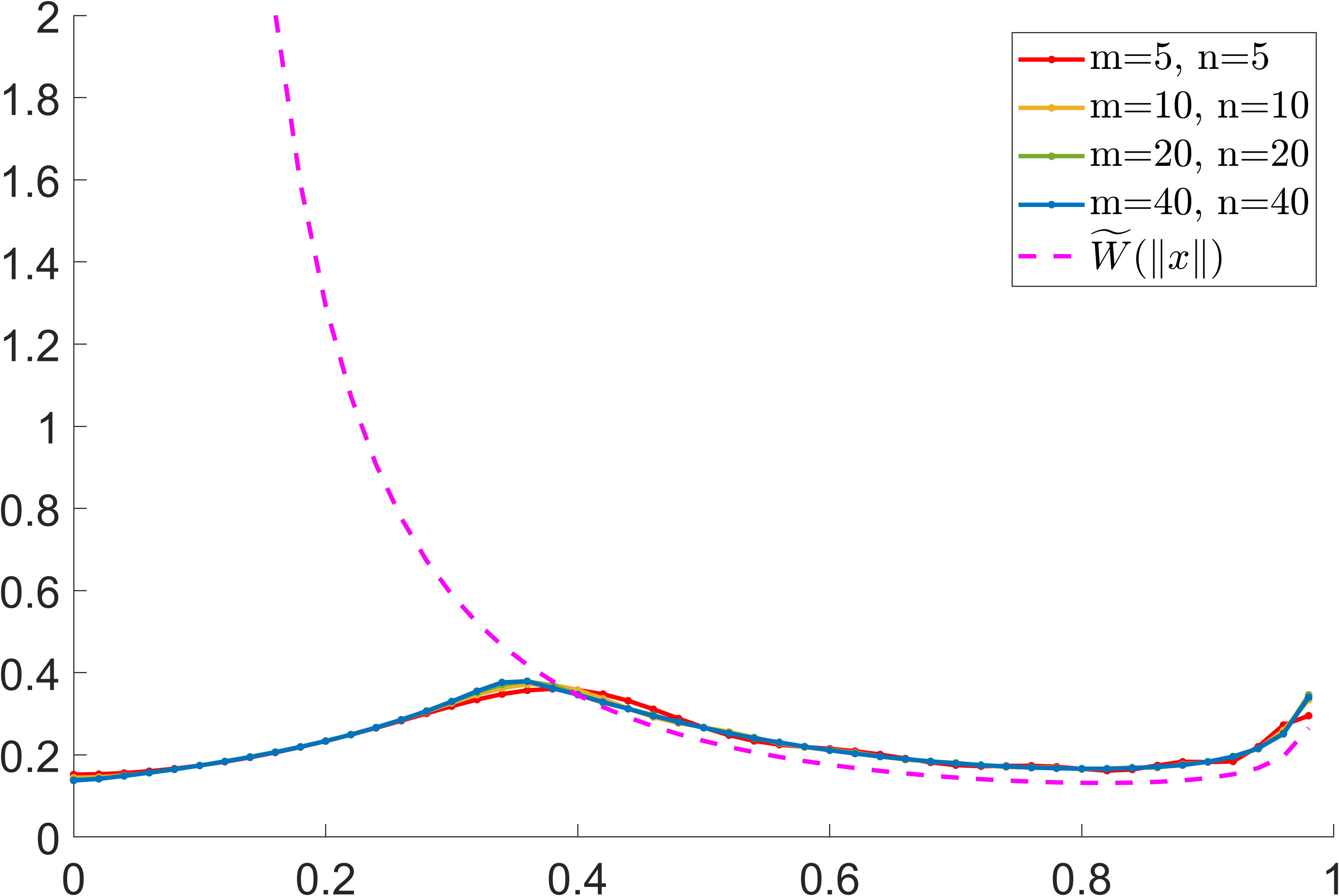}
	\\\includegraphics[scale=0.28,angle=0]{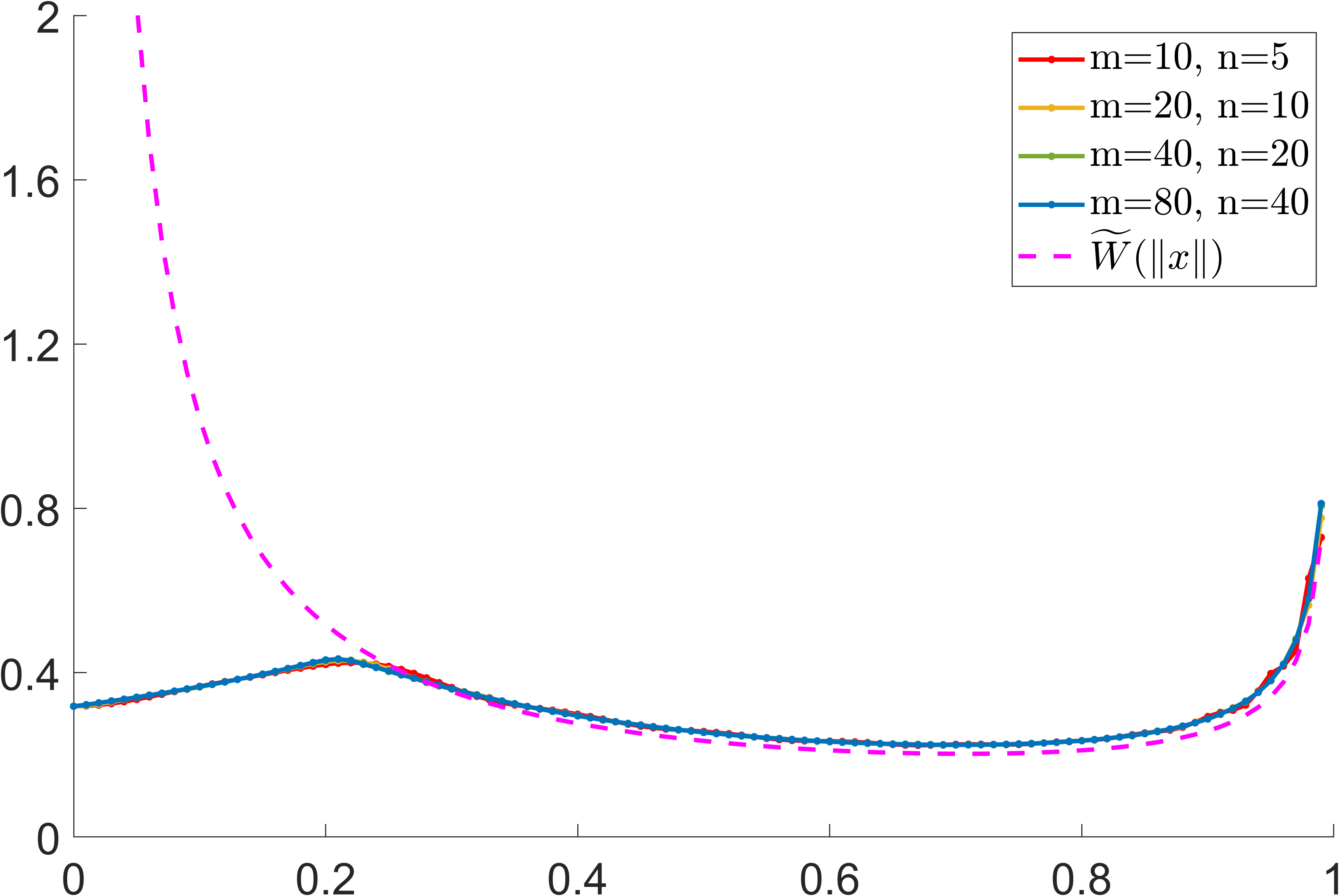}\quad	
	\quad
	\includegraphics[scale=0.28,angle=0]{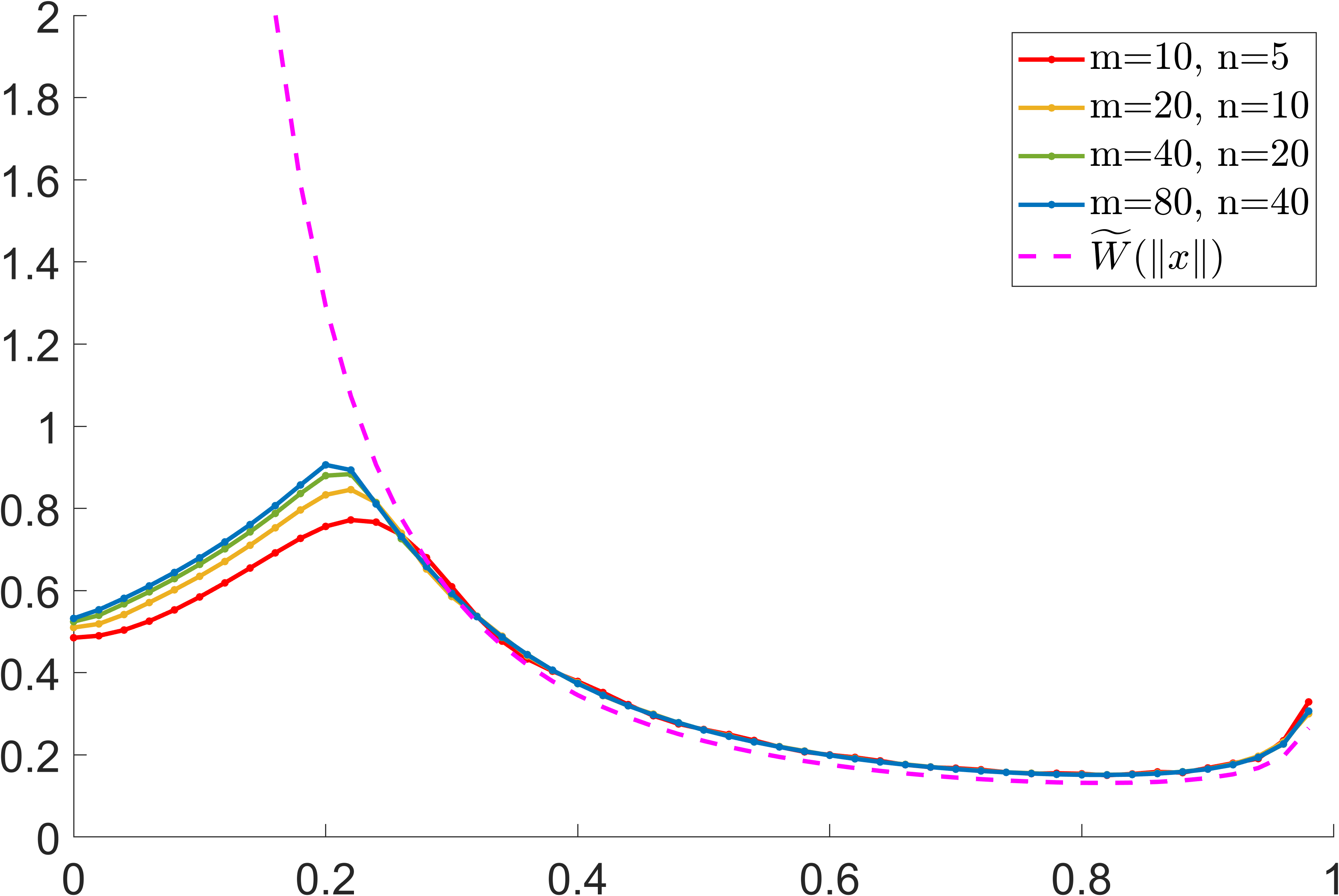}
	\\\includegraphics[scale=0.28,angle=0]{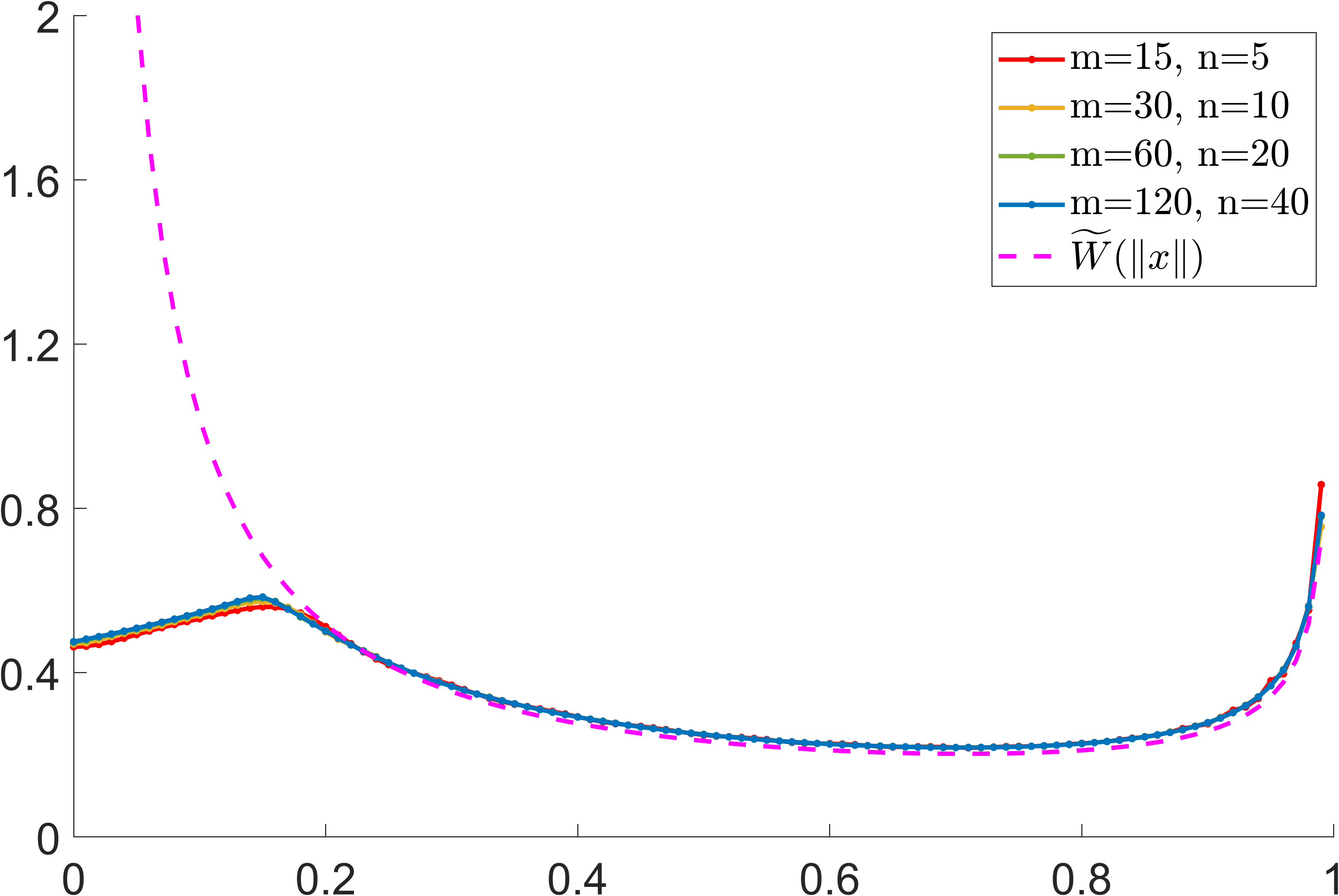}\quad	
	\quad
	\includegraphics[scale=0.28,angle=0]{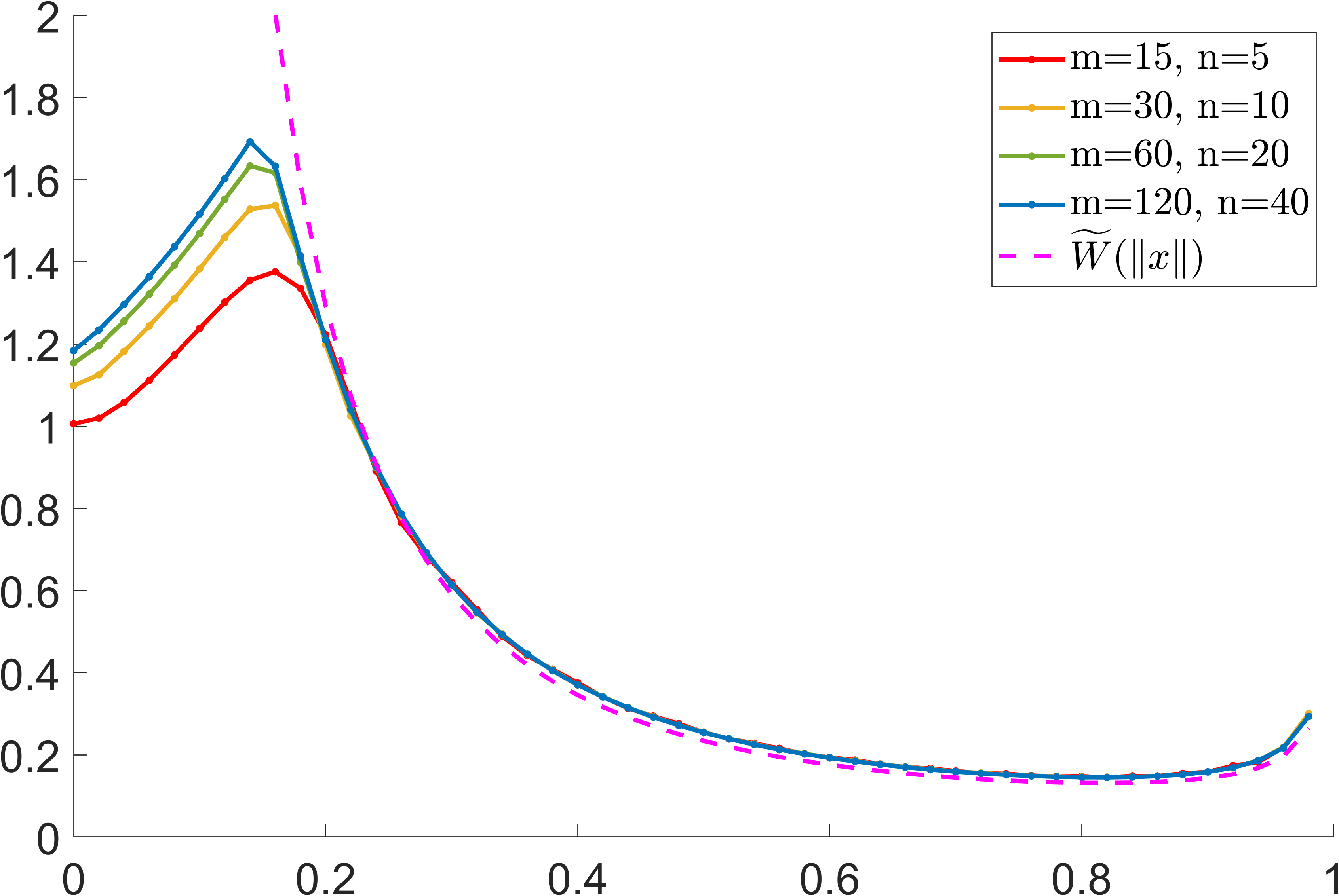}
	\caption{Illustration of $\widetilde{\mathcal{K}}_{m, n}(x,x)/\nmd$ for different ratios $\kappa$ that link $n=\kappa m$.  Left column: dimension $d=2$, right column: dimension $d=3$. The top, middle and bottom rows correspond to $\kappa=1$, $\kappa=\frac{1}{2}$ and $\kappa=\frac{1}{3}$, respectively. The dashed purple lines indicate the reference $\widetilde{W}$ from \eqref{eqn:wtilde1} (which differs for dimension $d=2$ and $d=3$).} 
	\label{fig:W_FB_Zernikes}
\end{figure}

We now examine the asymptotic behavior of the reproducing kernel diagonal via the normalized quantity
\begin{align}
	\frac{\widetilde{\mathcal{K}}_{m_k, n_k}(x,x)}{\widetilde{\mathcal{N}}_{m_k,n_k}^d}.
\end{align}
Figure \ref{fig:W_FB_Zernikes} illustrates this quantity for dimensions $d=2,3$ and several ratios $\kappa$.

Our numerical experiments support two key observations. First, the normalized kernel appears to converge to a limit function $\widetilde{W}_{\kappa}(x)$, which varies nontrivially with $\kappa$. Second, $\widetilde{W}_{\kappa}$ closely resembles the sequential-limit function $\widetilde{W}$ from \eqref{eqn:wtilde1}, with significant deviations primarily occurring near the origin. Moreover, as $\kappa$ decreases, $\widetilde{W}_{\kappa}$ converges toward $\widetilde{W}$.


\begin{remark}
	Based on the above observations, one may come up with the following conjectures: Let $x$ be in the interior of $\BB^d\setminus\{0\}$ (i.e., $0<\|x\|<1$). If $m_k$ and $n_k$ are linearly coupled via \eqref{eqn:coulingkappa} with $\kappa>0$, then a limiting function $\widetilde{W}_{\kappa}$ that satisfies $\widetilde{W}_{\kappa}(x)=\lim_{k\to\infty}\widetilde{\mathcal{K}}_{m_k, n_k}(x,x)/{\widetilde{\mathcal{N}}_{m_k,n_k}^d}$ exists and depends nontrivially on $\kappa$ (we do not have a guess for a closed-form expression of $\widetilde{W}_{\kappa}$). For the limiting cases $\kappa=0$ and $\kappa=\infty$, it holds
	\begin{align}
		\lim_{k\to\infty} \frac{\widetilde{\mathcal{K}}_{m_k,n_k}(x,x)}{\widetilde{\mathcal{N}}_{m_k,n_k}} &= \widetilde{W}(x),&\textnormal{if }n_k=o(m_k),
		\\\lim_{k\to\infty} \frac{\widetilde{\mathcal{K}}_{m_k,n_k}(x,x)}{\widetilde{\mathcal{N}}_{m_k,n_k}} &= 0,&\textnormal{if }m_k=o(n_k),
	\end{align}
	with $\widetilde{W}$ given by \eqref{eqn:wtilde1}. The statements in this section come purely from numerical observations, we have neither a proof of the existence of $\widetilde{W}_{\kappa}(x)$, nor a mathematical justification of its pattern. However, for the setup of a Fourier-Bessel basis instead of the Fourier-Jacobi basis used here, similar phenomena have been more rigorously analyzed in \cite{huang23}.
\end{remark}

\begin{remark}
	It is to note that the limiting behaviour of ${\widetilde{\mathcal{K}}_{m_k, n_k}(x,x)}/{\widetilde{\mathcal{N}}_{m_k,n_k}^d}$ may not only depend on the coupling between $m_k$ and $n_k$ but also on the choice of the parameter ${\supscriptInCFourOl}_j$ in the definition of the Fourier-Jacobi functions $\ZernikePolynomials_{i,j,\degb}$ from \eqref{def:ONS}, which we have chosen to be ${\supscriptInCFourOl}_j=j$ in the main body of this paper. For the case of a constant ${\supscriptInCFourOl}_j$, this is discussed in the appendix \ref{app:arbrhoj}.
\end{remark}


\section{Relation to indexing schemes of Zernike polynomials}\label{sec:zernike}

The 2-d Zernike polynomials, as an orthogonal function system on the disc, are widely used in optics, like wavefront analysis and lens aberration correction problems. In related applications, various indexing schemes are used to arrange the countably many basis functions. Naturally, Zernike polynomials can be indexed using 2-dimensional arrays. For example,  in \cite{Noll:76}, the Zernike polynomials are denoted by
\begin{align*}
	\nollZ_n^m:=	\begin{cases}
		\begin{rcases}  
			\sqrt{n+1}\nollR_n^m(r)\sqrt{2}\cos m\theta  & \\
			\sqrt{n+1}\nollR_n^m(r)\sqrt{2}\sin m\theta  &
		\end{rcases}    & m\neq 0 \\
		\sqrt{n+1} \nollR_n^{0}(r), & m=0
	\end{cases}
\end{align*}
with $n,m$ being non-negative integers satisfying $m\leq n$ and $n-m$ being even, and
\begin{align*}
	\nollR_n^m(r)=\sum_{s=0}^{(n-m)/2}\frac{(-1)^s(n-s)!}{s![(n+m)/2-s]![(n-m)/2-s]!}r^{n-2s}.
\end{align*}
Notice that, in this definition, every $\nollZ_n^m$ with $m\neq 0$ contains two functions: the azimuthal sine and cosine variants. The index $m$ stands for the azimuthal frequency and $n$ is called the radial degree (but note that $m$ and $n$ are actually coupled in the radial contribution). This definition corresponds to our definition of ${\ZernikePolynomials}_{i,j,\degb}$ from \eqref{def:ONS} with ${\supscriptInCFourOl}_j=j$, where $\ell$ distinguishes between the sine and cosine variants, $j$ reflects the azimuthal frequency, and $i$ denotes the radial degree.

However, for very practical reasons, 1-dimensional schemes are frequently used to order Zernike polynomials. For example, Noll himself has described the notation $\nollZ_j=\nollZ_n^m$ for $j\in\mathbb{N}$, where the transfer between $j$ and $m,n$ is given by
\begin{align*}
	& n=\lfloor \sqrt{(2j+\tfrac{1}{2})}-1\rfloor, \\
	&m=	\begin{cases}
		2\,\lfloor \frac{2j+1-n(n+1)}{4}\rfloor, &\textnormal{if $n$ is even},	\\
		2\,\lfloor \frac{2j+1-n(n+1)}{4}\rfloor	-1, &\textnormal{if $n$ is odd},
	\end{cases}
\end{align*}
and where $\nollZ_j$ with an even $j$ takes the azimuthal sine variant and  $\nollZ_j$ with an odd $j$ takes the cosine variant. Such a 1-dimensional indexing scheme naturally provides a way to define low-passes/bandlimits. In applications, only finitely many basis functions can be employed, and typically the first several functions of an indexing system would be used. The low-order Zernike polynomials in a indexing system can then be interpreted as more stable or significant modes, just like low-frequency terms in Fourier-transform-based signal processing. 

To formally express the relation between indexing schemes and our notions of bandlimit, we call a sequence of function spaces $\mathcal{F}_n$, $n\in\N_{0}$, and an indexing scheme $\nollZ_j$, $j\in\N$, \emph{compatible} if
\begin{align}
	\mathcal{F}_n=\mathrm{span}\{\nollZ_j:1\leq j\leq \mathrm{dim}(\mathcal{F}_n)\}, \quad \text{for all } n\in\N_{0}. 
\end{align}
It is, e.g., easy to verify that Noll's 1-dimensional indexing scheme from above is compatible with the polynomial spaces $\widehat{\Pi}_{n}^2$ from \eqref{eqn:hatpi}, \eqref{eqn:polyeq}.

In the review \cite{Niu2022}, the author listed six indexing schemes of 2-dimensional Zernike polynomials used in optics. There are several technical differences between these six indexing schemes, e.g, whether the functions are normalized, the direction of $\theta$ is clockwise/anticlockwise or the indexing starts from 0/1. In particular, the so-called Fringe/U of Arizona indexing system comprises 37 elements. However, if we focus on the notion of bandlimit implied by these schemes, it turns out that they can be grouped into two categories: those compatible with $\widehat{\Pi}_n^2$ from \eqref{eqn:hatpi}, \eqref{eqn:polyeq}, and those compatible with $\widecheck{\Pi}_n^2$ from \eqref{eqn:ckeckpi}. More specifically, the classification is given in the following diagram:
\begin{align*}
	\text{Indexing scheme of} \begin{cases}
		\begin{rcases}  \text{Noll} &  \\
			\text{OSA/ANSI} &  \\
			\text{Born and Wolf} &  \\
			\text{Malacara}
		\end{rcases} &   \text{compatible with } \widehat{\Pi}_n^2 \\[5ex]
		\begin{rcases} 
			\text{ISO-14999} &  \\
			\text{Fringe/U of Arizona} &
		\end{rcases} &   \text{compatible with } \widecheck{\Pi}_n^2.
	\end{cases}
\end{align*}
For the explicit definitions of the above mentioned indexing schemes and the discussions on their differences and applications, we refer the reader to the topical reviews \cite{Niu2022}, \cite{Jim17}, and references therein. 

\begin{figure}
	\centering
	\includegraphics[scale=0.5]{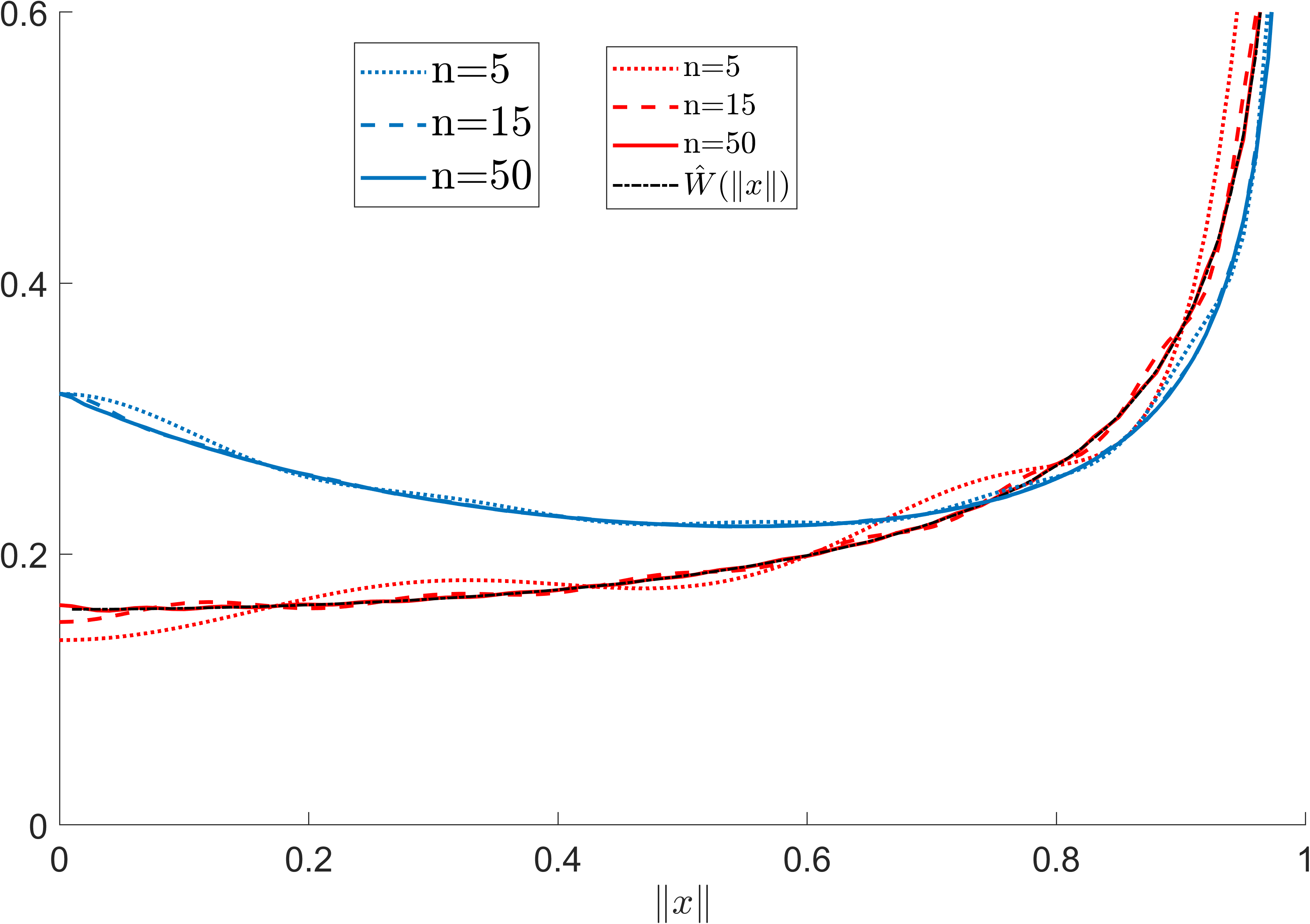}
	\caption{Illustration of ${\widehat{\mathcal{K}}_{n}(x,x)}/\widehat{\mathcal{N}}_n^2$ (red) and ${\widecheck{\mathcal{K}}_{n}(x,x)}/\widecheck{\mathcal{N}}_n^2$ (blue).} 
	\label{fig:W_optic_Zernikes}
\end{figure}

A comparison of ${\widehat{\mathcal{K}}_{n}(x,x)}/\widehat{\mathcal{N}}_n^2$ and ${\widecheck{\mathcal{K}}_{n}(x,x)}/\widecheck{\mathcal{N}}_n^2$ for different bandlimits $n$ is shown in Figure \ref{fig:W_optic_Zernikes}. While the figure suggests that both quantities converge to some limit function as $n$ tends to infinity, currently this has been proven rigorously only for ${\widehat{\mathcal{K}}_{n}(x,x)}/\widehat{\mathcal{N}}_n^2$ (cf. Proposition \ref{convergenceofChristoffelfunction}). Nonetheless, one can clearly see the different patterns of the graphs for the different notions of bandlimit provided by the spaces $\widehat{\Pi}_n^2$ and $\widecheck{\Pi}_n^2$, respectively. The two groups of graphs in Figure \ref{fig:W_optic_Zernikes} intersect at approximately $\|x\|\approx 0.75$, with the graphs corresponding to ${\widehat{\mathcal{K}}_{n}(x,x)}/\widehat{\mathcal{N}}_n^2$ taking lower values than those corresponding to ${\widecheck{\mathcal{K}}_{n}(x,x)}/\widecheck{\mathcal{N}}_n^2$ for $\|x\|\lessapprox 0.75$ and larger values for $\|x\|\gtrapprox 0.75$. One could interpret the consequences of that behaviour in the following way: when we take the same amount of (low-order) Zernike polynomials up to some bandwidth $n$ in the two categories of indexing schemes, the schemes that are compatible with $\widehat{\Pi}_n^2$ (e.g., Noll, OSA/ANSI, Born and Wolf, and Malacara) are expected to be more efficient for the representation of structures closer to the center of the disc (i.e., for $\|x\|\lessapprox 0.75$), while the schemes that are compatible with $\widecheck{\Pi}_n^2$ (e.g., ISO-14999 and Fringe/U of Arizona) are expected to be more efficient for the representation of structures closer to the boundary of the disc (i.e., for $\|x\|\gtrapprox 0.75$). 


	\section{Numerical examples of the eigenvalue distribution}\label{sec:num}

\begin{figure}[!t]
	\centering
	\includegraphics[scale=0.65]{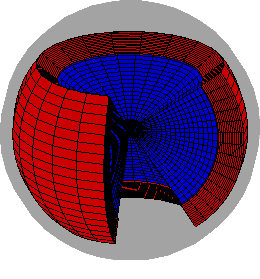}
	\caption{Illustration of the subdomains $D_1$ (blue) and $D_2$ (red) used in numerical experiments.}
	\label{fig:domains}
\end{figure}

In this section, we want to numerically illustrate 
the distribution features regarding the eigenvalues of the concentration operators discussed in Section \ref{sec:poly} and Section \ref{sec:seq} for the 3-d setup. 
For the spatial concentration regions, we select tesseroids
\begin{align*}
	D=\left\{x\in\BB^3 : r_1\leq r_x\leq r_2 \textnormal{ and } \xi_x=\begin{pmatrix}\sin\theta\cos\phi\\\sin\theta\sin\phi\\\cos\theta
	\end{pmatrix}, \theta_1\leq \theta \leq \theta_2,\phi_1\leq\phi \leq\phi_2\right\},
\end{align*}
with two different sets of parameters $r_1$, $r_2$, $\theta_1$, $\theta_2$, $\phi_1$, $\phi_2$. The first region $D_1$ (localized closer to the center of the ball) is represented by parameters $r_1=0.1$, $r_2=0.8$, $\theta_1=0.3\pi$, $\theta_2= 0.9\pi$, $\phi_1=-0.6\pi$, $\phi_2=0.9\pi$, the second region $D_2$ (localized closer to the surface) by parameters $r_1=0.7$, $r_2=0.9$, $\theta_1=0.3\pi$, $\theta_2= 0.9\pi$, $\phi_1=-0.6\pi$, $\phi_2=0.9\pi$. Thus, $D_1$ and $D_2$ probe concentration properties near the center and the boundary, respectively. Both tesseroids are illustrated in Figure \ref{fig:domains}.

\subsection{Comparison of $\widehat{\Pi}_n^3$ and $\widetilde{\Pi}_{m,n}^3$}

We consider the two bandlimited spaces $\widehat{\Pi}_n^3$ from \eqref{eqn:hatpi}, \eqref{eqn:polyeq} and $\widetilde{\Pi}_{m,n}^3$ from \eqref{eqn:tildekappapi} in dimension $d=3$. We want to empirically illustrate the behaviour of the corresponding relative Shannon numbers by
\begin{align}\label{eqn:numrelshan}
	\frac{\sum_{i=1}^{\widehat{\mathcal{N}}_n^3}\widehat{\lambda}_i(D;n)}{\widehat{\mathcal{N}}_n^3}\quad\textnormal{and}\quad\frac{\sum_{i=1}^{\widetilde{\mathcal{N}}_{m,n}^3}\widetilde{\lambda}_i(D;m,n)}{\widetilde{\mathcal{N}}_{m,n}^3},
\end{align}
for various bandlimits $m,n$. The asymptotic limits of the relative Shannon numbers for this setup have been shown in Theorems \ref{thm:1} and \ref{thm:eigendistributionfortildePi} to be $\int_{D}\widehat{W}(x)\,\diffsymbol x$ and $\int_{D}\widetilde{W}(x)\,\diffsymbol x$, with
\begin{align}\label{eqn:repwnum}
	\widehat{W}(x)=\frac{1}{\pi^2} \frac{1}{\sqrt{1-\|x\|^2}}\quad\textnormal{and}\quad	\widetilde{W}(x)=\frac{1}{2\pi^2} \frac{1}{\|x\|^2\sqrt{1-\|x\|^2}}.
\end{align}
As a reminder, the different behaviour of $\widehat{W}$ and $\widetilde{W}$ has already been shown in Figure \ref{fig:tildew1}. 

\begin{figure}[!t]
	\includegraphics[scale=0.28,angle=0]{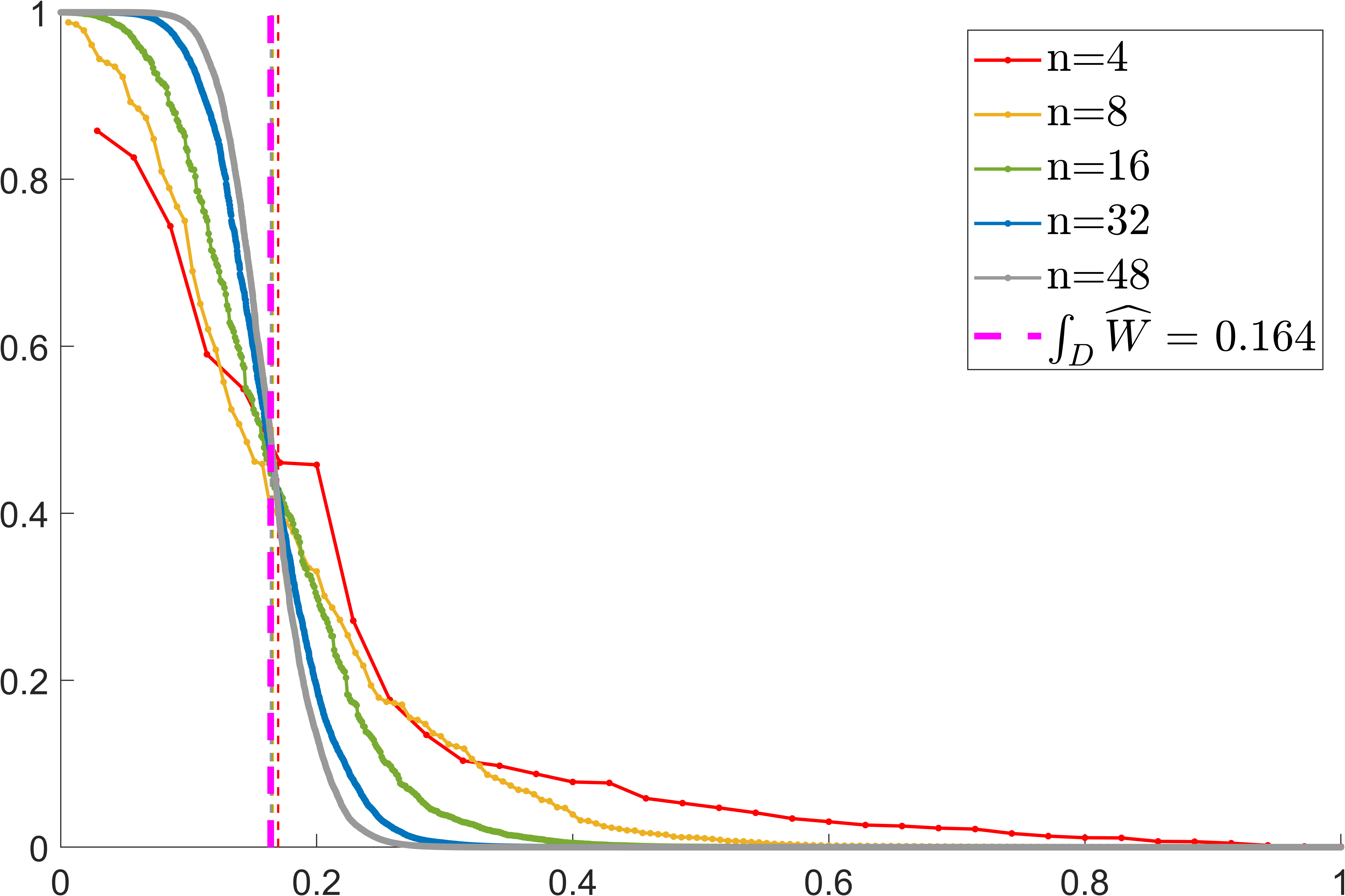}\quad\quad
	\includegraphics[scale=0.28,angle=0]{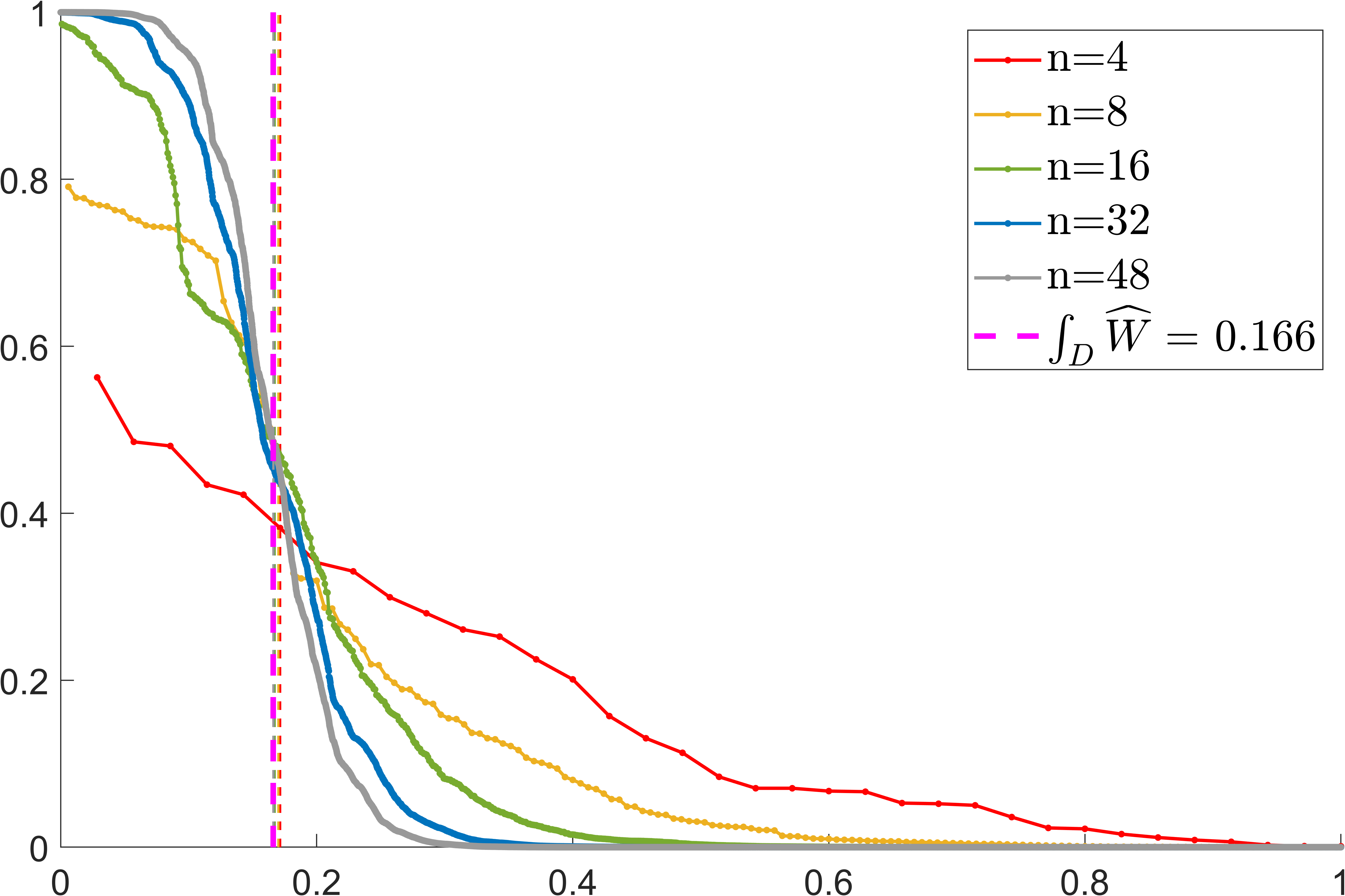}
	\\\includegraphics[scale=0.28,angle=0]{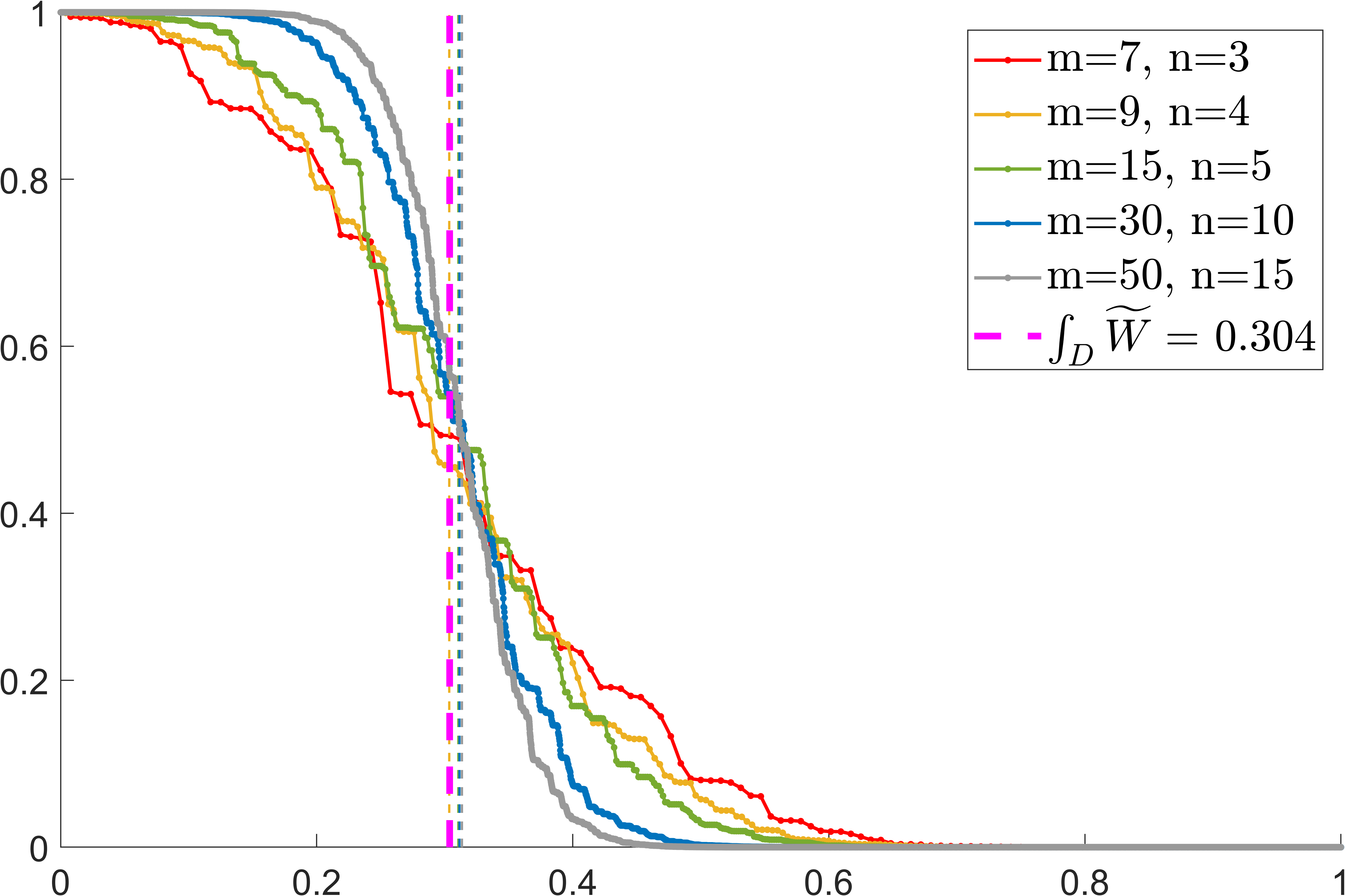}\quad	\quad
	\includegraphics[scale=0.28,angle=0]{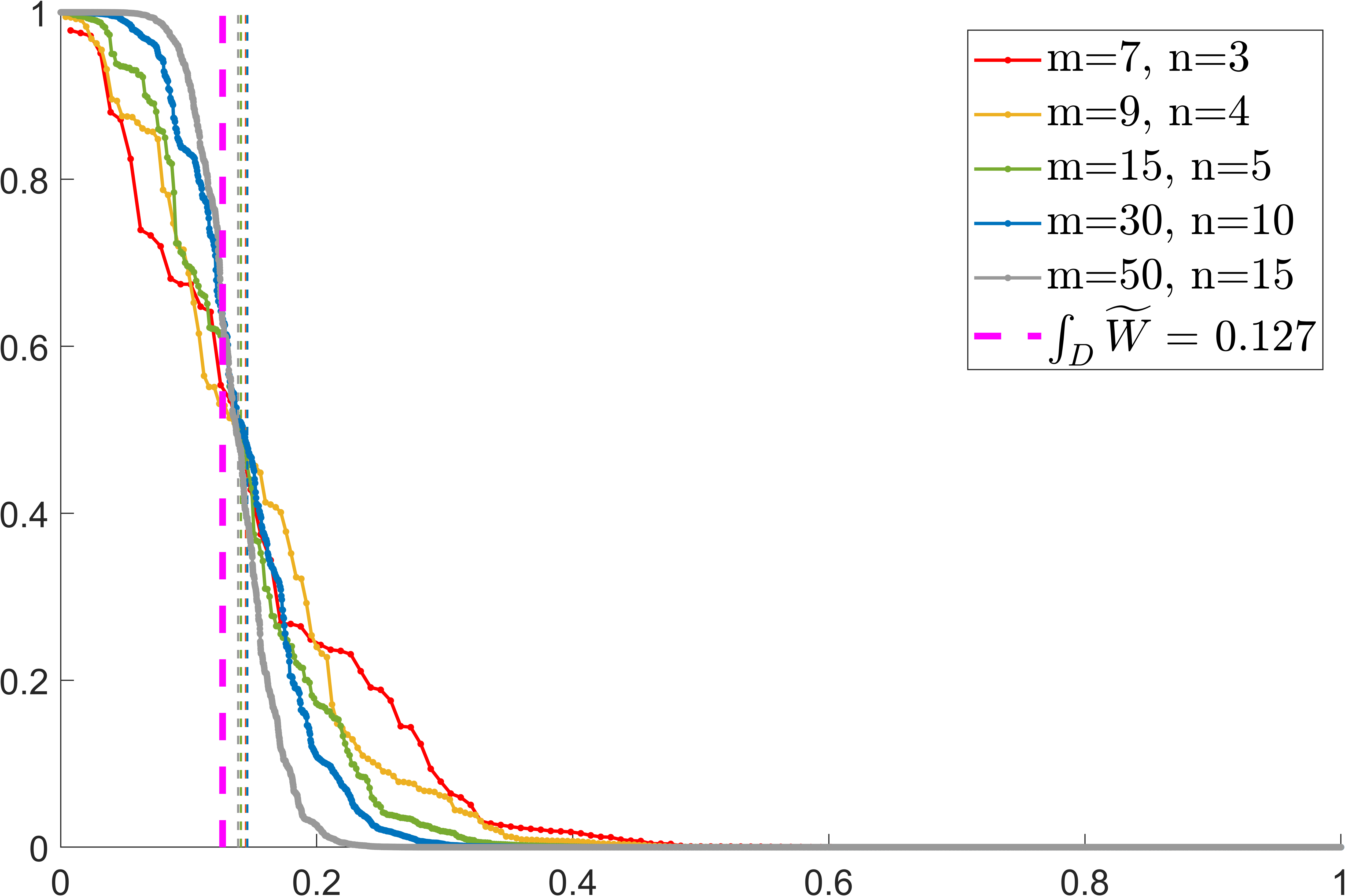}
	\caption{Illustration of the eigenvalue distributions $\widehat{\lambda}_i(D;n)$ (top row) and $\widetilde{\lambda}_i(D;m,n)$ (bottom) for the two subdomains $D=D_1$ (left column) and $D=D_2$ (right column).
	The solid curves demonstrate the distribution of eigenvalues in the relative x-coordinate, i.e., every single curve is composed of $(\frac{i}{\ndl},\lambda_i)$ corresponding to one concentration operator. 
	 Dotted vertical lines indicate the empirical relative Shannon numbers ${\sum_{i=1}^{\widehat{\mathcal{N}}_n^3}\widehat{\lambda}_i(D;n)}/{\widehat{\mathcal{N}}_n^3}$ and ${\sum_{i=1}^{\widetilde{\mathcal{N}}_{m,n}^3}\widetilde{\lambda}_i(D;m,n)}/{\widetilde{\mathcal{N}}_{m,n}^3}$. Dash vertical lines indicate the asymptotic relative Shannon numbers (asymptotic rSN) $\int_D \widehat{W}(x)\diffsymbol x$ and $\int_D\widetilde{W}(x)\diffsymbol x$. } 
	\label{fig:fig_eigen_distribution}
\end{figure}

If we interpret the relative Shannon number as the relative number of basis functions that is required to adequately represent a signal that is spatially localized in the domain $D$, then a smaller relative Shannon number indicates a more efficient representation of this localized signal for the underlying notion of bandlimit. The illustration of $\widehat{W}$ and $\widetilde{W}$ in Figure \ref{fig:tildew1} intuitively suggests that signals localized in $D=D_1$ (i.e., signals localized around the origin of the ball) should be more efficiently represented by the notion of bandlimit given by $\widehat{\Pi}_n^3$, while signals localized in $D=D_2$ (i.e., signals localized towards the surface of the ball) should be more efficiently represented by the notion of bandlimit given by $\widetilde{\Pi}_{m,n}^3$. The numerical results presented in Figure \ref{fig:fig_eigen_distribution} underpin this intuition: the numerically computed (empirical) relative Shannon numbers ${\sum_{i=1}^{\widehat{\mathcal{N}}_n^3}\widehat{\lambda}_i(D;n)}/{\widehat{\mathcal{N}}_n^3}$ are smaller than ${\sum_{i=1}^{\widetilde{\mathcal{N}}_{m,n}^3}\widetilde{\lambda}_i(D;m,n)}/{\widetilde{\mathcal{N}}_{m,n}^3}$ for $D=D_1$ while they are slightly larger for the domain $D=D_2$. This suggests that, if the approximate location of the signal of interest is known a priori, an adequate choice of the notion of bandlimit would potentially allow for a more efficient representation of that signal. And in order to judge which notion of bandlimit is particularly suitable for a given situation, the knowledge of the weight function $W$ is crucial. This was a main motivation for the study presented in this paper. 


\begin{remark}\label{rem:nujmdesc}
	The eigenvalues $\widehat{\lambda}_i(D;n)$ and $\widetilde{\lambda}_i(D;m,n)$ of the two spatiospectral localization operators $\widehat{\mathcal{B}\mathcal{S}\mathcal{B}}_{D,n}$ and $\widetilde{\mathcal{B}\mathcal{S}\mathcal{B}}_{D,m,n}$ are obtained by computing the eigenvalue decomposition of the matrices 
	\begin{align*}
		\widehat{T}_{D,n}=\left(\int_{D}{\ZernikePolynomials}_{i_1,j_1,\degb_{1}}(x){\ZernikePolynomials}_{i_2,j_2,\degb_{2}}(x)\diffsymbol x\right)_{0\leq2i_1+j_1\leq n, 1\leq \degb_1\leq 2j_1+1\atop0\leq 2i_2+j_2\leq n,1\leq \degb_2\leq 2j_2+1}\in \R^{[\frac{1}{6}(n+3)(n+2)(n+1)]^2},
	\end{align*}
	and
	\begin{align*}
		\widetilde{T}_{D,m,n}=\left(\int_{D}{\ZernikePolynomials}_{i_1,j_1,\degb_{1}}(x){\ZernikePolynomials}_{i_2,j_2,\degb_{2}}(x)\diffsymbol x\right)_{0\leq i_1\leq m,\, 0\leq j_1\leq n,\, 1 \leq \degb_1\leq 2j_1+1 \atop
			0\leq i_2\leq m,\, 0\leq j_2\leq n,\, 1 \leq \degb_2\leq 2j_2+1 }\in \R^{[(m+1)(n+1)^2]^2},
	\end{align*}
	respectively. For the case of $\widehat{\mathcal{B}\mathcal{S}\mathcal{B}}_{D,n}$ we ran the computations for various bandlimits up to a maximum of $n=48$, for the case of  $\widetilde{\mathcal{B}\mathcal{S}\mathcal{B}}_{D,m,n}$ we chose various bandlimits up to a maximum of $m=50$, $n=15$. Since we only consider the domains $D$ with a tesseroid shape, the integrals of the matrix entries are simplified via the separation of $r,\theta,\phi$ variables. All computations have been done in Matlab (mainly the Matlab default \emph{integral} and \emph{eig} functions for computing the integration and matrix eigendecomposition).
\end{remark}

\begin{remark} 
	Figure \ref{fig:fig_eigen_distribution} additionally illustrates that the empirical relative Shannon numbers converge to the derived asymptotic relative Shannon numbers as the bandlimit increases. However, for the case $\widetilde{\Pi}_{m,n}^3$, the empirical results seem to be slightly off from the indicated asymptotic result. This can be explained since the rigorous result for $\widetilde{\Pi}_{m,n}^3$ from Section \ref{sec:seq} only holds for sequential limits with respect to $m$ and $n$, which cannot be reasonably realized numerically. The bandlimits $m$ and $n$ that we have chosen for our illustrations are related by a ratio $\kappa$ varying roughly from $2$ to $3$. The numerical examples from Section \ref{sec:kappa} indicate that, within the regions $D_1$ and $D_2$, the asymptotic weight function $\widetilde{W}$ behaves very similar to the weight functions $\widetilde{W}_\kappa$ for this range of ratios $\kappa$. Due to these slight deviations, it is expected that the empirical relative Shannon numbers for $\widetilde{\Pi}_{m,n}^3$  in Figure \ref{fig:fig_eigen_distribution} are close to the indicated asymptotic relative Shannon number but can be    slightly off, opposed to the case $\widehat{\Pi}_{n}^3$ where the empirical and asymptotic results show a very strong agreement.
\end{remark}

	\section{Conclusion and outlook}\label{sec:con}
	This paper has introduced a general framework for Slepian-type spatiospectral concentration within the ball, centered on the concept of a spectral shape $\Omega$ for defining bandlimited spaces. Leveraging classical results on orthogonal algebraic polynomials, we have provided rigorous asymptotic analyses for two special cases: the polynomial space $\polynd$ and the sequential limit case of $\widetilde{\Pi}_{m,n}^d$. Our results establish the expected eigenvalue clustering and demonstrate that the Shannon number---the number of well-concentrated eigenfunctions---is characterized by distinct weight functions $\widehat{W}$ and $\widetilde{W}$, respectively. These analyses unequivocally reveal that different definitions of bandlimit lead to markedly different spatiospectral concentration properties, and further numerical experiments underpin the spectral shape's pivotal role in this characterization.
	
	From an applied perspective, the weight function $W$ in the Shannon number estimation indicates that, unlike the Euclidean and spherical cases, bandlimited spaces associated with the Fourier-Jacobi system possess a spatially varying resolution density. By altering the spectral shape, one can tune this spatial density pattern to a certain degree. We have explicitly connected this theory to practice by showing that common indexing schemes for Zernike polynomials correspond to specific spectral shapes, thereby providing a mathematical principle to study their varying efficacy for features of different spatial origin.
	
	In summary, the spectral shape $\Omega$ emerges as a fundamental design parameter for various Fourier-Jacobi function based signal processing on the ball. Future work will extend our analysis to arbitrary spectral shapes and develop efficient computational methods to leverage these tailored bases for real signal processing and inverse problems.

\section*{Acknowledgments}
 The authors are grateful to Dr. Alexander Kegeles and Dr. Laurent Baratchart for constructive comments on various drafts of this manuscript. We are also grateful to Prof. Volker Michel for pointing out the Fourier-Jacobi functions as well as for his hospitality during the second author's visit to the University of Siegen. Furthermore, we want to thank Prof. Heping Wang and Prof. Yongping Liu for valuable suggestions. This first author is supported by the Deutsche Forschungsgemeinschaft (DFG, German Research Foundation) under grant no. 551712198.

\bibliography{example}
\bibliographystyle{plain}

\appendix 

\section{General asymptotic eigenvalue distribution}\label{app:eigdist}

Assuming that a statement on the trace and Hilbert-Schmidt norm for the spatiospectral concentration operator of the form \eqref{asytr1}, \eqref{asyhs1} holds true, one can directly show the desired asymptotic eigenvalue distributions \eqref{eigdis11}, \eqref{eigdis21}. More precisely, the following statement holds.

\begin{theorem}\label{thm:geneigdist}
	Let $D$ be a subset of $\B^d$ and let there exist a function $W=W_\Omega$ in $L^2(\B^d)$, depending on the spectral shape $\Omega$, such that
	\begin{align}
		\limifl \frac{\mathrm{tr}(\tbtl)}{\ndl}&=\int_{D}W(x)\diffsymbol x, \label{eqn:atr2}
		\\\limifl \frac{\|\tbtl\|_{\HS}^{2}}{\ndl}&=\int_{D}W(x)\diffsymbol x.  \label{eqn:ahs2}
	\end{align}
	Then it holds for any $0<\varepsilon<1/2$ and any $0<\tau<1$ that
	\begin{align}
		\limifl \frac{\sharp\{i:\varepsilon<\lambda_i(D;L)<1-\varepsilon\}}{\ndl}&=0,\label{eqn:eig0}
		\\\limifl \frac{\sharp\{i:\tau<\lambda_i(D;L)\leq 1\}}{\ndl}&=\int_{D}W(x)\,\diffsymbol x. \label{eqn:eigw}
	\end{align}
\end{theorem}

\begin{proof}
	For any $0\leq a\leq b\leq1$, we use the abbreviation
	\begin{align}
		N_{L}(a,b)=\sharp\{i:a<\lambda_i(D;L)\leq b\}.
	\end{align}
	We begin with the proof of \eqref{eqn:eig0} and first observe that, for any $0<\varepsilon<1/2$ and $\varepsilon<t<1-\varepsilon$, it holds
	\begin{align}
		t-t^2\geq \varepsilon-\varepsilon^2=:C_{\varepsilon}>0.
	\end{align}
	Thus, by use of \eqref{eqn:atr2} and \eqref{eqn:ahs2}, we get
	\begin{align}\label{eqnC4:proofthm1labeledeq1} 
		C_{\varepsilon}\lim_{L\to\infty}\frac{N_L(\varepsilon,1-\varepsilon)}{\ndl}&\leq \lim_{L\to\infty}\sum_{i\leq\ndl}\frac{\lambda_i(D;L)-\lambda_i^2(D;L)}{\ndl} \\ \nonumber
		&=\lim_{L\to\infty}\frac{\mathrm{tr}(\tbtl)}{\ndl}-\lim_{L\to\infty}\frac{\|\tbtl\|_{\HS}^2}{\ndl}=0,
	\end{align}
	which is precisely the desired statement. Next, we turn to the proof of \eqref{eqn:eigw}. For any $\varepsilon>0$, we have
	\begin{align}\label{eqnC4:thm12est}
		(1-\varepsilon) N_{L}(1-\varepsilon,1)\leq \sum_{i\leq\ndl}\lambda_i(D;L)=\mathrm{tr}(\tbtl) \leq N_{L}(\varepsilon,1)+ \varepsilon\, \ndl
	\end{align}
	For $\varepsilon<\tau<1-\varepsilon$, one can further split $N_{L}(1-\varepsilon,1)$ and $N_{L}(\varepsilon,1)$ to reformulate \eqref{eqnC4:thm12est} in the following form:
	\begin{align}\label{eqnC4:proofthm1labeledeq2}
		(1-\varepsilon) \left(N_{L}(\tau,1)-N_{L}(\tau,1-\varepsilon)\right)\leq \mathrm{tr}(\tbtl)\leq \left(N_{L}(\tau,1)+N_{L}(\varepsilon,\tau)\right)+ \varepsilon\, \ndl.
	\end{align}
	Since \eqref{eqnC4:proofthm1labeledeq1} implies $\limifl{N_{L}(\varepsilon,\tau)}/{\ndl}=0$ and $\limifl{N_{L}(\tau,1-\varepsilon)}/{\ndl}=0$, the previous estimate and \eqref{eqn:atr2} yield 
	\begin{align}
		(1-\varepsilon)\limifl \frac{N_L(\tau,1)}{\nd}&\leq \limifl \frac{\mathrm{tr}(\tbtl)}{\ndl} =\int_{D}W(x)\,\diffsymbol x\leq \limifl \frac{N_L(\tau,1)}{\ndl}+\varepsilon.
	\end{align}
	Taking the limit $\varepsilon\to0$ concludes the proof.
\end{proof}

\section{Proofs from Section \ref{sec:seq}}\label{app:univlimit}


\begin{proof}[Proof of Proposition \ref{prop:convergenceofChristoffelfunctionoftildePi}]
	We set $r=\|x\|$ and $\xi=x/\|x\|$. Substituting $z=2r^2-1$, \eqref{eqnC4:reproducingkerneloftildePi} and \eqref{eqnC4:reproducingkernelofSH}, as well as $P_{j}^{(d)}(\xi\cdot\xi)=P_{j}^{(d)}(1)=1$ leads to
	\begin{align}\label{eqnC4:decompseofChristoffelfunctionoftildePi}
		&{\Kmn}(x,x)\\
		&=\nonumber \scalemath{0.95}{\sum_{j\leq n} \left(\sum_{i\leq m} \frac{4i+2\parzernike+d}{2^{{\supscriptInCFourOl}_j+\frac{d}{2}+1}} P_i^{0,\parzernike+\tfrac{d-2}{2}}(2r^2-1) P_i^{0,\parzernike+\tfrac{d-2}{2}}(2r^2-1)\right)2^{\parzernike+\frac{d}{2}+1}(r^2)^{\parzernike}\frac{\mathrm{dim}(H_j^d)}{\mathrm{vol}(\Sphered)}  }  \nonumber\\
		&=\sum_{j\leq n}R_m^j(z,z)\,(1+z)^{\parzernike+\tfrac{d-2}{2}} 4r^{2-d} \frac{\mathrm{dim}(H_j^d)}{\mathrm{vol}(\Sphered)},\nonumber
	\end{align}	
	with the auxiliary function
	\begin{align}\label{eqnC4:rmj}
		R_m^j(u,v)=\sum_{i\leq m} \frac{4i+2\parzernike+d}{2^{\parzernike+\frac{d}{2}+1}} P_i^{0,\parzernike+\tfrac{d-2}{2}}(u) \,P_i^{0,\parzernike+\tfrac{d-2}{2}}(v).
	\end{align}
	From \cite[Eq. (4.3.3)]{szego1975} we know that 
	\begin{align}
		\int_{-1}^1P_i^{0,\parzernike+\tfrac{d-2}{2}}(z)P_i^{0,\parzernike+\tfrac{d-2}{2}}(z)\,(1+z)^{\parzernike+\tfrac{d-2}{2}}\diffsymbol z= \frac{2^{\parzernike+\tfrac{d-2}{2}+1}}{2i+\parzernike+\tfrac{d-2}{2}+1}.
	\end{align}
	The above is precisely the inverse of the prefactor in the sum of \eqref{eqnC4:rmj}, which means that $R_m^j$ is the reproducing kernel of the Hilbert space $\Pi_m^1$ of univariate polynomials on the interval $[-1,1]$ of degree at most $m$, equipped with the inner product $\langle \cdot,\cdot\rangle_{w}$ for $w(z)=(1+z)^{\parzernike+\nicefrac{d-2}{2}}$. Hence, from \cite[Theorem 2.1, Remark (c)]{Lubinsky2009} and the regularity of the corresponding measure $\sigma_W$ (e.g., \cite[P. 101]{stahl_totik_1992}), we get
	\begin{align}\label{eqnC4:convergencofunivariateChristoffelfunction}
		\lim_{m\to\infty}\frac{1}{m+1} R_m^j(z,z)\,(1+z)^{\parzernike+\tfrac{d-2}{2}} = \frac{1}{\pi\sqrt{1-z^2}},
	\end{align}
	where the convergence holds pointwise for $-1<z<1$ and uniformly for $z$ in any compact subset of $(-1,1)$. The latter corresponds to pointwise convergence for $0<r<1$ and uniform convergence for $r$ in any compact subset of $(0,1)$.
	
	Combining \eqref{eqnC4:decompseofChristoffelfunctionoftildePi} and \eqref{eqnC4:convergencofunivariateChristoffelfunction}, and remembering \eqref{eqnC4:nmd}, we finally obtain
	\begin{align}
		\lim_{m\to\infty}\frac{{\Kmn}(x,x)}{\nmd}& =\sum_{j\leq n}
		\frac{4}{r^{d-2}}\frac{1}{\pi \sqrt{1-z^2}} \frac{\mathrm{dim}(H_j^d)}{ \nds \,\mathrm{vol}(\Sphered)}\\ \nonumber
		& = \frac{2}{\pi \mathrm{vol}(\Sphered)}\frac{1}{r^{d-1}\sqrt{1-r^2}}=\widetilde{W}(x).
	\end{align}
	The conditions on pointwise and uniform convergence in the statement of the Proposition follow directly from the corresponding conditions after \eqref{eqnC4:convergencofunivariateChristoffelfunction}.
\end{proof}

\begin{proof}[Proof of Proposition \ref{prop:unilimit2}]
	Similar to \eqref{eqnC4:reproducingkerneloftildePi} and \eqref{eqnC4:decompseofChristoffelfunctionoftildePi}, we can express ${\Kmn}$ as
	\begin{align}\label{eqnC4:decompseofChristoffelfunctionoftildePi2}
		&{\Kmn}(x,y)=\sum_{j\leq n}\frac{\mathrm{dim}(H_j^d)}{\mathrm{vol}(\Sphered)}P_{j}^{(d)}( \xi_x\cdot \xi_y)\, R_m^j(x,y)\left[(1+z_x)(1+z_y)\right]^{\frac{\parzernike}{2}+\tfrac{d-2}{4}} 4(r_xr_y)^{\tfrac{2-d}{2}},
	\end{align}	
	with $z_x=2r_x^2-1$, $z_y=2r_y^2-1$. A direct computation shows that
	\begin{align}\label{eqnC4:differenceofzandr}
		z_x-z_y=2r_x^2-2r_y^2=2(r_x-r_y)(r_x+r_y)=4r_x(r_x-r_y)-2(r_x-r_y)^2.
	\end{align}
	To keep the upcoming calculations a bit clearer, we write $x_{t,\xi}^m=y$ (noting that $y$ now depends on $t$, $\xi$, and $m$, without stating this explicitly via indices). In particular, it holds $r_y=r_x+\frac{t}{m+1}$. For a fixed $j$, the relations \eqref{eqnC4:convergencofunivariateChristoffelfunction} and \eqref{eqnC4:differenceofzandr} yield 
	\begin{align}
		(z_y-z_x)\left((1+z_x)^{\parzernike+\tfrac{d-2}{2}}R^j_m(z_x,z_x)\right) = \frac{4r_x t}{\pi\sqrt{1-z_x^2}}+\mathcal{O}(m^{-1}),
	\end{align}
	which holds uniformly for $x$ in $D$ and for $t$ in a compact subset of $\R$, as $m$ tends to infinity. This implies that $z_y=z_x+b/[(1+z_x)^{\parzernike+\nicefrac{d-2}{2}}R^j_m(z_x,z_x)]$, for $b=4r_x t/(\pi\sqrt{1-z_x^2})+\mathcal{O}(m^{-1})$, and that such $b$ lies in a compact subset of $\R$ if $x$ is in $D$ and $t$ is in a compact subset of $\R$.
	Application of the univariate universality limit from \cite[Thm. 1.1]{Lubinsky2009} (choosing $a=0$ and $b$ as above) leads to
	\begin{align}\label{eqnC4:universallimittildeK1}
		&\lim_{m\to\infty}\frac{[(1+z_x)(1+z_y)]^{\frac{\parzernike}{2}+\tfrac{d-2}{4}}R^j_m(z_x,z_y)}{(1+z_x)^{\parzernike+\tfrac{d-2}{2}}R^j_m(z_x,z_x)}
		\\& = \lim_{m\to\infty} \textnormal{sinc}\left( \pi \left( \frac{4r_x t}{\pi\sqrt{1-z_x^2}}+\mathcal{O}(m^{-1}) \right)\right)
		= \textnormal{sinc}\left(  \frac{4r_xt}{\sqrt{1-z_x^2}} \right)=\textnormal{sinc}\left( \frac{2t}{\sqrt{1-r_x^2}}\right),\nonumber
	\end{align}
	where the second equality follows from the continuity of the sinc function and the last equality comes from the substitution $z_x=2r_x^2-1$. Furthermore, we can invoke \eqref{eqnC4:convergencofunivariateChristoffelfunction} to obtain, for any $j,j'\in\N_0$, that
	\begin{align}\label{eqnC4:relatchristoffel}
		\lim_{m\to\infty}\frac{(1+z_x)^{{\parzernike}^{'}+\tfrac{d-2}{2}}R^{j'}_m(z_x,z_x)}{(1+z_x)^{\parzernike+\tfrac{d-2}{2}}R^{j}_m(z_x,z_x)}=1,
	\end{align}
	uniformly for $z_x$ in a compact subset of $(-1,1)$, i.e., uniformly for $x$ in $D$.
	
	We again want to point out that, due our notation $x_{t,\xi}^m=y$, the variable $y$ actually depends on $m$. In particular, it holds $r_y=r_x+\nicefrac{t}{m+1}$  and $\lim_{m\to\infty}r_y/r_x=1$, uniformly for $x$ in $D$ and for $t$ in a compact subset of $\R$. The latter together with \eqref{eqnC4:decompseofChristoffelfunctionoftildePi2}, \eqref{eqnC4:universallimittildeK1}, and \eqref{eqnC4:relatchristoffel} eventually yields 
	\begin{align}\label{eqnC4:universallimittildeK2}
		&\lim_{m\to\infty}	\frac{{\Kmn}(x,y)}{{\Kmn}(x,x)} \\
		&= \lim_{m\to\infty} \left(\frac{r_y}{r_x}\right)^{\tfrac{2-d}{2}} \frac{\sum_{j\leq n}[(1+z_x)(1+z_y)]^{\tfrac{\parzernike}{2}+\tfrac{d-2}{4}}R^j_m(z_x,z_y) \mathrm{dim}(H_j^d) P_j^{(d)}(\xi_x\cdot \xi)}{\sum_{j\leq n} (1+z_x)^{\parzernike+\tfrac{d-2}{2}}R^j_m(z_x,z_x) \mathrm{dim}(H_j^d) }  \nonumber 
		\\&= \lim_{m\to\infty} \left(\frac{r_y}{r_x}\right)^{\tfrac{2-d}{2}}\sum_{j\leq n}\Bigg\{ \frac{[(1+z_x)(1+z_y)]^{\tfrac{\parzernike}{2}+\tfrac{d-2}{4}}R^j_m(z_x,z_y) }{ (1+z_x)^{\parzernike+\tfrac{d-2}{2}}R^j_m(z_x,z_x) } \mathrm{dim}(H_j^d) P_j^{(d)}(\xi_x\cdot \xi)  \nonumber
		\\&\qquad\qquad\qquad\qquad\quad\quad\times\left(\sum_{j'\leq n}\frac{(1+z_x)^{{\parzernike}^{'}+\tfrac{d-2}{2}}R^{j'}_m(z_x,z_x)}{(1+z_x)^{\parzernike+\tfrac{d-2}{2}}R^{j}_m(z_x,z_x)} \mathrm{dim}(H_{j'}^d) \right)^{-1}\Bigg\}\nonumber
		\\&= \textnormal{sinc}\left(\frac{ 2t}{\sqrt{1-r_x^2}}\right)\frac{\sum_{j\leq n}\mathrm{dim}(H_{j}^d)\,P_j^{(d)}(\xi_x\cdot \xi)}{\sum_{j'\leq n}\mathrm{dim}(H_{j'}^d)} \nonumber
		\\&= \textnormal{sinc}\left(\frac{ 2t}{\sqrt{1-r_x^2}}\right)\frac{\mathrm{vol}(\Sphered)}{\nds}\mathcal{K}_{\mathrm{Harm}_{n}}(\xi_x,\xi).\nonumber
	\end{align}
	The above is exactly the desired relation \eqref{unilimit2}, with $r_x=\|x\|$ and $\xi_x=x/\|x\|$.
\end{proof}

\section{Auxiliary proofs}\label{app:auxproofs}

\begin{proof}[Proof of \eqref{eqn:polyeq}]
	Since, for ${\supscriptInCFourOl}_j=j$, it holds that ${\ZernikePolynomials}_{i,j,\degb}$ is a polynomial of degree $2i+j$, we first get $\text{span}\{{\ZernikePolynomials}_{i,j,\degb}:2i+j\leq n\}\subset \Pi_n^d$.  To prove the opposite inclusion, we only need to show that the dimensions of the spaces are the same. We know that the dimension of $\Pi_n^d$ is $\binom{n+d}{n}$. The dimension of $\text{span}\{{\ZernikePolynomials}_{i,j,\degb}:2i+j\leq n\}$ can be calculated as 
	\begin{align}\label{eqnC4:dimensionCountingOfTildePi}
		&\mathrm{dim}(\text{span}\{{\ZernikePolynomials}_{i,j,\degb}:2i+j\leq n\}) =\sum_{i=0}^n \sum_{j=0}^{n-2i}\mathrm{dim}(H_j^d) \\  \nonumber
		&=\sum_{j=0}^{n} \left( \left\lfloor \frac{n-j}{2}\right\rfloor+1\right) \left[\binom{j+d-1}{j}-\binom{j+d-3}{j-2}\right] \\\nonumber
		&=\sum_{j=0}^{n} \binom{j+d-1}{j} =\binom{n+d}{n}, \nonumber
	\end{align}
	which concludes the proof.
\end{proof}

\begin{proof}[Proof of \eqref{asytr}]
	We start with the following estimate from below
	\begin{align}
		\liminfif \frac{\mathrm{tr}(\widehat{\mathcal{SBS}}_{D,n})}{\nnd}&\geq\liminfif \int_{D\cap \BBCFour{1-\epsilon}}\frac{\Kn^{\mu}(x,x)}{\nnd}\diffsymbol \sigma_\mu(x)\label{eqnC4:lethm1}\\
		&=\int_{D\cap \BBCFour{1-\epsilon}}\frac{w_0(x)}{w_\mu(x)}w_\mu(x)\diffsymbol x = \int_{D\cap \BBCFour{1-\epsilon}}w_0(x)\diffsymbol x,\nonumber
	\end{align}
	and the following estimate from above
	\begin{align}
		&\limsupif \frac{\mathrm{tr}(\widehat{\mathcal{SBS}}_{D,n})}{\nnd}=\limsupif \int_{D}\frac{\Kn^{\mu}(x,x)}{\nnd}\diffsymbol \sigma_\mu(x)\label{eqnC4:uethm1}\\
		&\leq\limif \left(\int_{\BB^d}\frac{\Kn^{\mu}(x,x)}{\nnd}\diffsymbol \sigma_\mu(x)-\int_{\BBCFour{1-\epsilon}}\frac{\Kn^{\mu}(x,x)}{\nnd}\diffsymbol \sigma_\mu(x)+\int_{D\cap \BBCFour{1-\epsilon}}\frac{\Kn^{\mu}(x,x)}{\nnd}\diffsymbol \sigma_\mu(x)\right)\nonumber\\
		&= 1 -\int_{\BBCFour{1-\epsilon}}w_0(x)\diffsymbol x + \int_{D\cap \BBCFour{1-\epsilon}}w_0(x)\diffsymbol x.\nonumber
	\end{align}
	The exchange of integral and limit is allowed due to the uniform convergence of ${\Kn^{\mu}(x,x)}/{\nnd}$ on $\BBCFour{1-\epsilon}$, and the first term on the right-hand side of \eqref{eqnC4:uethm1} holds because
	\begin{align}	
		\int_{\BB^d}\Kn^{\mu}(x,x)\diffsymbol \sigma_\mu(x)=\sum_{k\leq n}\sum_{i\leq \nkkMone}\langle p_{i}^{k,\mu},p_{i}^{k,\mu}\rangle_{w_\mu} = \mathrm{dim}(\polynd)=\nnd.
	\end{align}
	Letting $\epsilon$ tend to zero, the estimates \eqref{eqnC4:lethm1}, \eqref{eqnC4:uethm1}, and the observation  $\int_{\BB^d}w_0(x)\diffsymbol x=1$ lead to \eqref{asytr}. 
\end{proof}

\begin{proof}[Proof of \eqref{eqnC4:ednorm}]
	Applying a polar coordinate transformation and using $J_{d/2}^{*}(0)=2^{-d/2}\Gamma^{-1}(\tfrac{d}{2}+1)$, we directly get
	\begin{align}\label{valueofs}
		e_d&=\int_{\mathbb{R}^d}\left|\frac{J_{d/2}^{*}(\|x\|)}{J_{d/2}^{*}(0)}\right|^2\diffsymbol x\\ \nonumber
		&=\mathrm{vol}(\Sphered)\int_{0}^{\infty}\left|t^{-\frac{d}{2}}2^{\frac{d}{2}}\Gamma(\tfrac{d}{2}+1)J_{d/2}(t)\right|^2t^{d-1}\diffsymbol t\\ \nonumber
		&=\mathrm{vol}(\Sphered)2^d|\Gamma(\tfrac{d}{2}+1)|^2\int_{0}^{\infty}t^{-1}|J_{d/2}(t)|^2\diffsymbol t\\ \nonumber
		&=\mathrm{vol}(\Sphered)2^{d-1}\Gamma(\tfrac{d}{2})\Gamma(\tfrac{d}{2}+1),
	\end{align}
	where the equality (see, e.g., \cite[Chap. 13.42, Eq. (1)]{Watson}, \cite[Eq. (2.7)]{Glasser1994})
	\begin{equation}\label{vals}
		\int_{0}^{\infty}t^{-1}|J_{d/2}(t)|^2 dt=\frac{\Gamma(\frac{d}{2})}{2\Gamma(\frac{d}{2}+1)}=\frac{1}{d}
	\end{equation} 
	has been used in the last line of \eqref{valueofs}.
	From \eqref{valueofs}, we already get that $e_d$ is finite, but using the Legendre duplication formula and some basic computations for the Gamma function one get the desired estimation~\eqref{eqnC4:ednorm}
\end{proof}

\begin{proof}[Proof of \eqref{eqnC4:concentrationonsphere}]
	We first gather some preliminary properties. Using $\lambda=(d-3)/2$ as an abbreviation, $\mathcal{K}_{\mathrm{Harm}_{n}}$ is known to have the closed form representation (e.g., \cite[P. 565]{Marzo2007}) 
	\begin{align}\label{eqnC4:closedformofreproducingkernelofSH}
		\mathcal{K}_{\mathrm{Harm}_{n}}(\xi_x,\xi)=\frac{C_{d,n}}{\mathrm{vol}(\Sphered)}P_n^{1+\lambda,\lambda}(\xi_x\cdot\xi),
	\end{align}
	for $\xi_x,\xi\in\Sphered$ and $C_{d,n}=\binom{n+d-2}{n}/\binom{n+\frac{d-3}{2}}{n}$. Asymptotically, it holds $C_{d,n}\simeq n^{(d-1)/2}$. Letting $c>0$ be some fixed constant, we furthermore get from \cite[P. 198]{szego1975} that, for $\frac{c}{n}\leq \theta \leq \pi-\frac{c}{n}$,
	\begin{align}\label{eqnC4:estimateofJacobiP1}
		P^{\lambda+1,\lambda}_n(\cos \theta)=\frac{k(\theta)}{\sqrt{n}}\left(\cos\left((n+\lambda+1)\theta-\frac{\pi}{2}\left(\lambda+\frac{3}{2}\right)\right)+\frac{\mathcal{O}(1)}{n\sin\theta}\right)\textnormal{ for }n\to\infty,
	\end{align}
	where $k(\theta)=\pi^{-1/2}(\sin\frac{\theta}{2})^{-\lambda-3/2}(\cos\frac{\theta}{2})^{-\lambda-1/2}$. Finally, by the symmetry relation of Jacobi polynomials and the Mehler-Heine formula (e.g., \cite[Thm. 8.1.1]{szego1975}), it holds
	\begin{align}\label{eqnC4:estimateofJacobiP2}
		\lim_{n\to\infty}(-1)^{n}n^{-\lambda}P_n^{\lambda+1,\lambda}\left(\cos\left(\pi-\frac{z}{n}\right)\right)&=\lim_{n\to\infty}n^{-\lambda}P_n^{\lambda,\lambda+1}\left(\cos\left(\frac{z}{n}\right)\right)\\ \nonumber
		&=\left(\frac{z}{2}\right)^{-\lambda}J_{\lambda}(z)= 2^{\lambda}J_{\lambda}^{*}(z),
	\end{align}
	uniformly for $z$ in a bounded subset of $\R$. 
	
	Now we can begin with the actual proof. To remind the reader: the equality to prove (i.e., \eqref{eqnC4:concentrationonsphere}) states
	\begin{align}\label{eqnC4:profa2}
		\lim_{n\to\infty}\int_{\mathcal{C}_{\epsilon}(\xi_x)}|\mathcal{K}_{\mathrm{Harm}_{n}}(\xi_x,\xi)|^2 \frac{\mathrm{vol}(\Sphered) }{\nds} \diffsymbol \omega(\xi)=1.
	\end{align}
	From the reproducing kernel property of $\mathcal{K}_{\mathrm{Harm}_{n}}$, we get
	\begin{align*}
		\int_{\Sphered}|\mathcal{K}_{\mathrm{Harm}_{n}}(\xi_x,\xi)|^2 \frac{\mathrm{vol}(\Sphered) }{\nds}\diffsymbol \omega(\xi)=\frac{\mathrm{vol}(\Sphered) }{\nds}\mathcal{K}_{\mathrm{Harm}_{n}}(\xi_x,\xi_x)=1.
	\end{align*} 
	Thus, in order to prove \eqref{eqnC4:profa2}, it suffices to show 
	\begin{align}\label{eqnC4:equivalentEqconcentrationonsphere0}
		\lim_{n\to\infty} \int_{\Sphered\setminus \mathcal{C}_{\epsilon}(\xi_x)}|\mathcal{K}_{\mathrm{Harm}_{n}}(\xi_x,\xi)|^2 \frac{\mathrm{vol}(\Sphered) }{\nds} \diffsymbol \omega(\xi)=0.
	\end{align}
	This is what we will do next. Observing \eqref{eqnC4:nsd} and the asymptotic behaviour of $C_{d,n}$, we get $C_{d,n}^2/\nds=\mathcal{O}(1)$, as $n$ tends to infinity. The equality \eqref{eqnC4:closedformofreproducingkernelofSH} and a parametrization of the integral in terms of spherical coordinates then lead us to
	\begin{align}\label{eqnC4:equivalentEqconcentrationonsphere}
		&\int_{\arccos (1-\epsilon)}^{\pi}  |P^{1+\lambda,\lambda}_n(\cos\theta)|^2(\sin\theta)^{d-2}\, \diffsymbol \theta \nonumber \\
		&\quad\quad\simeq  c_d\int_{\Sphered\setminus \mathcal{C}_{\epsilon}(\xi_x)}|\mathcal{K}_{\mathrm{Harm}_{n}}(\xi_x,\xi)|^2 \frac{\mathrm{vol}(\Sphered) }{\nds} \diffsymbol \omega(\xi),
	\end{align}
	for $n\to\infty$. Throughout this proof, $c_d>0$ denotes some generic constant that only depends on the dimension $d$ and that may change at every appearance. The constant $c>0$, however, that has already be mentioned in the beginning, is fixed throughout. We split the integration over $[\arccos (1-\epsilon),{\pi}]$ on the left-hand side of \eqref{eqnC4:equivalentEqconcentrationonsphere} into three subintervals $[\arccos(1-\epsilon), \frac{\pi}{2}]$, $[\frac{\pi}{2}, \pi-\frac{c}{n}]$, and $[\pi-\frac{c}{n},\pi]$. For $\theta$ in $[\arccos(1-\epsilon), \frac{\pi}{2}]$, the function $k(\theta)$ in \eqref{eqnC4:estimateofJacobiP1} is bounded and, therefore, it holds uniformly that $|P^{1+\lambda,\lambda}_n(\cos(\theta))|^2= \mathcal{O}(n^{-1})$, as $n$ tends to infinity. We get
	\begin{align}\label{eqnC4:equivalentEqconcentrationonsphere1}
		\lim_{n\to\infty}\int_{\arccos (1-\epsilon)}^{\frac{\pi}{2}}  |P^{1+\lambda,\lambda}_n(\cos\theta)|^2(\sin\theta)^{d-2}\, \diffsymbol \theta =0.
	\end{align}
	On the interval  $[\frac{\pi}{2}, \pi-\frac{c}{n}]$, we have $k(\theta)\leq c_d (\cos\frac{\theta}{2})^{-\lambda-\frac{1}{2}}\leq c_d (\pi-\theta)^{-\lambda-\frac{1}{2}}$ and $\sin \theta\leq \pi-\theta$. Thus, 
	\begin{align}\label{eqnC4:equivalentEqconcentrationonsphere2}
		\lim_{n\to\infty}\int_{\frac{\pi}{2}}^{\pi-\frac{c}{n}}  |P^{1+\lambda,\lambda}_n(\cos\theta)|^2(\sin\theta)^{d-2} \diffsymbol \theta&\leq c_d\lim_{n\to\infty}  \int_{\frac{\pi}{2}}^{\pi-\frac{c}{n}} \frac{(\pi-\theta)^{-2\lambda-1}}{n} (\pi-\theta)^{d-2} \diffsymbol \theta\nonumber
		\\ &=c_d\lim_{n\to\infty}  \int_{\frac{\pi}{2}}^{\pi-\frac{c}{n}}  \frac{1}{n}\, \diffsymbol \theta=0.
	\end{align}
	On the interval $[\pi-\frac{c}{n},\pi]$, we substitute $\theta=\pi-\frac{z}{n}$ and get
	\begin{align*}
		\int_{\pi-\frac{c}{n}}^{\pi}  |P^{1+\lambda,\lambda}_n(\cos\theta)|^2(\sin\theta)^{d-2}  \diffsymbol \theta &= \int_0^c \left|P_n^{1+\lambda,\lambda}\left(\cos\left(\pi-\frac{z}{n}\right)\right)\right|^2\left(\sin\frac{z}{n}\right)^{d-2}\frac{1}{n}\diffsymbol z\\ \nonumber
		&\leq \frac{1}{n^2}\int_0^c \frac{1}{n^{d-3}}\left|P_n^{1+\lambda,\lambda}\left(\cos\left(\pi-\frac{z}{n}\right)\right)\right|^2 z^{d-2}\diffsymbol z .
	\end{align*}
	The relation \eqref{eqnC4:estimateofJacobiP2} yields that the integrand uniformly converges to $2^{2\lambda}z|J_{\lambda}(z)|^2$ as $n$ tends to infinity. The latter is integrable by \cite[Chap. 13.42, Eq. (1)]{Watson}
	. Thus, we obtain
	\begin{align}\label{eqnC4:equivalentEqconcentrationonsphere3}
		\lim_{n\to\infty}\int_{\pi-\frac{c}{n}}^{\pi}  |P^{1+\lambda,\lambda}_n(\cos\theta)|^2(\sin\theta)^{d-2} \diffsymbol \theta =0.
	\end{align}
	Combining the results \eqref{eqnC4:equivalentEqconcentrationonsphere}-\eqref{eqnC4:equivalentEqconcentrationonsphere3} leads to the desired statement \eqref{eqnC4:equivalentEqconcentrationonsphere0}, which eventually concludes the proof.
\end{proof}

\section{Numerical illustration of the eigenvalue distribution: comparison of $\Pi_{\bw,\Omega}^{d}$ for two nonstandard spectral shapes 
}

In \eqref{eqn:pilom} we have provided a fairly general definition of the bandlimited function space $\Pi_{\bw,\Omega}^{d}$, where the notion of bandlimit is based on the spectral shape $\Omega$. So far, we have restricted our investigations to the special cases \eqref{eqn:nmspace}-\eqref{eqn:tildekappapi}. Here, we briefly want to provide some illustrations for two less common choices of $\Omega$, namely,
\begin{align}
	\Omega_1&=\{(x,y)\in[0,1]^2:x^2+y^2\leq 1\},\label{eqn:shape1}
	\\\Omega_2&=\{(x,y)\in[0,1]^2:(1-x)^2+(1-y)^2\geq 1\}.\label{eqn:shape2}
\end{align}

\begin{figure}
	\centering
	\includegraphics[scale=0.25]{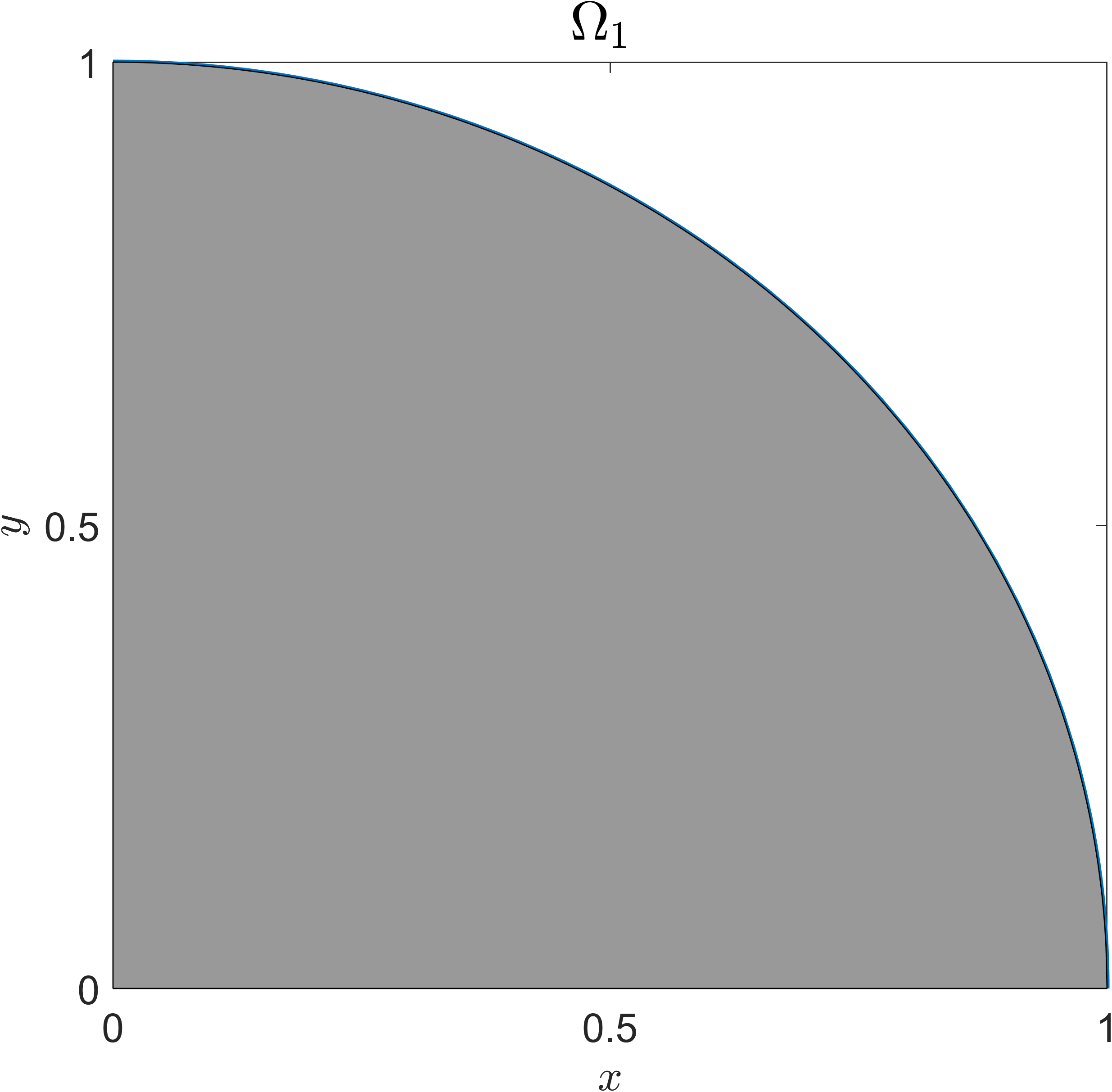}\qquad\qquad\includegraphics[scale=0.25]{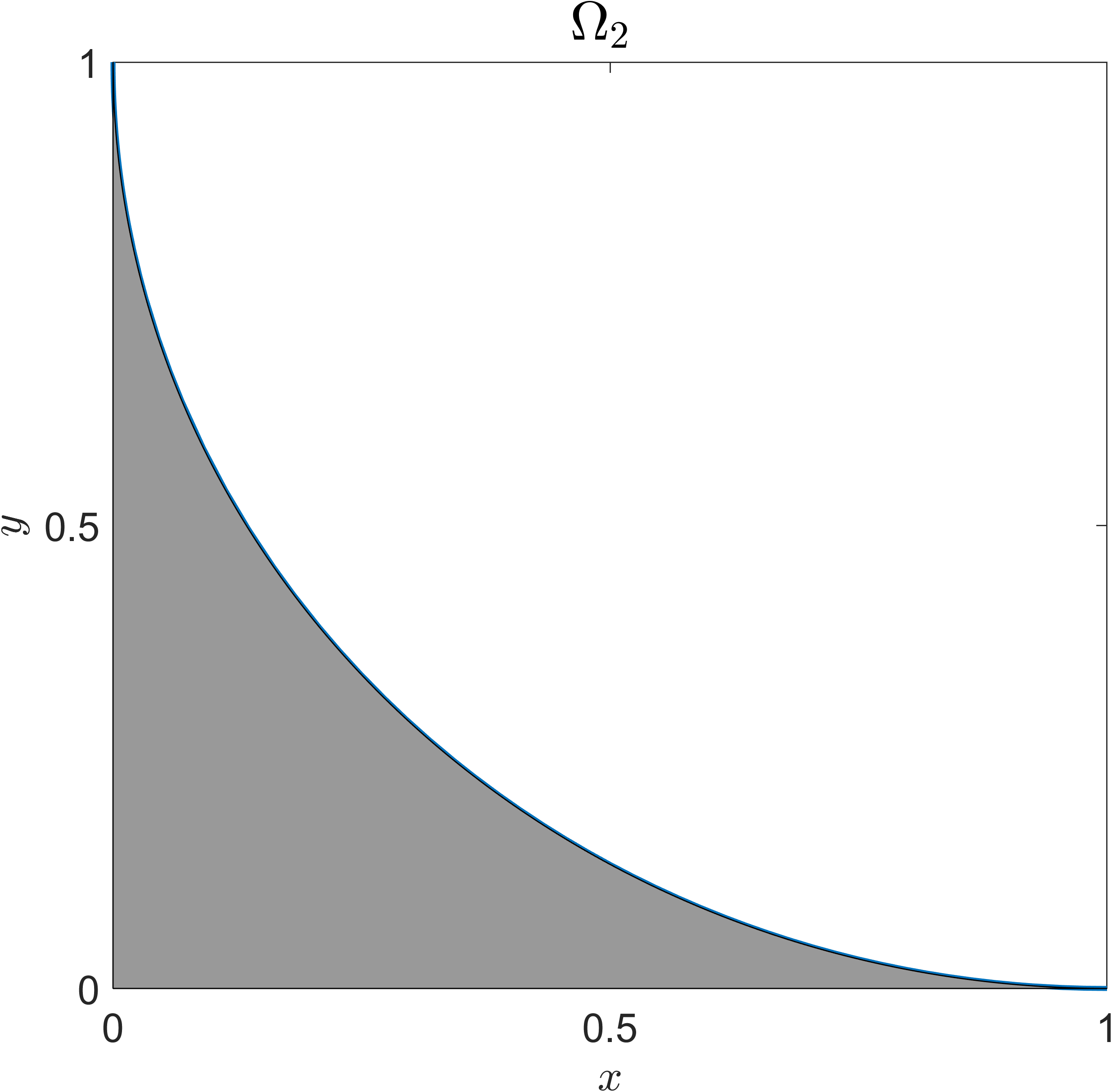}
	\caption{Illustration of the two spectral shapes $\Omega_1$ and $\Omega_2$ from \eqref{eqn:shape1} and \eqref{eqn:shape2}.}\label{fig:omegaillust}
\end{figure}

\begin{figure}[!t]
	\centering
	\includegraphics[scale=0.27,angle=0]{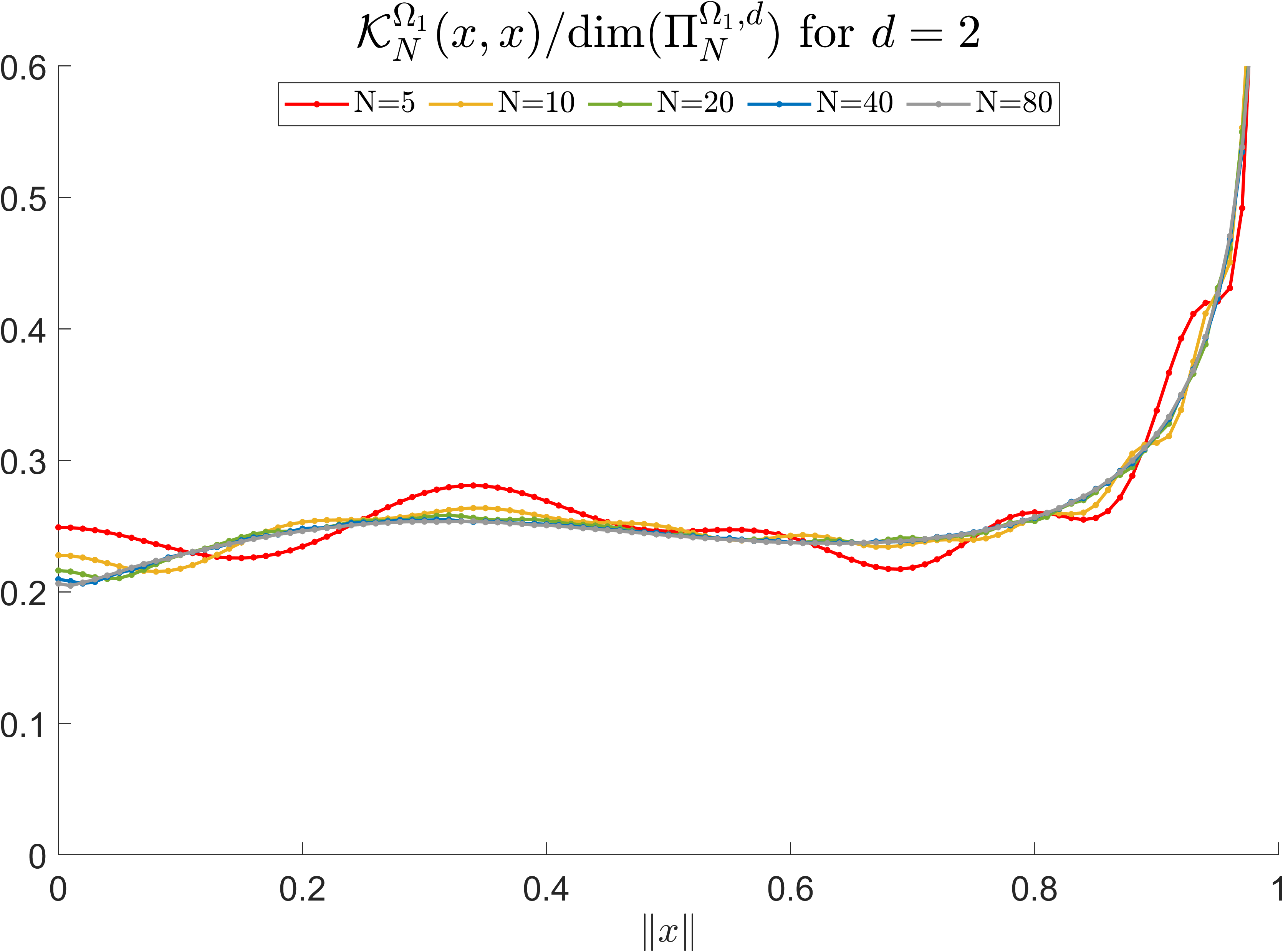}\quad\quad\quad
	\includegraphics[scale=0.27,angle=0]{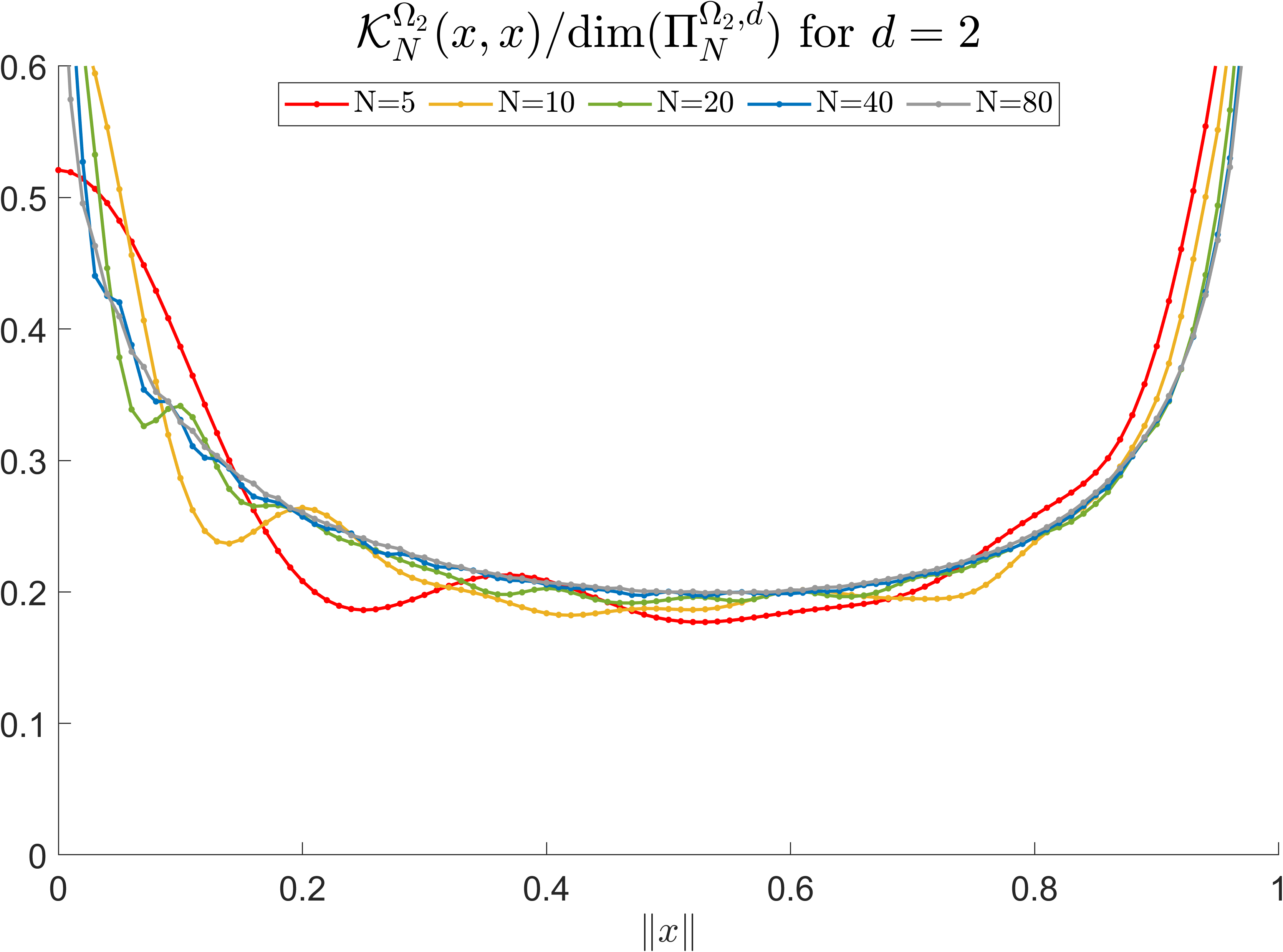}
	\\\includegraphics[scale=0.27,angle=0]{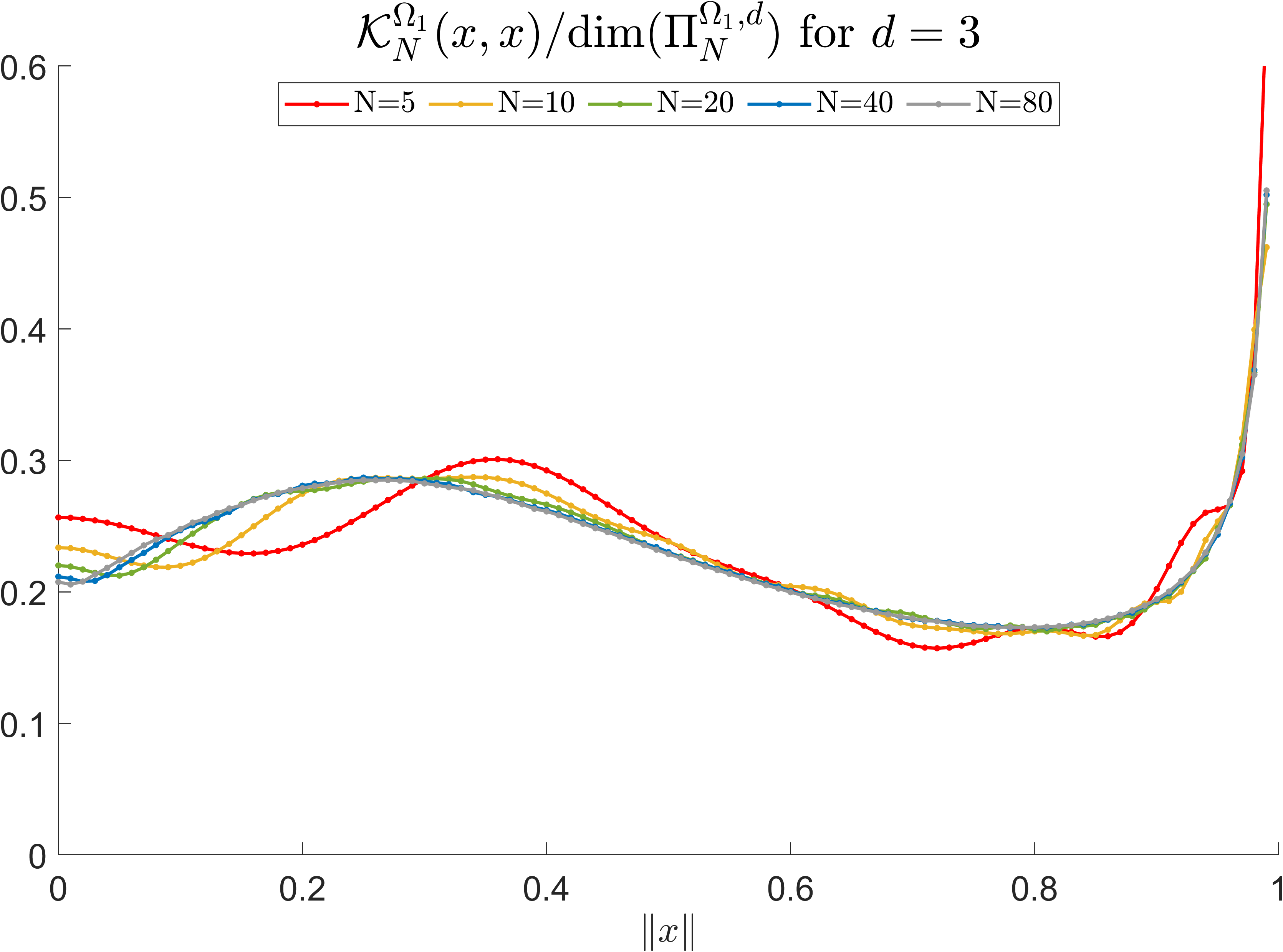}\quad\quad\quad
	\includegraphics[scale=0.27,angle=0]{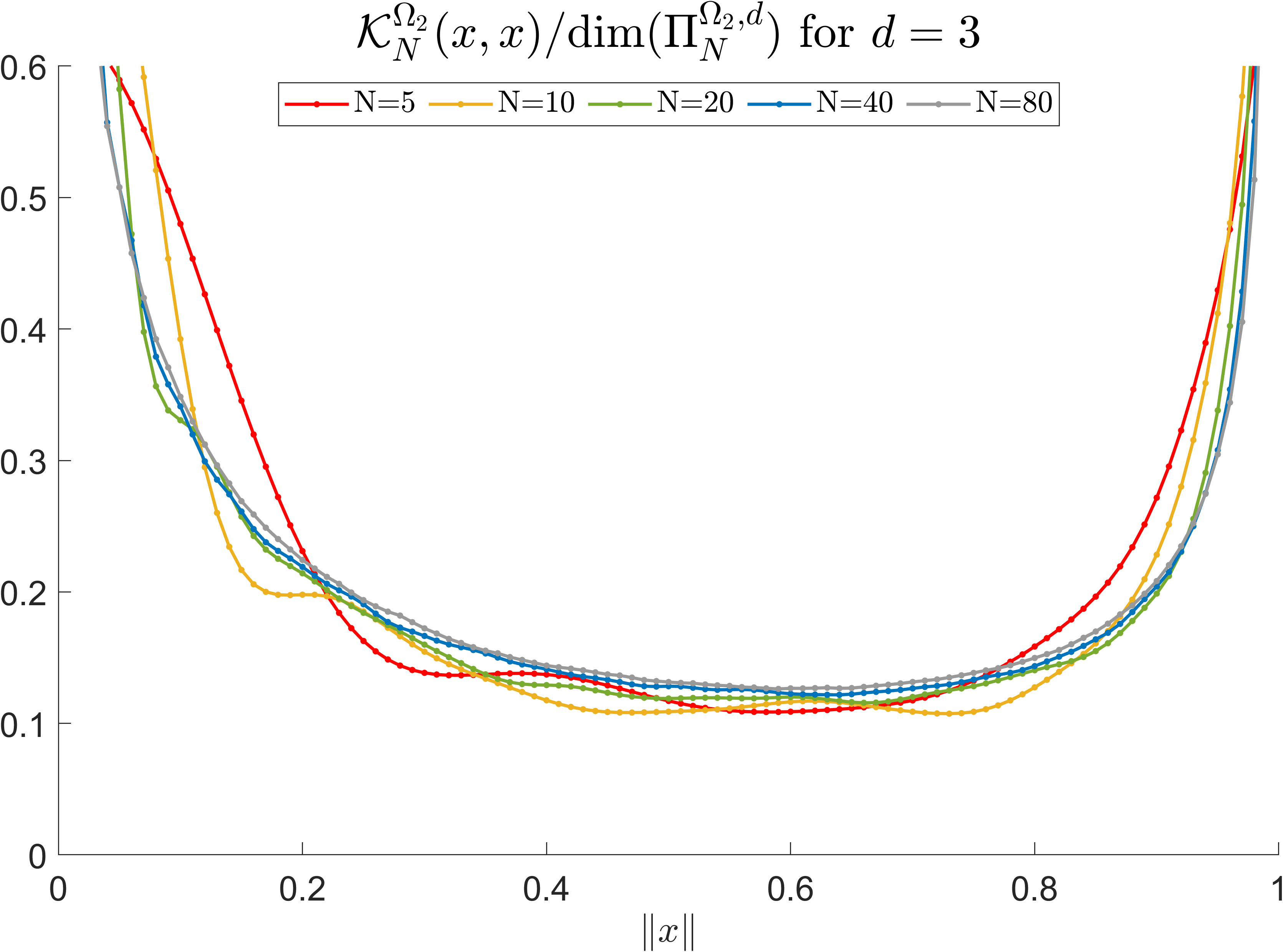}
	\caption{Illustration of $\frac{\mathcal{K}_{\bw,\Omega}(x,x)}{\mathrm{dim}(\Pi_{\bw,\Omega}^{d})}$ for two different spectral shapes $\Omega$. Left column: $\Omega=\Omega_1$, right column: $\Omega=\Omega_2$. Top row: dimension $d=2$, bottom row: dimension $d=3$.} 
	\label{fig:W_Omega_Zernikes}
\end{figure}

These shapes are illustrated in Figure \ref{fig:omegaillust}. Figure \ref{fig:W_Omega_Zernikes} indicates the quantity ${\Klo(x,x)}/{\ndl}$ for both choices of spectral shapes, for dimensions $d=2,3$, and for various bandwidths $L$. Again, it can be seen that in both cases ${\Klo(x,x)}/{\ndl}$ seems to converge to weight functions $W_\Omega$ as $L$ increases. However, we could not prove this and do not have a closed-form expression of the asymptotic weight functions. Nonetheless, one can clearly see the different behaviour of the left and right graphs in Figure \ref{fig:W_Omega_Zernikes}, especially near the origin and near the boundary (similar to the differences between the rigorously derived weight functions $\widehat{W}$ and $\widetilde{W}$).

\begin{figure}[!t]
	\centering
	\includegraphics[scale=0.27,angle=0]{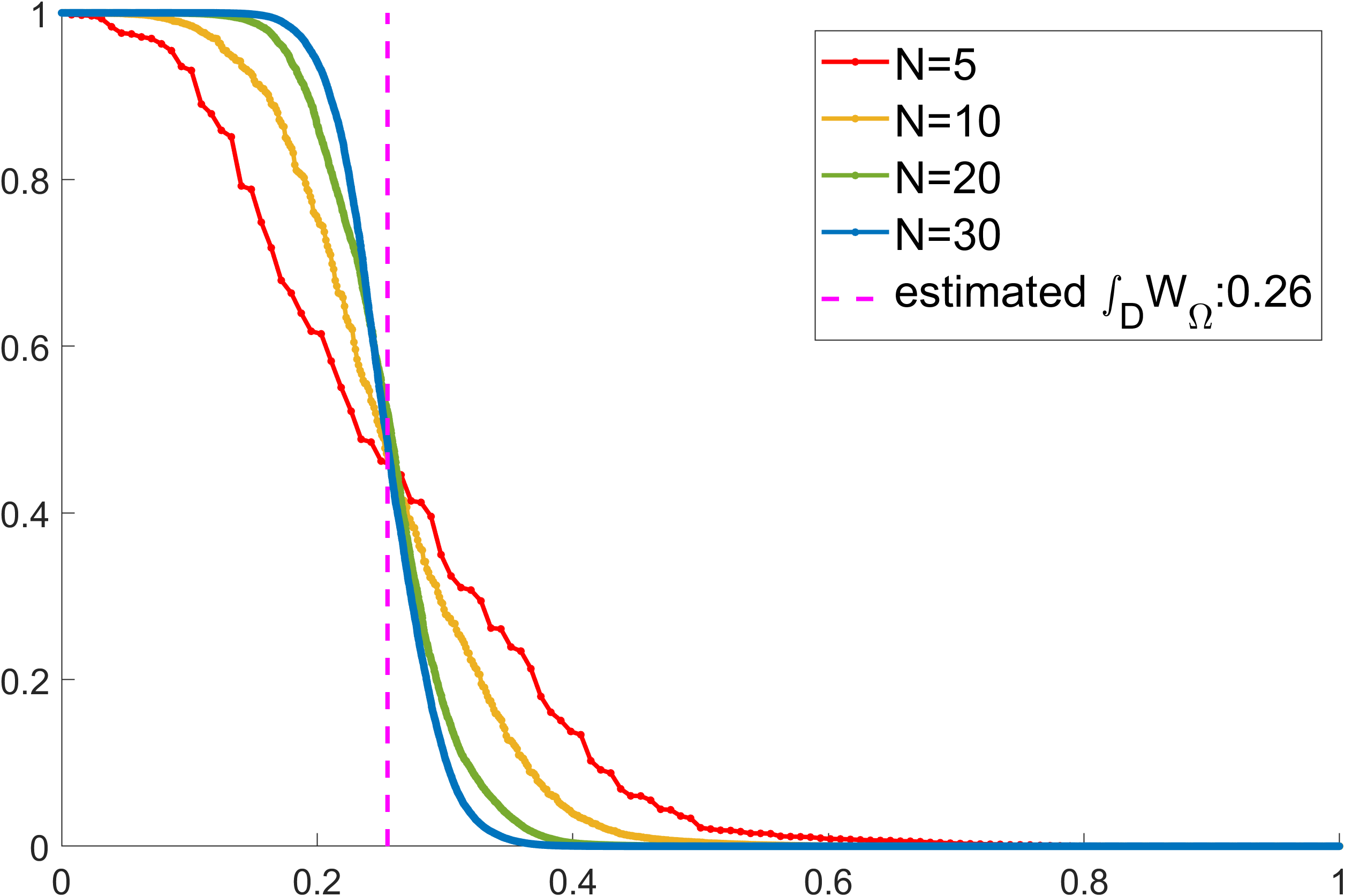}\quad	\quad
	\includegraphics[scale=0.27,angle=0]{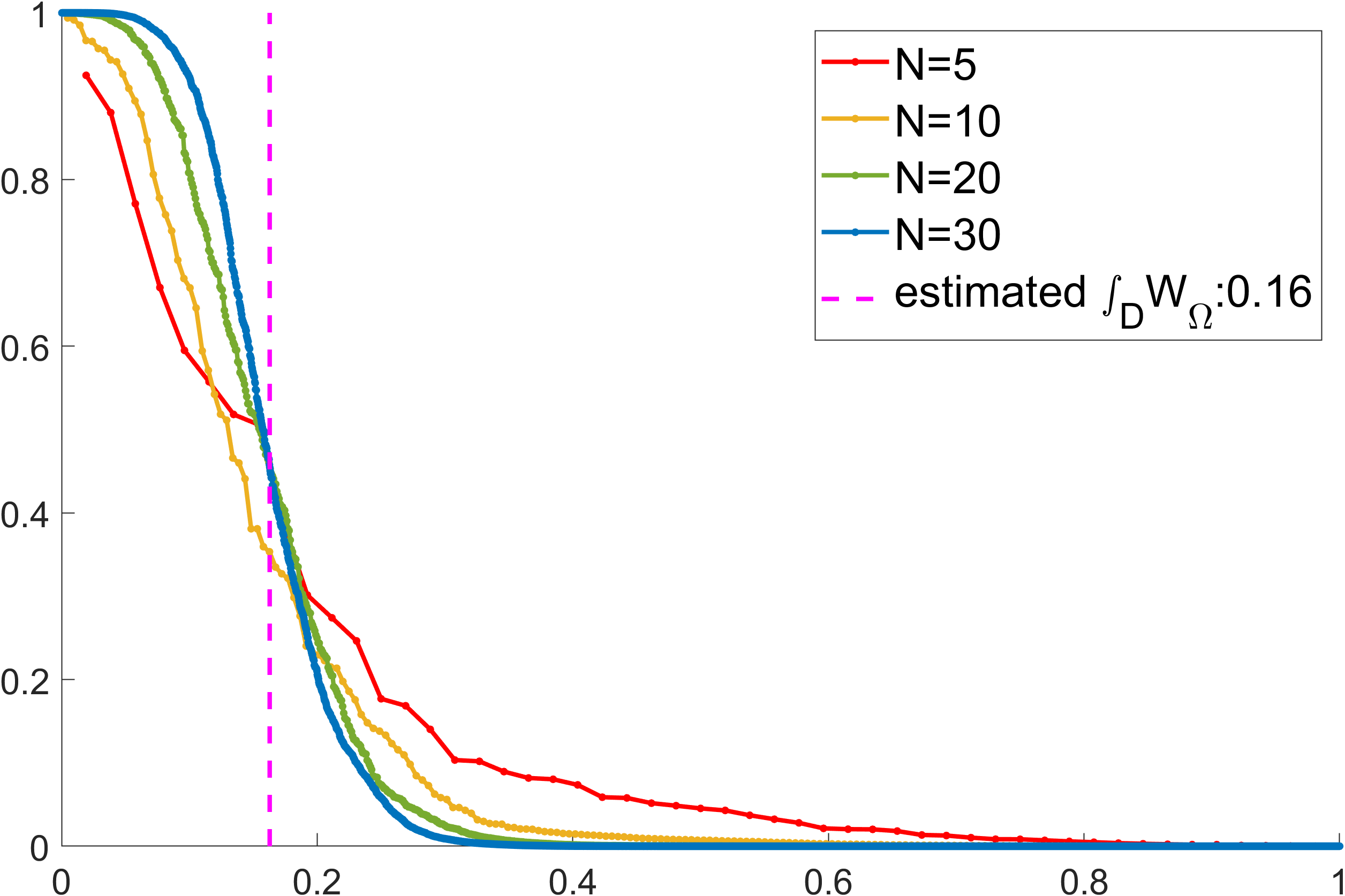}
	\caption{Illustration of the eigenvalue distributions $\lambda_i^{\Omega}(D_1,\bw)$ for $\Omega=\Omega_1$ (left) and $\Omega=\Omega_2$ (right). Dashed vertical lines indicate the empirical relative Shannon numbers $\sum_{i=1}^{\mathcal{N}_{L,\Omega}^d}\lambda_i^{\Omega}(D_1,\bw)/\mathcal{N}_{L,\Omega}^d$.}
	\label{fig:ED_Omega_Zernikes}
\end{figure}

The eigenvalue distribution $\lambda_i^\Omega(D,\bw)$ of the spatiospectral localization operator $\tbtl^\Omega$, for dimension $d=3$, is demonstrated in Figure \ref{fig:ED_Omega_Zernikes}. This time only for the spatial domain $D_1$ from the Section \ref{sec:num} (which is localized towards the origin of the ball). The numerical evaluation has been done analogous to the procedure described in Remark \ref{rem:nujmdesc}. One can see that the empirical relative Shannon numbers are smaller for the notion of bandlimit determined by $\Pi_{\bw,\Omega_1}^{d}$ than for the notion of bandlimit determined by $\Pi_{\bw,\Omega_2}^{d}$. This is in agreement with the expectations that one would get from the graphs in Figure \ref{fig:W_Omega_Zernikes}.

\section{On the influence of the parameter ${\supscriptInCFourOl}_j$ in Fourier-Jacobi bases}\label{app:arbrhoj}

If we choose ${\supscriptInCFourOl}_j$ to be constant with respect to $j$, the particular choice of sequences $(m_k)_{k\in\mathbb{N}},(n_k)_{k\in\mathbb{N}}$ has no influence on $W$. More precisely, the following statements hold true.

\begin{proposition}\label{prop:convergenceofChristoffelfunctionoftildePiconst}
	Let ${\supscriptInCFourOl}_j=\supscriptInCFourOl$ for all $j\in\N_0$, with $\supscriptInCFourOl\in\R$ some constant such that $\supscriptInCFourOl+\tfrac{d-2}{2}>-1$. Additionally, let  $(m_k)_{k\in\mathbb{N}},(n_k)_{k\in\mathbb{N}}\subset\N_0$ be two sequences with $\lim_{k\to\infty}m_k=\lim_{k\to\infty}n_k=\infty$. Then it holds that, for any $x$ in the interior of $\BB^d\setminus\{0\}$ (i.e., for $0<\|x\|<1$), the following pointwise limit holds true:
	\begin{align}
		\lim_{k\to\infty} \frac{\widetilde{\mathcal{K}}_{m_k,n_k}(x,x)}{\widetilde{\mathcal{N}}_{m_k,n_k}} = \widetilde{W}(x),
	\end{align}
	with $\widetilde{W}$ given by \eqref{eqn:wtilde1}. Furthermore, the above limit holds uniformly for any compact subset $D$ in the interior of $\BB^d\setminus\{0\}$. 
\end{proposition}

\begin{proposition}\label{prop:unilimit2const}
	Let ${\supscriptInCFourOl}_j=\supscriptInCFourOl$ for all $j\in\N_0$, with $\supscriptInCFourOl\in\R$ some constant such that $\supscriptInCFourOl+\tfrac{d-2}{2}>-1$. And let  $(m_k)_{k\in\mathbb{N}},(n_k)_{k\in\mathbb{N}}\subset\N_0$ be two sequences with $\lim_{k\to\infty}m_k=\lim_{k\to\infty}n_k=\infty$. Additionally, let $D$ be a compact set contained in the interior of $\BB^d\setminus\{0\}$. For some given $x\in D$, we define $x_{t,\xi}^m=(\|x\|+\frac{t}{m+1})\xi$, with $\xi\in\Sphered$ and $t\in\R$. Then, uniformly for $x\in D$, for $\xi\in\Sphered$, and for $t$ in a compact subset of $\R$, it holds
	\begin{equation}\label{unilimit2const}
		\lim_{k\to\infty}\frac{\widetilde{\mathcal{K}}_{m_k,n_k}(x,x_{t,\xi}^{m_k})}{\widetilde{\mathcal{K}}_{m_k,n_k}(x,x)}=\textnormal{sinc}\left(\frac{2t}{\sqrt{1-\|x\|^2}}\right)\,\,\lim_{k\to\infty}\,\frac{\mathrm{vol}(\Sphered)}{\mathcal{N}_{n_k}(\Sphered)}\mathcal{K}_{\mathrm{Harm}_{n_k}}\left(\frac{x}{\|x\|},\xi\right).
	\end{equation}
\end{proposition}

\begin{theorem}\label{thm:eigendistributionfortildePiconst}
	Let ${\supscriptInCFourOl}_j=\supscriptInCFourOl$ for all $j\in\N_0$, with $\supscriptInCFourOl\in\R$ some constant such that $\supscriptInCFourOl+\tfrac{d-2}{2}>-1$. And let  $(m_k)_{k\in\mathbb{N}},(n_k)_{k\in\mathbb{N}}\subset\N_0$ be two sequences with $\lim_{k\to\infty}m_k=\lim_{k\to\infty}n_k=\infty$. Additionally, let $D\subset\B^d$ be a Lipschitz domain. Then, for any $0<\varepsilon<1/2$ and any $0<\tau<1$, it holds
	\begin{align}\label{eqnC4:eigdis1_tildePiconst}
		\lim_{k\to\infty} \frac{\sharp\{i:\varepsilon<\widetilde{\lambda}_i(D;m_k,n_k)<1-\varepsilon\}}{\widetilde{\mathcal{N}}_{m_k,n_k}}&=0,
		\\\label{eqnC4:eigdis2_tildePiconst}
		\lim_{k\to\infty}\frac{\sharp\{i:\tau<\widetilde{\lambda}_i(D;m_k,n_k)\leq 1\}}{\widetilde{\mathcal{N}}_{m_k,n_k} }&=\int_{D}\widetilde{W}(x)\diffsymbol x,
	\end{align}
	with $\widetilde{W}$ given by \eqref{eqn:wtilde1}.
\end{theorem} 

The proofs of Propositions \ref{prop:convergenceofChristoffelfunctionoftildePiconst}, \ref{prop:unilimit2const}, and Theorem \ref{thm:eigendistributionfortildePiconst} follow completely analogous to the corresponding proofs in Section \ref{sec:seq} once one observes that the choice of a constant ${\supscriptInCFourOl}_j$ fully decouples the summation indices of the radial and spherical contributions of $\widetilde{\mathcal{K}}_{m_k,n_k}$. Thus, we do not further elaborate them here.

\end{document}